\documentclass[12pt]{amsart}
\usepackage[utf8]{inputenc}

\usepackage{amsmath,amsthm,amssymb,mathscinet,bbm}
\usepackage{parskip}
\usepackage{geometry}
\usepackage{xcolor} 
\usepackage{xfrac, nicefrac}
\usepackage{exscale, relsize}

\usepackage{mathtools} 
\usepackage[hidelinks]{hyperref}
\author{Akash Singha Roy}
\address{Department of Mathematics \\ University of Georgia \\ Athens, GA 30602}
\email{akash01s.roy@gmail.com}


\renewcommand\phi\varphi
\usepackage{accents}

\renewcommand{\pod}[1]{\allowbreak\mathchoice
  {\if@display \mkern 18mu\else \mkern 8mu\fi (#1)}
  {\if@display \mkern 18mu\else \mkern 8mu\fi (#1)}
  {\mkern4mu(#1)}
  {\mkern4mu(#1)}
}
\usepackage{graphicx}
\DeclareMathAlphabet{\curly}{U}{rsfs}{m}{n}

\newcommand{\1}{\mathbbm{1}}

\newcommand{\F}{\mathbb{F}}
\newcommand\Z{\mathbb{Z}}
\allowdisplaybreaks

\usepackage{cleveref}
\crefname{section}{§}{§§}
\Crefname{section}{§}{§§}

\newcommand\NatNos{\mathbb N}

\newcommand\Q{\mathbb{Q}}

\newtheorem{thm}{Theorem}[section]
\newtheorem{cor}[thm]{Corollary}
\newtheorem{prop}[thm]{Proposition}
\newtheorem{lem}[thm]{Lemma}

\newtheorem*{thmN}{Theorem N}
\theoremstyle{remark}

\newcommand\ord{\mathrm{ord}}

\newcommand\Ree{\mathrm{Re}}

\newcommand\bbm{\mathbbm 1}
\newcommand\alphatil{\widetilde{\alpha}}

\newcommand\condofchi{\cond(\chi)}
\newcommand\cond{\mathfrak f}

\newtheorem{innercustomgeneric}{\customgenericname}
\providecommand{\customgenericname}{}
\newcommand{\newcustomtheorem}[2]{
  \newenvironment{#1}[1]
  {
   \renewcommand\customgenericname{#2}
   \renewcommand\theinnercustomgeneric{##1}
   \innercustomgeneric
  }
  {\endinnercustomgeneric}
}

\newcustomtheorem{customthm}{Theorem}
\newcustomtheorem{customlemma}{Lemma}

\newcommand\Dmin{D_{\mathrm{min}}}

\newcommand\sm\setminus

\newcommand\VqmPr{\mathcal V_{q, m}'}
\newcommand\VQZmPr{\mathcal V_{Q_0, m}'}

\newcommand\largesum{\sum}
\newcommand\qtil{\widetilde q}

\newcommand\ii{\text{i}}
\newcommand\Wiv{W_{i, v}}
\newcommand\Wik{W_{i, k}}
\newcommand\Wikset{\{\Wik\}_{1 \le i \le K}}

\newcommand\Wivfullset{\{\Wiv\}_{\substack{1 \le i \le K\\1 \le v \le V}}}
\newcommand\Wivfullfam{(\Wiv)_{\substack{1 \le i \le K\\1 \le v \le V}}}
\newcommand\prodik{\prod_{i=1}^K}
\newcommand\sumik{\sum_{i=1}^K}
\newcommand\alphavq{\alpha_v(q)}
\newcommand\Fiset{\{F_i\}_{i=1}^K}
\newcommand\subZT{\subset \Z[T]}
\newcommand\Gjset{\{G_j\}_{j=1}^M}
\newcommand\omegaFprod{\omega(F_1 \cdots F_K)}
\newcommand\Fifam{(F_i)_{i=1}^K}
\newcommand\Filist{F_1, \dots, F_K}
\newcommand\aifam{(a_i)_{i=1}^K}
\newcommand\betaFi{\beta\left(\Filist\right)}
\newcommand\qlelogxKZ{q \le (\log x)^{K_0}}
\newcommand\degWik{\deg \Wik}
\newcommand\degWiv{\deg \Wiv}
\newcommand\finaimodq{f_i(n) \equiv a_i \pmod q}
\newcommand\Wiklist{W_{1, k}, \dots, W_{K, k}}
\newcommand\Wivlist{W_{1, v}, \dots, W_{K, v}}
\newcommand\betaWiklist{\beta(\Wiklist)}

\newcommand\IFHFiBZ{IFH(\Filist; B_0)}
\newcommand\IFHWikBZ{IFH(\Wiklist; B_0)}

\newcommand\quadxtoinfty{\quad\text{as $x\to\infty$}}
\newcommand\forallifinaimodq{(\forall i) ~ \finaimodq}
 
\newcommand\phiQZphiq{\frac{\phi(Q_0)}{\phi(q)}}
\newcommand\phiQZphiqK{\left(\phiQZphiq\right)^K} 
\newcommand\finaimodQZ{f_i(n) \equiv a_i \pmod{Q_0}}
\newcommand\sumnxfnqfiQZ{\largesum_{\substack{n \le x: ~ (f(n), q)=1\\(\forall i) ~ \finaimodQZ}} 1} 
\newcommand\sumnxfnq{\largesum_{\substack{n \le x\\ (f(n), q)=1}} 1}
\newcommand\osumnxfnq{o\Bigg(\sumnxfnq\Bigg)}
\newcommand\sumnxfnqgcd{\largesum_{\substack{n \le x\\ \gcd(f(n), q)=1}} 1}
\newcommand\ophiqKsumnxfq{o\Bigg(\frac1{\phi(q)^K}\sumnxfnq\Bigg)}

\newcommand\CkQZ{\mathcal C_k(Q_0)}
\newcommand\bbmfnqOne{\bbm_{(f(n), q)=1}}
\newcommand\chiituplist{(\chi_1, \dots, \chi_K)}
\newcommand\psiituplist{(\psi_1, \dots, \psi_K)}
\newcommand\chiZerotuplist{(\chi_0, \dots, \chi_0)}

\newcommand\chiZelltuplist{(\chi_{0, \ell}, \dots, \chi_{0, \ell})}
\newcommand\iskfull{\text{ is }k\text{-full}}
\newcommand\iskfree{\text{ is }k\text{-free}}
\newcommand\nopgrystpkOnedivn{p>y \implies p^{k+1} ~ \nmid ~ n}
\newcommand\nopgrYstpkOnedivn{p>Y \implies p^{k+1} ~ \nmid ~ n}
\newcommand\nopgrystpkOnedivm{p>y \implies p^{k+1} ~ \nmid ~ m}
\newcommand\issqfree{\text{ squarefree}}
\newcommand\inconv{\text{ inconvenient}}
\newcommand\conv{\text{ convenient}}
\newcommand\aifimfam{(a_i f_i(m)^{-1})_{i=1}^K}
\newcommand\VJKqaifim{\mathcal V_{J, K}^{(k)} \left(q; \aifimfam\right)}

\newcommand\wifam{(w_i)_{i=1}^K}
\newcommand\wilist{w_1, \dots, w_K}
\newcommand\VNKqwi{\mathcal V_{N, K}^{(k)} \left(q; \wifam\right)}

\newcommand\VNKQZwi{\mathcal V_{N, K}^{(k)} \left(Q_0; \wifam\right)}

\newcommand\alphakratioN{\frac{\alpha_k(q)^N}{\alpha_k(Q_0)^N}}

\newcommand\FellT{\F_\ell[T]}

\newcommand\FiderivFj{F_i' \prod_{\substack{1 \le j \le K\\j \ne i}} F_j}
\newcommand\FiTderivFjT{F_i'(T) \prod_{\substack{1 \le j \le K\\j \ne i}} F_j(T)}
\newcommand\Fiderivset{\left\{\FiderivFj\right\}_{i=1}^K}
\newcommand\C{\mathbb C}
\newcommand\Res{\text{Res}}
\newcommand\Ftil{\widetilde F}
\newcommand\FtilT{\widetilde F(T)}
\newcommand\FellUnits{\F_\ell^\times}
\newcommand{\diag}{\text{diag}}
\newcommand\Fiprod{F_1 \cdots F_K}
\newcommand\kappaell{{\kappa(\ell)}}
\newcommand\betatil{\widetilde\beta}
\newcommand\disc{\text{disc}}

\newcommand\Fij{F_{i, j}}
\newcommand\Gir{G_{i, r}}
\newcommand\GirlistiOK{G_{1, r}, \dots, G_{K, r}}
\newcommand\GirsetiOK{\{\Gir\}_{1 \le i \le K}}
\newcommand\GirsetOKL{\{\Gir\}_{\substack{1 \le i \le K\\1 \le r \le L}}}
\newcommand\GirsetOKLmat{{(\Gir)_{\substack{1 \le i \le K\\1 \le r \le L}}}}
\newcommand\FijsetOKN{\{\Fij\}_{\substack{1 \le i \le K\\1 \le j \le N}}}
\newcommand\FijVectiOK{(\Fij)_{i=1}^K}

\newcommand\alphajtilq{\widetilde\alpha_j(q)} 
\newcommand\alphaStNq{\alpha_N^*(q)}
\newcommand\VNKtilqwi{\widetilde{\mathcal V}_{N, K}\left(q; \wifam\right)}
\newcommand\VNKtilQwi{\widetilde{\mathcal V}_{N, K}\left(Q; \wifam\right)}
\newcommand\IFHGirBZ{IFH(\GirlistiOK; B_0)}
\newcommand\VNKtilQZwi{\widetilde{\mathcal V}_{N, K}\left(Q_0; \wifam\right)}
\newcommand\alphaNStQZ{\alpha_N^*(Q_0)}
\newcommand\alphaNStqRatio{\frac{\alphaStNq}{\alphaNStQZ}}

\newcommand\chiibarwiprod{\overline\chi_1(w_1) \cdots \overline\chi_K(w_K)}
\newcommand\prodjN{\prod_{j=1}^N}
\newcommand\FijlistiOK{F_{1, j}, \dots, F_{K, j}}

\newcommand\expOomegaq{\exp\big(O(\omega(q))\big)}

\newcommand\OloglogqBddPower{O\big((\log_2 (3q))^{O(1)}\big)}
\newcommand\alphajtil{\widetilde \alpha_j}
\newcommand\psieZ{\psi_{e_0}}
\newcommand\GijPrVectiOK{(G_{i, j'})_{i=1}^K}
\newcommand\GijPrlistiOK{G_{1, j'}, \dots,G_{K, j'}}
\newcommand\lcm{\text{lcm}}
\newcommand\VNKtilqtilwi{\widetilde{\mathcal V}_{N, K}\left(\qtil; \wifam\right)} \newcommand\largesumchiiOKmodqtil{\largesum_{\chi_1, \dots, \chi_K \bmod \qtil}} 
\newcommand\SPr{\mathcal S'}
\newcommand\SDPr{\mathcal S''}
\newcommand\lcmcondchilist{\lcm[\cond(\chi_1), \dots, \cond(\chi_K)]}
\newcommand\iskappaPlusOneFull{\text{ is }(\kappa+1)\text{-full}}
\newcommand\critpoly{\mathcal C}
\newcommand\critpts{\mathcal A}
\newcommand\multtheta{\mu_\theta}

\newcommand\VNOnetilqtilw{\widetilde{\mathcal V}_{N, 1}\left(\qtil; w\right)}
\newcommand\VNOnetilQZw{\widetilde{\mathcal V}_{N, 1}\left(Q_0; w\right)}

\newcommand\VTwoOne{\mathcal V_{2, 1}(\ell; w)}
\newcommand\VThreeTwo{\mathcal V_{3, 2}(\ell; u, w)}
\newcommand\Fellbar{\overline{\F}_\ell}
\newcommand\FXFYFZminusu{F(X) F(Y) F(Z) - u}
\newcommand\Fhat{\widehat{F}}
\newcommand\Ghat{\widehat{G}}
\newcommand\FhatXYZ{\widehat{F}(X, Y, Z)}
\newcommand\GhatXYZ{\widehat{G}(X, Y, Z)}
\newcommand\CZeroTilde{\widetilde C_0}
\newcommand\omegaParalleln{\omega_{\parallel}(n)} 
\newcommand\omegaStn{\omega^*(n)}
\newcommand\foralljPjvjmodq{(\forall j)~ P_j \equiv v_j \pmod q}

\newcommand\expOlogThreexSqPluslogTwoqSq{\exp\left(O\big((\log_3 x)^2 + (\log_2 (3q))^{O(1)}\big)\right)}
\newcommand\omegakn{\omega_k(n)}
\newcommand\Sigmatilrs{\widetilde\Sigma_{r, s}}
\newcommand\NtiltauOnetaus{\widetilde{\mathcal N}_{r, s}(\tau_1, \dots, \tau_s)}
\newcommand\xonekphiqKlogxOMalpha{\frac{x^{1/k}}{\phi(q)^K (\log x)^{1-2\alpha_k/3}}}
\newcommand\xonekphiqlogxOMalpha{\frac{x^{1/k}}{\phi(q) (\log x)^{1-2\alpha_k/3}}}

\newcommand\xoneklogxOMalphaHalf{\frac{x^{1/k}}{(\log x)^{1-\alpha_k/2}}}
\newcommand\Sigmars{\Sigma_{r, s}}
\newcommand\mcTOne{\mathcal T_1}
\newcommand\mcTwo{\mathcal T_2}

\newcommand\Mrstil{\widetilde{\mathcal M}_{r,s}}
\newcommand\MrstilcOnecs{\Mrstil(c_1, \dots, c_s)}
\newcommand\cjfam{(c_j)_{j=1}^s}
\newcommand\Vsrkqcjw{\mathcal V_{r+s, K}(q; ~\cjfam; ~\wifam)}
\newcommand\VrsKPrqck{V_{r+s, K}'(q; ~\cjfam)}
\newcommand\mcI{\mathcal I}
\newcommand\VRKkqai{\mathcal V_{R, K}^{(k)}\big(q; \aifam\big)}
\newcommand\VRKMinusOnekqai{\mathcal V_{R, K-1}^{(k)}\big(q; (a_i)_{i=1}^{K-1}\big)}

\newcommand\largesummqSmfmqOne{\largesum_{\substack{m \le x: ~ P(m) \le q\\\gcd(f(m), q)=1}} \frac1{m^{1/k}}}
\newcommand\xkonekloglogxlogx{\frac{x^{1/k} (\log_2 x)^{O(1)}}{\log x}}
\newcommand\xonekqKloglogxlogx{\frac{x^{1/k} (\log_2 x)^{O(1)}}{q^K\log x}}
\newcommand\expOsqrtlogq{\exp\big(O(\sqrt{\log q})\big)}
\newcommand\Fchis{F_\chi(s)}
\newcommand\FchiTils{\widetilde{F}_\chi(s)}
\newcommand\FchiTil{\widetilde{F}_\chi}
\newcommand\Hchis{H_\chi(s)}
\newcommand\HchiTils{\widetilde{H}_\chi(s)}
\newcommand\HchiTil{\widetilde{H}_\chi}
\newcommand\GchiOnes{G_{\chi, 1}(s)}
\newcommand\GchiOne{G_{\chi, 1}}
\newcommand\GchiTwos{G_{\chi, 2}(s)}
\newcommand\GchiTwo{G_{\chi, 2}}
\newcommand\prodchiifin{\prod_{i=1}^K \chi_i(f_i(n))}
\newcommand\bbmfnQOne{\bbm_{(f(n), Q)=1}}
\newcommand\ElQ{\mathcal L_Q}
\newcommand\ElQt{\mathcal L_Q(t)}

\newcommand\FOnesk{F_1(sk)}
\newcommand\gsk{g(sk)}
\newcommand\chihat{\widehat\chi}
\newcommand\cchihat{c_{\chihat}}
\newcommand\FOneskcchihat{F_1(sk)^{\cchihat}}
\newcommand\gskcchihat{g(sk)^{\cchihat}}
\newcommand\OnekOneMinuscOneLQtk{\frac1k \left(1-\frac{c_1}{\ElQt}\right)}
\newcommand\prodQOnemidQ{\prod_{Q_1 \mid Q}}
\newcommand\prodPrimChar{\prod_{\substack{\psi \bmod Q_1\\ \psi\text{ primitive}}}}
\newcommand\Lskpsi{L(sk, \psi)}
\newcommand\LskpsiGammapsi{\Lskpsi^{\gamma(\psi)}}
\newcommand\OneminPsiellksGamPsi{\left(1-\frac{\psi(\ell)}{\ell^{ks}}\right)^{\gamma(\psi)}}
\newcommand\OneminOnepks{1-\frac1{p^{ks}}}
\newcommand\OneminOnepksMincchi{\left(\OneminOnepks\right)^{-\cchihat}}
\newcommand\chibar{\overline\chi}
\newcommand\OneByk{\frac1k}
\newcommand\betaeByk{\frac{\beta_e}k}
\newcommand\LDgsk{\frac{g'(sk)}{\gsk}}
\newcommand\AbsLDgsk{\left|\LDgsk\right|}
\newcommand\LDGchiOnes{\frac{G_{\chi, 1}'(s)}{\GchiOnes}}
\newcommand\AbsLDGchiOnes{\left|\LDGchiOnes\right|}
\newcommand\GchiTwosDer{G_{\chi, 2}'(s)}
\newcommand\sumbInUQWkbInUQ{\sum_{\substack{b \in U_Q\\W_k(b) \in U_Q}}}
\newcommand\prodbInUQWkbInUQ{\prod_{\substack{b \in U_Q\\W_k(b) \in U_Q}}}
\newcommand\sumpbmodQ{\sum_{p \equiv b \pmod Q}}
\newcommand\prodpbmodQ{\prod_{p \equiv b \pmod Q}}
\newcommand\prodpMidQWkpInUQ{\prod_{\substack{p \mid Q\\W_k(p) \in U_Q}}}
\newcommand\sumpMidQWkpInUQ{\sum_{\substack{p \mid Q\\W_k(p) \in U_Q}}}
\newcommand\bbmfpvQOne{\bbm_{f(p^v), Q)=1}}
\newcommand\bbmWkpQOne{\bbm_{W_k(p), Q)=1}}
\newcommand\pvs{p^{vs}} 
\newcommand\pks{p^{ks}}
\newcommand\pvks{p^{vks}}
\newcommand\prodchiifipv{\prod_{i=1}^K \chi_i(f_i(p^v))}
\newcommand\OneminOneWkpQOnepks{1-\frac{\bbmWkpQOne}{\pks}}
\newcommand\OneminOneWkpQOnepkscchihat{\left(\OneminOneWkpQOnepks\right)^{\cchihat}}

\newcommand\sumPsimodQ{\sum_{\psi \bmod Q}}
\newcommand\psibar{\overline\psi} 
\newcommand\psibarb{\psibar(b)}
\newcommand\sMinFracOnek{s-\frac1k}
\newcommand\sMinOnekAlphcchi{\left(\sMinFracOnek\right)^{\alpha_k(Q) \cchihat}} 
\newcommand\sMinFracBetaek{s-\frac{\beta_e}k}
\newcommand\sMinBetaekAlphcchiGam{\left(\sMinFracBetaek\right)^{-\alpha_k(Q) \cchihat \gamma(\psi_e)}}
\newcommand\HchiOnek{H_\chi(1/k)}
\newcommand\HchiFracOnek{H_\chi\left(\frac1k\right)}
\newcommand\AbsHchiOnek{|\HchiOnek|} 
\newcommand\AbsHchiFracOnek{\left|H_\chi\left(\frac1k\right)\right|}
\newcommand\AbsHchiTilFracOnek{\left|\widetilde H_\chi\left(\frac1k\right)\right|}
\newcommand\logxalphakepsFiv{(\log x)^{\alpha_k(Q)\epsilon_1/5}}
\newcommand\logxalphakepsFour{(\log x)^{\alpha_k(Q)\epsilon_1/4}}
\newcommand\ScaledOneMincOneTwoElqt{\frac1k\left(1-\frac{c_1}{2 \ElQt}\right)}
\newcommand\ScaledOneMincOneFourElqt{\frac1k\left(1-\frac{c_1}{4 \ElQt}\right)}
\newcommand\TextIm{\mathrm{Im}}
\newcommand\BetaFactor{(1-\beta_e)^{-\alpha_k(Q)}}
\newcommand\BetaFactorTwice{(1-\beta_e)^{-2\alpha_k(Q)}}
\newcommand\ds{\mathrm{d}s}
\newcommand\ScaledOnekOnelogX{\frac1k\left(1+\frac1{\log X}\right)}
\newcommand\GammajBar{\overline{\Gamma_j}}
\newcommand\MinGammajBar{-\GammajBar}
\newcommand\MinGammaFourStBar{-\overline{\Gamma_4^*}}

\newcommand\MinGammaThrBar{-\overline{\Gamma_3}}
\newcommand\MinGammaTwoBar{-\overline{\Gamma_2}}
\newcommand\WUDkAdmSet{\mathcal Q(k; f_1, \cdots, f_K)}
\newcommand\sumnolimitsSt{\sum\nolimits^*}

\numberwithin{equation}{section}
\begin{document} 
\title[Equidistribution of families of multiplicative functions]{Joint distribution in residue classes of families of polynomially-defined multiplicative functions} 
\begin{abstract} 
We study the distribution of families of multiplicative functions among the coprime residue classes to moduli varying uniformly in a wide range, obtaining analogues of the Siegel--Walfisz Theorem for large classes of multiplicative functions. We extend a criterion of Narkiewicz for families of multiplicative functions that can be controlled by values of polynomials at the first few prime powers, and establish results that are completely uniform in the modulus as well as optimal in most parameters and hypotheses. This also significantly generalizes and improves upon previous work done for a single such function in specialized settings. Our results have applications for most interesting multiplicative functions, such as the Euler  totient function $\phi(n)$, the sum-of-divisors function $\sigma(n)$, the coefficients of the Eisenstein series, etc., and families of these functions. For instance, an application of our results shows that for any fixed $\epsilon>0$, the functions $\phi(n)$ and $\sigma(n)$ are jointly asymptotically equidistributed among the reduced residue classes to moduli $q$ coprime to $6$ varying uniformly up to $(\log x)^{(1-\epsilon)\alpha(q)}$, where $\alpha(q) \coloneqq \prod_{\ell \mid q} (\ell-3)/(\ell-1)$; furthermore, the coprimality restriction is necessary and the range of $q$ is essentially optimal. One of the primary themes behind our arguments is the quantitative detection of a certain mixing (or ergodicity) phenomenon in multiplicative groups via methods belonging to the `anatomy of integers', but we also rely heavily on more pure analytic arguments (such as a suitable modification of the Landau--Selberg--Delange method), -- whilst using several tools from arithmetic and algebraic geometry, and from linear algebra over rings as well.
\end{abstract}
\maketitle 
\section{Introduction}
We say that an integer-valued arithmetic function $g$ is \textsf{uniformly distributed} (or \textsf{equidistributed}) modulo $q$ if $\#\{n \le x: g(n)\equiv b\pmod q\} \sim x/q$ \text{as $x\to\infty$}, for each residue class $b$ mod $q$. 
This definition generalizes naturally to families of arithmetic functions, and has been well-studied for (integral-valued) additive functions, -- with work of Delange \cite{delange69}, \cite{delange74} characterizing when a family of such functions is equidistributed to a fixed modulus $q$. These results have also been partially extended in \cite{PSR}, \cite{akande23} and \cite{SR23}, where the modulus $q$ itself has been allowed to vary up to a certain threshold depending on the stopping point $x$ of inputs. 

However, for multiplicative functions, there are indications that uniform distribution is not the correct notion to consider. For instance, it can be shown that the Euler totient function $\phi(n)$ is almost always divisible by any fixed integer $q>1$, and hence is not equidistributed to any given modulus.  
Motivated by this, Narkiewicz in \cite{narkiewicz66} introduces the notion of weak uniform distribution: Given an integer-valued arithmetic function $f$ and a positive integer $q$, we say that $f$ is \textsf{weakly uniformly distributed} (or \textsf{weakly equidistributed} or \textsf{WUD}) modulo $q$ if there are infinitely many positive integers $n$ for which $\gcd(f(n), q)=1$, and if 
\begin{equation*}
\#\{n \le x: f(n)\equiv a\pmod q\} \sim \frac1{\phi(q)}\#\{n \le x: \gcd(f(n), q)=1\}, \quadxtoinfty,
\end{equation*}
for each coprime residue class $a\bmod q$. This definition extends naturally to families of arithmetic functions: we say that 
the integer-valued arithmetic functions $f_1, \dots, f_K$ are \textsf{jointly weakly equidistributed} (or \textsf{jointly WUD}) modulo $q$ if there are infinitely many $n$ for which $\gcd(f_1(n)$ $ \cdots f_K(n), q)=1$, and if for all coprime residue classes $a_1, \dots, a_K \bmod q$, we have
\begin{multline}\label{eq:jtWUDdef}
\#\{n \le x: \forall i \in [K], ~ ~ f_i(n)\equiv a_i\pmod q\} \sim \frac1{\phi(q)^K}\#\{n \le x: \gcd(f_1(n) \cdots f_K(n), q)=1\}
\end{multline}
as $x \to \infty$. (Here and below, $[K]$ denotes the set $\{1, \dots, K\}$.)

The phenomenon of weak uniform distribution has captured significant interest for specific as well as for general classes of multiplicative functions. Narkiewicz \cite{narkiewicz66} shows that $\phi(n)$ is weakly equidistributed precisely to those moduli $q$ that are coprime to $6$, while \'Sliwa \cite{sliwa73} shows that the sum of divisors function $\sigma(n) = \sum_{d \mid n} d$ is weakly equidistributed mod $q$ exactly when $q$ is not a multiple of $6$. Generalizations of \'Sliwa's result to Fourier coefficients of Eisenstein series, namely the functions $\sigma_r(n) \coloneqq \sum_{d \mid n} d^r$, has been studied in great depth by Narkiewicz, Rayner, Dobrowolski, Fomenko and others; see \cite{fomenko80}, \cite{narkray82}, \cite{narkiewicz83}, \cite[Theorem 6.12]{narkiewicz84}, \cite{rayner88}, \cite{rayner88a}. In fact in \cite[Theorem 1]{narkiewicz66}, Narkiewicz gives a general criterion for deciding weak equidistribution for a single  ``polynomially-defined" multiplicative function $f$, one that can be controlled by the values of polynomials at the first few powers of all primes. While the exact statement requires some set-up, the general idea is that such a function $f$ is weakly equidistributed modulo a fixed positive integer $q$ precisely when for every nontrivial Dirichlet character mod $q$ that acts trivially on a special subgroup of the unit group mod $q$, a certain ``local factor" (or Euler factor) associated to this Dirichlet character vanishes. Narkiewicz dedicates a significant portion of his monograph \cite{narkiewicz84} to give more explicit sufficient conditions that guarantee weak uniform distribution, and to obtain algorithms characterizing all the moduli to which a given ``polynomially-defined" multiplicative function is weakly equidistributed. 

In all these results, the modulus $q$ has been assumed to be fixed. A natural and interesting question is whether weak equidistribution continues to hold as $q$ varies uniformly in a suitable range depending on the stopping point $x$ of inputs, for instance, whether it is possible to obtain analogues of the Siegel--Walfisz Theorem for primes in arithmetic progressions, but with primes replaced by values of multiplicative functions. To this end, given a constant $K_0>0$, we shall say that 
integer-valued arithmetic functions $f_1, \dots, f_K$ are \textsf{jointly weakly equidistributed} (or \textsf{jointly WUD}) mod $q$, \textsf{uniformly for $q \le (\log x)^{K_0}$,} if:
\begin{enumerate}
    \item[(i)] For every such $q$, $\prodik f_i(n)$ is coprime to $q$ for infinitely many $n$, and 
    \item[(ii)] The relation \eqref{eq:jtWUDdef} holds as $x \rightarrow \infty$, uniformly in moduli $\qlelogxKZ$ and in coprime residue classes $a_1, \dots, a_K$ mod $q$. Explicitly, this means that for any $\epsilon>0$, there exists $X(\epsilon)>0$ such that the ratio of the left hand side of \eqref{eq:jtWUDdef} to the right hand side lies in $(1-\epsilon, 1+\epsilon)$ for all $x>X(\epsilon)$, $q \le (\log x)^{K_0}$ and coprime residues $a_1, \dots, a_K$ mod $q$.
\end{enumerate}
If $K=1$ and $f_1=f$, we shall simply say that $f$  is \textsf{weakly equidistributed} (or \textsf{WUD}) mod $q$, \textsf{uniformly for $q \le (\log x)^{K_0}$}. 

The question of weak equidistribution to varying moduli seems to have been first studied in \cite{LPR21},  \cite{PSR22} and 
\cite{PSR23}, 
which made some partial progress towards obtaining a uniform analogue of Narkiewicz's aforementioned criterion for a single ``polynomially-defined" multiplicative function. However, the settings in these papers were highly special instances of the setting in Narkiewicz's original criterion in \cite{narkiewicz66}, so much so that they could not be used to obtain satisfactory  uniform analogues of the aforementioned results on $\sigma_r(n)$. 

As a special case of our results in this manuscript, we are able to extend Narkiewicz's criterion in its full generality to obtain results that are completely uniform in the modulus $q$ and have optimal arithmetic restrictions on $q$. Certain special cases of our results also yield uniform extensions of the aforementioned results on $\sigma_r(n)$. For instance, we get all the following uniform analogues of \'Sliwa's result in \cite{sliwa73}: the sum of divisors function $\sigma(n)$ is weakly equidistributed uniformly to all odd moduli $\qlelogxKZ$ as well as to all even $q$ not divisible by $3$ that are either no more than a small power of $\log x$ or are squarefree without too many distinct prime factors. In addition, uniformity is restored to \textit{all} (resp. to \textit{squarefree}) even $\qlelogxKZ$ that are not multiples of $3$, provided we restrict to inputs $n$ having \textit{six} (resp. \textit{four}) large prime factors counted with multiplicity. By examples constructed in \cite{SR23A}, most of these restrictions are optimal. Applications of our main theorems also yield generalizations of these results for the functions $\sigma_r(n)$, thus obtaining complete uniform extensions of the aforementioned results of Narkiewicz, Rayner, Dobrowolski, Fomenko and others (see the discussion following the statement of Theorem \ref{thm:IFHNecessity}). 

All of these results and improvements are only for a single multiplicative function. In \cite{narkiewicz82}, Narkiewicz generalizes his aforementioned criterion to decide joint weak equidistribution for  \textit{families} of ``polynomially defined" multiplicative functions to a fixed modulus $q$; he uses this generalized criterion in \cite{narkiewicz81} to characterize those fixed $q$ to which the Euler totient $\phi(n)$ and sum of divisors $\sigma(n)$ are jointly weakly equidistributed. However, several arguments in the aforementioned papers (\cite{LPR21}, \cite{PSR22}, \cite{PSR23}) investigating uniform analogues of his previous criterion are all strictly constrained to a single multiplicative function and do not generalize to families. Our main results in this manuscript give complete uniform extensions of Narkiewicz's general criterion in \cite{narkiewicz82} for families of multiplicative functions to a single modulus $q$, extensions that are optimal in both the range of uniformity and the arithmetic restrictions on $q$, as well as in various other parameters. 

The qualitative summary of our main results is as follows. Under certain (provably) unavoidable conditions, a given family $f_1, \dots, f_K$ of polynomially-defined multiplicative functions is jointly weakly equidistributed \textit{exactly} to those moduli $q$ that satisfy Narkiewicz's criterion, and are also allowed to vary uniformly up to small powers of $\log x$, where these powers are all essentially optimal as well. In addition, weak equidistribution is restored in the full ``Siegel-Walfisz range" $\qlelogxKZ$ provided we restrict to inputs $n$ having sufficiently many large prime factors counted with multiplicity. This threshold can be reduced and optimized (thus ensuring equidistribution among larger sample spaces of inputs) whenever $q$ is squarefree, or whenever some reasonable additional control is available on the factorization of $n$ or on the behavior of the functions $f_i$ at some higher prime powers. 

The intuitive explanation for such constraints on our inputs $n$ comes from a certain `mixing' phenomenon in the unit group mod $q$, which can be interpreted as a quantitative ergodicity phenomenon for random walks on multiplicative groups. For example, let $q$ be an odd positive integer. From the set of units $u$ mod $q$ for which $u+1$ is also a unit, choose uniformly at random $u_1, u_2, u_3,\dots$, and construct the sequence of partial products $u_1+1, (u_1+1)(u_2+1), (u_1+1)(u_2+1) (u_3+1), \dots$. Then as we go further into the sequence, each unit mod $q$ is roughly equally likely to appear as one of the products $(u_1+1) \cdots (u_J+1)$. This particular example lies at the core of the weak equidistribution of $\sigma(n)$ to odd moduli. The phenomenon for $\sigma(n)$ to even moduli not divisible by $3$ is analogous, except that we work with units $u$ mod $q$ for which $u^2+u+1$ is also a unit mod $q$.

In general, for certain collections of $K$ multivariate polynomials, the values taken by them that are coprime to $q$ are jointly equidistributed among the unit group mod $q$ whenever the number of variables is significantly larger compared to $K$: having a large number of variables amplifies the power savings in certain character sum bounds, thus ensuring that any $K$ congruences (coming from the $K$ polynomials) maximally ``cut down" the ambient space of tuples. This is reminiscent of a very common phenomenon occurring in several applications of the circle method, such as in Waring's problem. It is to have this large number of variables that it becomes necessary to restrict our inputs $n$ to those having sufficiently many large prime factors, so as to restore weak equidistribution in the full ``Siegel-Walfisz range" $\qlelogxKZ$.  

Our arguments for the main results require ideas comprising a variety of themes. One of the central themes is the exploitation of the aforementioned mixing phenomenon in the multiplicative group via methods belonging to the `anatomy of integers'. In addition, we  crucially require several ``pure analytic" ideas, where we work with certain ``pretentious distances", and we also suitably modify the Landau--Selberg--Delange method to give strong estimates on the mean values of various multiplicative functions taking values in the unit disk. (To this end, we need to analyze a product of $L$-functions raised to complex powers.) Linear algebra over rings comes into play throughout the paper, -- most prominently in combination with methods from combinatorial number theory, -- in order to 
count solution tuples of multiple polynomial congruences in several variables. 
Furthermore, we need to understand the rational points of certain affine varieties over finite fields using tools from arithmetic and algebraic geometry. 

\section{The setting and the main results}
We say that an arithmetic function $f$ is \textsf{polynomially-defined} if there exists $V \ge 1$ and polynomials $\{W_v\}_{1 \le v \le V}$ with integer coefficients satisfying $f(p^v) = W_v(p)$ for all primes $p$ and all $v \in [V]$. To set up for Narkiewicz's general criterion in \cite{narkiewicz82}, we consider $K, V \ge 1$ and polynomially-defined multiplicative functions $f_1, \dots,$ $f_K \colon \NatNos \to \Z$, with defining polynomials $\Wivfullset \subZT$ satisfying $f_i(p^v) = \Wiv(p)$ for any prime $p$, and any $i \in [K], v \in [V]$. 
For any $q$ and $v \in [V]$, define $R_v(q) \coloneqq \{u \in U_q: \prod_{i=1}^K \Wiv(u) \in U_q\};$ here $U_q \coloneqq (\Z/q\Z)^\times$ denotes the multiplicative group mod $q$, so that saying ``$r \in U_q$" for an integer $r$ is synonymous with saying that ``$\gcd(r, q)=1$". Fix $k \in [V]$ and assume that $\Wikset$ are all nonconstant. We say that a positive integer $q$ is \textsf{$k$-admissible} (with respect to the family $\Wivfullfam$) if the set $R_k(q)$ is nonempty but the sets $R_v(q)$ are empty for all $v < k$. We define $\WUDkAdmSet$ to be the set of all $k$-admissible integers $q$ such that for every tuple $(\chi_1, \dots, \chi_K) \ne \chiZerotuplist$ of Dirichlet characters\footnote{Here $\chi_0$ or $\chi_{0, q}$ denotes, as usual, the trivial or principal character mod $q$.} mod $q$ for which $\prodik \chi_i(W_{i, k}(u)) = 1$ for all $u \in R_k(q)$, there exists a prime $p$ satisfying 
\begin{equation}\label{eq:JtWUD_Crit}
\largesum_{j \ge 0} ~ \frac{\chi_1(f_1(p^j)) \cdots \chi_K(f_K(p^j))}{p^{j/k}} = 0.
\end{equation}
Narkiewicz's criterion \cite[Theorem 1]{narkiewicz82} in this setting is then stated as follows.
\begin{thmN}\label{thm:NarkGen}
Fix a $k$-admissible integer $q$. The functions $f_1, \dots, f_K$ are jointly weakly equidistributed modulo $q$ if and only if $q \in \WUDkAdmSet$. 
\end{thmN}

As mentioned in the introduction, the first steps towards obtaining uniform analogues of a special case of Theorem N for a single multiplicative function were taken in \cite{LPR21}, \cite{PSR22} and \cite{PSR23}. However, several of the arguments in these papers cannot be generalized to families of multiplicative functions (i.e. the cases $K>1$), and even in the special case $K=1$, those results are far from being complete uniform analogues of Theorem N because they crucially need $q$ to be $1$-admissible and have sufficiently large prime factors, and also crucially need the defining polynomial $W_{1, 1}$ to be separable. 

In this work, we extend Narkiewicz's general criterion Theorem N to obtain results that are completely uniform in the modulus $q$ varying up to a fixed but arbitrary power of $\log x$. Our results will mostly not impose any additional restrictions, beyond those that can be \textit{proven} to be necessary and essentially optimal. These results are thus also new for a single multiplicative function as they address all the aforementioned shortcomings of \cite{LPR21}, \cite{PSR22} and \cite{PSR23}. For concrete and provably unavoidable reasons (see Theorems \ref{thm:MIHNecessity} and \ref{thm:IFHNecessity} below), we shall need to impose two additional hypotheses to get uniform analogues of Theorem N. First, we will need the polynomials $\Wikset$ to be multiplicatively independent. Here, we say that the polynomials $\{F_i\}_{1 \le i \le K} \subZT$ are \textsf{multiplicatively independent (over $\Z$)} if there is no nonzero tuple of integers $(c_i)_{i=1}^K$ for which the product $\prodik F_i^{c_i}$ is identically constant in $\Q(T)$. 

For the second hypothesis, we shall need the following set-up. Given nonconstant polynomials $\Fiset \subZT$, we factor $F_i \eqqcolon r_i \prod_{j=1}^M G_j^{\mu_{ij}}$ with $r_i \in \Z$, $\Gjset \subZT$ being pairwise coprime primitive\footnote{A polynomial in $\Z[T]$ is said to be \textsf{primitive} when the greatest common divisor of its coefficients is $1$.} irreducible polynomials and with $\mu_{ij}$ being nonnegative integers, such that each $G_j$ appears with a positive exponent $\mu_{ij}$ in some $F_i$.  We let $\omegaFprod \coloneqq M$ 
and define the \textsf{exponent matrix} of $\Fifam$ to be the $M \times K$ matrix 
$$E_0 \coloneqq E_0(\Filist) \coloneqq 
\begin{pmatrix}
\mu_{11} & \cdots & \mu_{K1}\\
\cdots & \cdots & \cdots\\
\cdots & \cdots & \cdots\\
\mu_{1M} & \cdots & \mu_{KM}
\end{pmatrix} \in \mathbb M_{M\times K}(\Z),$$
so that $E_0$ has a positive entry in each row. By the theory of modules over a principal ideal domain, $E_0$ has a Smith Normal Form given by the $M \times K$ diagonal matrix $\diag(\beta_1, \dots, \beta_r)$, where $r \coloneqq \min\{M, K\}$ and $\beta_1, \dots, \beta_r$ are integers (possibly zero) satisfying $\beta_j \mid \beta_{j+1}$ for each $1 \le j < r$ (for the moment, we accept the convention that $0 \mid 0$). The $\beta_j$ are often called the \textsf{invariant factors} of the matrix $E_0$.\footnote{In practice, it is usually the nonzero $\beta_j$ that are called the invariant factors but this terminology will be more convenient for us (and the possibility of any $\beta_j$ being zero shall soon become obsolete anyway).} We shall use $\betaFi$ to denote the last invariant factor $\beta_r$. (Here we fixed some ordering of the $G_j$ to define the exponent matrix $E_0(\Filist)$ but the invariant factors are independent of this ordering.) We now state our second hyothesis:

\textsf{Invariant Factor Hypothesis:} Given $B_0>0$, we shall say that a positive integer $q$ \textsf{satisfies $\IFHFiBZ$} if $\gcd(\ell-1, \beta(\Filist))=1$ for any prime $\ell \mid q$ satisfying $\ell>B_0$. 

For instance, if $\prodik F_i$ is separable over $\Q$ (or more generally, if the exponent matrix $E_0(F_1, \dots,$ $F_K)$ is equivalent to the diagonal matrix $\diag(1, \dots, 1)$), then $\beta(\Filist)$ $=1$, so any $q$ satisfies $IFH(F_1, \dots, $ $F_K; B_0)$ for any $B_0>0$. Note that the polynomials $\Fiset \subZT$ are multiplicatively independent if and only if the columns of $E_0(\Filist)$ are $\Q$-linearly independent. In this case, $\omega(\Fiprod) = M \ge K$ and $\beta(F_1, \dots $ $, F_K)= \beta_K \ne 0$ as the computation of the Smith normal form is a base-change over $\Z$.

We now state the main results of this manuscript, uniform analogues of Theorem N. The following set-up will be assumed in the main results below: Fix $K, V \ge 1$ and $K_0, B_0 >0$.
\begin{itemize}
\item Consider multiplicative functions $f_1, \dots,$ $f_K \colon \NatNos \to \Z$ and  polynomials $\Wivfullset \subZT$ satisfying $f_i(p^v) = \Wiv(p)$ for any prime $p$, any $i \in [K]$ and $v \in [V]$.
\item Consider the multiplicative function $f \coloneqq \prodik f_i$ and the polynomials $\{W_v\}_{1\le v \le V} \subZT$ given by $W_v \coloneqq \prodik \Wiv$, so that $f(p^v) = W_v(p)$ for all primes $p$ and all $v \in [V]$.
\item For any $q$ and $v \in [V]$, define 
$R_v(q)$ as before the statement of Theorem N so that $R_v(q) = \{u \in U_q: W_v(u) \in U_q\}.$ Let $\alphavq \coloneqq \frac1{\phi(q)} \#R_v(q)$. Also fix $k \in [V]$ and define \textsf{$k$-admissibility} and the set $\WUDkAdmSet$ as before Theorem N.
\item For each $v \in [V]$, let $D_v \coloneqq \deg W_v = \sumik \degWiv$, $D \coloneqq D_k = \sumik \degWik$, and $\Dmin \coloneqq \min_{1 \le i \le K} ~ \deg \Wik$.
Note that if $q$ is $k$-admissible, then $\alpha_v(q)=0$ for $1 \le v < k$, while $\alpha_k(q) \gg_{W_k} (\log \log (3q))^{-D}$ by the Chinese Remainder Theorem and a standard argument using Mertens' Theorem.
\item Assume that $\Wikset$ are multiplicatively independent. 
\end{itemize}
In Theorems \ref{thm:UnrestrictedInput_SqfrPoly} to \ref{thm:RestrictedInputSqfreeMod1} below, our implied constants depend only on $B_0$ and on the polynomials $\{\Wiv\}_{\substack{1 \le i \le K\\1 \le v \le k}}$, and are in particular independent of $V$ and of the polynomials $\{\Wiv\}_{\substack{1 \le i \le K\\k< v \le V}}$. 
\begin{thm}\label{thm:UnrestrictedInput_SqfrPoly} 
Fix $\epsilon \in (0, 1)$. The functions $f_1, \dots, f_K$ are jointly weakly equidistributed, uniformly to all moduli $\qlelogxKZ$ 
lying in $\WUDkAdmSet$ and satisfying $IFH(\Wiklist;$ $B_0)$, provided any one of the following holds.
\begin{enumerate}
    \item[(i)] \underline{Either} $K=1$ and $W_{1, k} = W_k$ is linear, \underline{or} $K \ge 2$, $q \le (\log x)^{(1-\epsilon) \alpha_k(q)/(K-1)}$ and at least one of $\Wikset$ is linear (i.e., $\Dmin=1$).
    \item[(ii)] $q$ is squarefree and $q^{K-1} \Dmin^{\omega(q)} \le (\log x)^{(1-\epsilon) \alpha_k(q)}$.
    \item[(iii)] $\Dmin>1$ and $q \le (\log x)^{(1-\epsilon) \alpha_k(q) (K-1/\Dmin)^{-1}}$. 
\end{enumerate}
\end{thm}  
By \cite[Theorem 1]{narkiewicz81}, the Euler totient $\phi(n)$ and the sum of divisors $\sigma(n)$ are jointly WUD modulo a fixed integer $q$ precisely when $q$ is coprime to $6$; in fact, $\mathcal Q(1; \phi, \sigma) = \{q: (q, 6)=1\}$. Theorem \ref{thm:UnrestrictedInput_SqfrPoly} shows that this joint weak equidistribution holds uniformly in $q \le (\log x)^{(1-\epsilon)\alpha(q)}$ coprime to $6$, where $\alpha(q) \coloneqq \prod_{\ell \mid q} (\ell-3)/(\ell-1)$ and $\epsilon>0$ is fixed but arbitrary. In subsection \cref{subsec:UnrestrictedInput_Optimality}, we will show that the ranges of $q$ in (i)--(iii) above are all essentially optimal, and that for $K \ge 2$, the range of $q$ under condition (i) is essentially optimal, even if $q$ is squarefree and $\Wikset$ are \textit{all} linear, for \textit{any} choice of (pairwise coprime) linear functions! In particular, this means that the range $(\log x)^{(1-\epsilon)\alpha(q)}$ is essentially optimal for the joint weak equidistribution of $\phi$ and $\sigma$, even if we restrict to squarefree $q$. 

Our constructions in \cref{subsec:UnrestrictedInput_Optimality} will reveal that obstructions to uniformity in $q$ come from inputs $n$ that are $k$-th powers of a prime $P$. Modifying these constructions, we can produce obstructions of the form $mP^k$ with $m$ fixed or growing slowly with $x$. It turns out that the problematic inputs are those with too few large prime factors. More precisely, complete uniformity in $q$ up to a fixed but arbitrary power of $\log x$ can be restored by restricting the set of inputs $n$ to those divisible by a sufficient number (say $R$) of primes exceeding $q$ (here and below, all prime factors are counted with multiplicity unless stated otherwise). A smaller value of $R$ suffices provided we assume that sufficiently many of these primes appear to a $k$-th power in $n$.   

To make these precise, we let $P(n)$ denote the largest prime divisor of $n$, with the convention that $P(1) \coloneqq 1$. Set $P_1(n) \coloneqq P(n)$, and inductively define $P_k(n) \coloneqq P_{k-1}(n/P(n))$. Thus, $P_k(n)$ is the $k$-th largest prime factor of $n$ (counted with multiplicity), with $P_k(n)=1$ if $\Omega(n) < k$. We also use $n_k$ to denote the largest positive integer such that $n_k^k$ is a unitary divisor of $n$; in other words, no prime divisor of the integer $n/n_k^k$ appears to an exponent divisible by $k$. (Informally, we may call $n_k$ the ``reduced $k$-th power part"of $n$; if $k=1$, then $n_1 = n$.) Since $D = 1$ forces $K=1$ and $W_k = W_{1, k}$ to be linear (a case in which Theorem \ref{thm:UnrestrictedInput_SqfrPoly}(i) already gives complete uniformity in $\qlelogxKZ$), we assume in Theorems \ref{thm:RestrictedInputGenMod1} to  \ref{thm:HigherPolyControl1} below that $D \ge 2$.  
\begin{thm}\label{thm:RestrictedInputGenMod1}
The following formulae hold as $x\to\infty$, uniformly in coprime residues $a_1, \dots, a_K$ to moduli $\qlelogxKZ$ lying in $\WUDkAdmSet$ and satisfying $\IFHWikBZ$.
\begin{enumerate}
\item[\textbf{(a)}] 
\begin{multline}\label{eq:WUDPRn}
\#\{n \le x: P_R(n)>q, ~ ~\forallifinaimodq\}\\ \sim \frac1{\phi(q)^K} \#\{n \le x: \gcd(f(n), q)=1\}
\sim \frac1{\phi(q)^K} \#\{n \le x: P_R(n)>q, \gcd(f(n), q)=1\} ,
\end{multline}
where 
\begin{equation*}
\begin{cases}
R = k(KD+1), & \text{ if } k<D\\
R \text{ is the least integer exceeding }k\left(1+(k+1)\left(K-1/D\right)\right), & \text{ if } k \ge D.
\end{cases}
\end{equation*}
\item[\textbf{(b)}]
\begin{multline}\label{eq:WUDPRnk}
\#\{n \le x: P_{KD+1}(n_k)>q, ~ ~\forallifinaimodq\}\\
\sim \frac1{\phi(q)^K} \#\{n \le x: \gcd(f(n), q)=1\} \sim \frac1{\phi(q)^K} \#\{n \le x: P_{KD+1}(n_k)>q, \gcd(f(n), q)=1\}.
\end{multline}
\end{enumerate}
\end{thm}
We remark here that since the inputs $n$ we work with satisfy $\gcd(f(n), q)=1$, the $k$-admissibility of $q$ guarantees that $n$ must differ by a bounded factor from a $k$-full integer (see Lemma \ref{lem:kfreepartbdd} below). This is what makes the anatomy of the reduced $k$-th power parts (i.e. the $n_k$), and hence also the kind of restriction in subpart (b), natural to consider. 
The two formulae \eqref{eq:WUDPRn} and \eqref{eq:WUDPRnk} coincide for $k=1$, and even in the special case $k=K=1$, either of them improves over Theorem 1.4(a) in \cite{PSR23}. The value of $R$ in Theorem \ref{thm:RestrictedInputGenMod1}(a) is optimal for the sum of divisors function $\sigma(n)$ to even moduli $q$; see the discussion on applications following the statement of Theorem \ref{thm:IFHNecessity}.  

For squarefree moduli $q$, it suffices to have much weaker restrictions (that are also exactly or nearly optimal) on the set of inputs $n$ so as to detect weak equidistribution. 
\begin{thm}\label{thm:RestrictedInputSqfreeMod1} 
The following hold as $x\to\infty$, uniformly in coprime residues $a_1, \dots, a_K$ modulo \underline{squarefree} $\qlelogxKZ$ 
lying in $\WUDkAdmSet$ and satisfying $\IFHWikBZ$.
\begin{enumerate} 
\item[\textbf{(a)}] The formulae \eqref{eq:WUDPRn} for $k \ge 2$, with\footnote{Here we write a polynomial $F \in \Z[T]$ as $F = r\prod_{j=1}^M H_j^{\nu_j}$ for some $\nu_j \in \NatNos$ and pairwise coprime primitive irreducibles $H_j \in \Z[T]$, and we say that $F$ is ``squarefull" in $\Z[T]$ if $(\prod_{j=1}^M H_j)^2 \mid F$. This condition is equivalent to saying that $\prod_{\substack{\theta \in \C\\F(\theta) = 0}} (T-\theta)^2 \mid F(T)$ in $\C[T]$, i.e., that every root of $F$ in $\C$ has multiplicity at least $2$.}
\begin{equation*}
R \coloneqq 
\begin{cases}
k(Kk + K - k)+1, & \text{ if at least one of }\Wikset \subZT\text{ is not squarefull.}\\
k(Kk + K - k+1)+1, & \text{ in general.}
\end{cases}
\end{equation*}
\item[\textbf{(b)}] The formulae \eqref{eq:WUDPRnk}, \underline{either} with ``$KD+1$" replaced by $2K+1$ for any $K \ge 1$, \underline{or} with ``$KD+1$" replaced by $2$ for $K = 1$ when $W_k = W_{1, k} \in \Z[T]$ is not squarefull.
\end{enumerate}
\end{thm}
Since $n_1=n$, the case $k=1$ missing in (a) is accounted for in (b).  
It is worthwhile to strive for the optimality of $R$ since doing so ensures weak equidistribution among the largest possible set of inputs $n$. In subsection \cref{subsec:Optimality_RestrictedSqfree}, we show that the restriction on the inputs $n$ in (a) is optimal in the sense that in order to have uniformity in $\qlelogxKZ$, it is not possible to reduce ``$k(Kk + K - k)+1$" to ``$k(Kk + K - k)$". Likewise, the restriction in (b) is nearly optimal in that it is not possible to reduce ``$2K+1$" to ``$2K-1$" for any $K \ge 2$, nor is it possible to reduce the ``$2$" to ``$1$" for $K=1$. (In fact, in all these examples, $\{\Wik\}_{i=1}^K$ will be pairwise coprime irreducibles, making $\prodik \Wik$ separable over $\Q$.) The restriction $k \ge 2$ and the nonsquarefullness condition in (a) are for technical reasons that will become clear from the arguments. 

Our constructions demonstrating the aforementioned optimality or near-optimality of the values of $R$ in Theorem \ref{thm:RestrictedInputSqfreeMod1} will come from multiplicative functions $f_i$ for which the polynomials $\Wikset$ are nonconstant (in fact multiplicatively independent), but for which the polynomials $\{W_{i, k+1}\}_{1 \le i \le K}$ or $\{W_{i, 2k}\}_{1 \le i \le K}$ are constant. In practice however, the $\Wiv$ are often nonconstant for many more values of $v$ (beyond a fixed threshold $k$); in fact, for many well-known arithmetic functions $f$ (such as the Euler totient and sums of divisor-powers $\sigma_r(n) \coloneqq \sum_{d \mid n} d^r$), the values $f(p^v)$ are controlled by nonconstant polynomials $W_v \in \Z[T]$ for \textit{all} $v \ge 1$. Hence, it is natural to ask whether the restriction on inputs $n$ in Theorems \ref{thm:RestrictedInputGenMod1} and \ref{thm:RestrictedInputSqfreeMod1} can be weakened when such additional control on the $f_i$ is available, or in other words, if $V$ (the number of powers of primes at which we are assuming the $f_i$ to be controlled by nonconstant polynomials $\Wiv$) can be taken to be sufficiently large. It turns out that we can almost always do this for squarefree $q$ and in several cases in general. 
Unlike the results stated so far, the implied constants in Theorem \ref{thm:HigherPolyControl1} below could depend on the full set of polynomials $\Wivfullset$. 
\begin{thm}\label{thm:HigherPolyControl1}
Assume that the polynomials $\{\Wiv\}_{1 \le i \le K} \subZT$ are multiplicatively independent for each $v$ satisfying $k \le v \le V$. Let $D_0 \coloneqq \max_{k \le v \le V} D_v = \max_{k \le v \le V} \sum_{i=1}^K \degWiv$.  
\begin{enumerate}
\item[\textbf{(a)}] If \underline{either} $V>k(K+1-1/\Dmin)-1$ and $R \coloneqq \max\{k(KD+1), (Kk-1)D_0+2\}$, \underline{or}
\item[\textbf{(b)}] If $q$ is squarefree, $V \ge Kk$, and $R \coloneqq k(2K+1)$,
\end{enumerate} 
then the relations \eqref{eq:WUDPRn} hold, uniformly in coprime residues $a_1, \dots, a_K$ modulo $\qlelogxKZ$ lying in $\WUDkAdmSet$ and satisfying $\IFHWikBZ$.
\end{thm} 
Notice that for any $K>2$, the result under (b) 
unconditionally improves over Theorem \ref{thm:RestrictedInputSqfreeMod1}(a) in terms of weakening the restriction on inputs $n$. On the other hand, the result under condition (a) 
improves over Theorem \ref{thm:RestrictedInputGenMod1}(a) whenever $k$ or $D$ is large enough compared to $D_0$. 

We now explain the necessity of the two key additional hypotheses that we have been assuming in our main results so far, namely the multiplicative independence of $\Wikset \subZT$ and the invariant factor hypothesis. It turns out that without the former condition, the $K$ congruences $\finaimodq$ (for $1 \le i \le K$) may degenerate to fewer congruences for sufficiently many inputs $n$, making weak equidistribution fail uniformly to \textit{all} sufficiently large $q \le (\log x)^{K_0}$. In this situation, weak equidistribution \textit{cannot} be restored \textit{no matter} how much we restrict the set of inputs $n$ to those having sufficiently many large prime factors. We make this explicit in the next result. 
\begin{thm}\label{thm:MIHNecessity}
Fix $R \ge 1$, $K>1$ and assume that $\{\Wik\}_{1 \le i \le K-1} \subZT$ are multiplicatively independent, with $\sum_{i=1}^{K-1} \degWik>1$. Suppose $W_{K, k}$ $= \prod_{i=1}^{K-1} \Wik^{\lambda_i}$ for some nonnegative integers $(\lambda_i)_{i=1}^{K-1} \ne (0, \dots, 0)$. There exists a constant $C \coloneqq C(W_{1, k}, \dots, W_{K-1, k})>0$ such that 
\begin{multline*} 
\#\{n \le x: P_{Rk}(n)>q, ~ (\forall i \in [K]) ~ \finaimodq\} \gg \frac1{\phi(q)^{K-1}} \cdot \frac{x^{1/k} (\log \log x)^{R-2}}{\log x}
\end{multline*}
as $x \rightarrow \infty$, uniformly in $k$-admissible $\qlelogxKZ$ supported on primes $\ell>C$ satisfying $\gcd(\ell-1, \beta(W_{1, k},$ $\dots,$ $W_{K-1, k}))=1$, and in $a_i \in U_q$ with $a_K \equiv \prod_{i=1}^{K-1} a_i^{\lambda_i} \pmod q$. 
\end{thm} 
The compatibility of the relations involving $W_{K. k}$ and $a_K$ suggests the aforementioned degeneracy from $K$ to $K-1$ congruences. Note that the above lower bound will in fact come from the $n$ which are supported on primes much larger than $q$. A similar lower bound holds for $K=1$ when $W_k = W_{1, k}$ is constant (see the remark preceding subsection \cref{subsec:MIHIFH_Egs}). Using the above theorem, 
we shall construct (in \cref{subsec:MIHIFH_Egs}) explicit examples of polynomials $\{W_{i, k}\}_{1 \le i \le K-1}$ and moduli $q \in \WUDkAdmSet$ where the above lower bound grows strictly faster than the expected proportion of $n \le x$ having $\gcd(f(n), q)=1$. This would demonstrate an overrepresentation of the coprime residues $(a_i \bmod q)_{i=1}^K$ by the multiplicative functions $f_1, \dots, f_K$, coming from inputs $n$ that have at least $Rk$ many prime factors exceeding $q$, showing the necessity of our hypothesis on the multiplicative independence of $\Wikset$. 

Turning to the invariant factor hypothesis, we show that the failure of this condition incurs an additional factor over the expected proportion of $n \le x$ satisfying $\gcd(f(n), q)=1$.  
For certain choices of $q$ and $\Wikset$, this factor can be made too large, once again leading to an overrepresentation of the tuple $(a_i \bmod q)_{i=1}^K$ by the multiplicative functions $f_1, \dots, f_K$. In what follows, $P^-(q)$ denotes the smallest prime dividing $q$.
\begin{thm}\label{thm:IFHNecessity}
Fix $R \ge 1$ and assume that $\Wikset \subZT$ are nonconstant, monic and multiplicatively independent, so that $\beta = \betaWiklist \in \Z\sm\{0\}$.  
There exists a constant $C \coloneqq C(\Wiklist) > 0$ such that
\begin{equation}\label{eq:IFHNecess_Estim}
\#\{n \le x: P_{Rk}(n)>q, ~ (\forall i \in [K]) ~ \finaimodq\} \gg \frac{2^{\#\{\ell \mid q: ~ \gcd(\ell-1, \beta) \ne 1\}}}{\phi(q)^K} \cdot \frac{x^{1/k} (\log \log x)^{R-2}}{\log x}
\end{equation}
as $x \rightarrow \infty$, uniformly in $k$-admissible $\qlelogxKZ$ having $P^-(q)>C$, and in coprime residues $(a_i)_{i=1}^K$ mod $q$ which are all congruent to $1$ modulo the largest squarefree divisor of $q$.
\end{thm}
Here, the restriction on the residues $a_i$ is imposed in order to  have a positive contribution of certain character sums modulo the prime divisors of $q$. In subsection \cref{subsec:MIHIFH_Egs}, we shall construct explicit examples of $q \in \WUDkAdmSet$ and $\Wikset$ for which the expression in the above lower bound is much larger than the expected proportion of $n \le x$ having $\gcd(f(n), q)=1$. 

We can give several applications of our main results to arithmetic functions of common interest. For instance,  recall \'Sliwa's \cite{sliwa73} result that $\sigma(n)$ is weakly equidistributed precisely to moduli that are not multiples of $6$; in fact, his result shows that $\mathcal Q(1; \sigma) = \{q: \gcd(q, 2)=1\}$ and $\mathcal Q(2; \sigma) = \{q: \gcd(q, 6)=2\}$. By Theorem \ref{thm:UnrestrictedInput_SqfrPoly}(i), 
$\sigma(n)$ is WUD uniformly to all odd moduli $\qlelogxKZ$. Calling the members of the set $\mathcal Q(2; \sigma)$ ``special", Theorem \ref{thm:UnrestrictedInput_SqfrPoly}(ii) and (iii) show that $\sigma(n)$ is WUD uniformly to all special $q \le (\log x)^{(2-\delta)\alphatil(q)}$ and also to all squarefree special $\qlelogxKZ$ satisfying $2^{\omega(q)} \le (\log x)^{(1-\epsilon)\alphatil(q)}$, where $\alphatil(q) \coloneqq \alpha_2(q) = \prod_{\substack{\ell \mid q\\ \ell \equiv 1 \pmod 3}} (1-2/(\ell-1))$. By the example constructed in \cite[subsection 7.1]{SR23A}, the latter restriction is optimal.  
Furthermore, by Theorems 
\ref{thm:RestrictedInputGenMod1}(a) or \ref{thm:HigherPolyControl1}(a) (resp. by Theorem \ref{thm:RestrictedInputSqfreeMod1}(a)), uniformity is restored to \textit{all} (resp. to squarefree) special $\qlelogxKZ$ by restricting to inputs $n$ with $P_6(n)>q$ (resp. $P_4(n)>q$); here we have noted that the condition $P_3(n)>q$ forces $P_4(n)>q$ since for $\sigma(n)$ to be coprime to the even number $q$, it is necessary for $n$ to be of the form $m^2$ or $2m^2$. By the examples constructed in \cite{SR23A}, both of these restrictions are optimal as well. 
Alternatively, by Theorem \ref{thm:RestrictedInputGenMod1}(b) (resp. \ref{thm:RestrictedInputSqfreeMod1}(b)), we may restrict to $n$ with $P_3(n_2)>q$ (resp. $P_2(n_2)>q$) to restore complete uniformity in all (resp. squarefree) special $\qlelogxKZ$. 

For another example, we saw using Theorem \ref{thm:UnrestrictedInput_SqfrPoly} that $\phi(n)$ and $\sigma(n)$ are jointly WUD modulo $q \le (\log x)^{(1-\epsilon)\alpha(q)}$ coprime to $6$, and that these two restrictions on $q$ are necessary and essentially optimal. By Theorem \ref{thm:RestrictedInputGenMod1}, 
complete uniformity is restored to all moduli $\qlelogxKZ$ coprime to $6$ by restricting to inputs $n$ with $P_5(n)>q$. 

We can give more applications of our main results to study the weak equidistribution of the Fourier coefficients of Eisenstein series; precisely, the functions $\sigma_r(n) \coloneqq \sum_{d \mid n} d^r$ (for $r>1$). An easy check shows that the polynomial $\sum_{0 \le j \le v} T^{rj} = \frac{T^{r(v+1)}-1}{T^r-1}$ shares no roots with its derivative, hence is separable. Calling the $q \in \mathcal Q(k; \sigma_r)$ as ``$k$-good", Theorem \ref{thm:UnrestrictedInput_SqfrPoly} thus shows that $\sigma_r$ is WUD uniformly to all $k$-good $q \le (\log x)^{(1-\epsilon)\alpha_k(q)(1-1/kr)^{-1}}$, and to all squarefree $k$-good $\qlelogxKZ$ having $\omega(q) \le (1-\epsilon) \alpha_k(q) \log \log x/\log(kr)$. Further, by Theorems \ref{thm:RestrictedInputGenMod1} to \ref{thm:HigherPolyControl1}, weak equidistribution is restored modulo all $k$-good $\qlelogxKZ$ by restricting to $n$ with $P_{k(kr+1)}(n)>q$, whereas it is restored modulo all squarefree $k$-good $\qlelogxKZ$ by restricting to $n$ with $P_{k+1}(n)>q$ or to $n$ with $P_2(n_k)>q$. An explicit characterization of the moduli $\qlelogxKZ$ to which a given $\sigma_r$ is weakly equidistributed thus reduces to an understanding of the possible $k$ and of the set $\mathcal Q(k; \sigma_r)$ for a given (fixed) $r$; both of these are problems of fixed moduli that (as mentioned in the introduction) 
have been studied in great depth in \cite{sliwa73}, \cite{fomenko80}, \cite{narkray82}, \cite{narkiewicz83}, \cite{narkiewicz84}, \cite{rayner88} and \cite{rayner88a}. In fact, the sets $\mathcal Q(k; \sigma_r)$ have been explicitly characterized for all odd $r \le 200$ and all even $r \le 50$, and partial results are known for general $r \ge 4$. For example, the only two possible $k$'s for $\sigma_3$ are $k=1, 2$, and $\mathcal Q(1; \sigma_3) = \{q: \gcd(q, 14)=1\}$ while $\mathcal Q(2; \sigma_3) = \{q: \gcd(q, 6)=2\}$.

For a general family $(f_1, \dots, f_K)$, Narkiewicz \cite{narkiewicz81, narkiewicz84} gives algorithms to determine the sets $\WUDkAdmSet$ for a fixed $k$. He shows (among other results) that in some of the most commonly occurring cases (which includes the cases of $\sigma_r$ for all $r>2$), the set of possible $k$ is finite, and that for each such $k$, the set $\mathcal Q(k; f_1, \dots , f_K)$ can be characterized by certain (finitely many) coprimality conditions that can be determined effectively. 

We conclude this section with the remark that although for the sake of simplicity of statements, we have been assuming 
that our multiplicative functions $\{f_i\}_{i=1}^K$ and polynomials $\Wivfullset$ are both fixed, our proofs will reveal that these results are also uniform in the 
$\{f_i\}_{i=1}^K$ as long as they are defined by the fixed polynomials $\Wivfullset$. 
\subsection*{Notation and conventions:} We do not consider the zero function as multiplicative (thus, if $f$ is multiplicative, then $f(1)=1$). 
Given $z>0$, we say that a positive integer $n$ is \textsf{$z$-smooth} if $P(n) \le z$, and \textsf{$z$-rough} if $P^-(n) > z$; by the \textsf{$z$-smooth part} (resp. \textsf{$z$-rough part}) of $n$, we shall mean the largest $z$-smooth (resp. $z$-rough) positive integer dividing $n$. For a ring $R$, we shall use $R^\times$ to denote the multiplicative group of units of $R$. We denote the number of primes dividing $q$ counted with and without multiplicity by $\Omega(q)$ and $\omega(q)$ respectively, and we write $U_q \coloneqq (\Z/q\Z)^\times$. For a Dirichlet character $\chi$ mod $q$, we use $\condofchi$ to denote the conductor of $\chi$. When there is no danger of confusion, we shall write $(a_1, \dots, a_k)$ in place of $\gcd(a_1, \dots, a_k)$. Throughout, the letters $p$ and $\ell$ are reserved for primes. For nonzero $H \in \Z[T]$, we use $\ord_\ell(H)$ to denote the highest power of $\ell$ dividing all the coefficients of $H$; for an integer $m \ne 0$, we shall sometimes use $v_\ell(m)$ in place of $\ord_\ell(m)$. We use $\mathbb M_{A \times B}(\Z)$ to refer to the ring of $A \times B$ matrices with integer entries, while $GL_{A \times B}(\Z)$ refers to the group of units of $\mathbb M_{A \times B}(\Z)$, i.e. the matrices with determinant $\pm 1$. Implied constants in $\ll$ and $O$-notation, as well as implicit constants in qualifiers like ``sufficiently large", may always depend on any parameters declared as ``fixed''; in particular, they will always depend on the polynomials $\{\Wiv\}_{\substack{1 \le i \le K\\1 \le v \le k}}$. Other dependence will be noted explicitly (for example, with parentheses or subscripts); notably, we shall use $C(\Filist)$, $C'(\Filist)$ and so on, to denote constants depending on the fixed polynomials $\Filist$. 
We write $\log_{k}$ for the $k$-th iterate of the natural logarithm. 
\section{Technical preparation: The number of $n \le x$ for which $\gcd(f(n), q)=1$}\label{sec:TechnicalPrep}
In this section, we shall provide a rough estimate on the count of $n \le x$ for which $f(n) = \prodik f_i(n)$ is coprime to the modulus $q$, uniformly in $\qlelogxKZ$. We aim to show the following estimate, which generalizes Proposition 2.1 in \cite{PSR23}. In the rest of the paper, we abbreviate $\alpha_v(q)$ to $\alpha_v$ for each $v \in [V]$.
\begin{prop}\label{prop:fnqcoprimecount}
For all sufficiently large $x$ and uniformly in $k$-admissible $\qlelogxKZ$, 
\begin{equation}\label{eq:fnqcoprimeEstimate}
\sumnxfnq ~ = \largesum_{\substack{n \le x\\\text{each} ~ (f_i(n),q)=1}} 1 = \frac{x^{1/k}}{(\log x)^{1-\alpha_k}} \exp(O((\log_2 (3q))^{O(1)})).
\end{equation}
\end{prop}
\subsection{Proof of the lower bound.}
Any $m \le x^{1/k}$ satisfying $\gcd(f(m^k), q) = 1$ is certainly counted in the left hand side of \eqref{eq:fnqcoprimeEstimate}. To estimate the number of such $m$, we apply Proposition \cite[Proposition 2.1]{PSR23}, with $f(n^k)$ and $x^{1/k}$ playing the roles of ``$f(n)$" and ``$x$" in the quoted proposition. 
This shows that the sum in \eqref{eq:fnqcoprimeEstimate} is bounded below by the right hand side. 
\subsection{Proof of the upper bound.}
We start by giving an upper bound on the count of $r$-full smooth numbers; 
here we consider any $n \in \NatNos$ to be $1$-full (and we consider $1$ as being $r$-full for any $r \ge 1$). The case $r=1$ of the lemma below is a known estimate on smooth numbers.
\begin{lem}\label{lem:rfullsmooth}
Fix $r \in \NatNos$. We have as $X, Z \to \infty$, 
$$\#\{n \le X: P(n) \le Z, ~ n \text{ is }r\text{-full}\} \ll X^{1/r}(\log Z) \exp\left(-\frac Ur \log U + O(U \log_2(3U))\right),$$
uniformly for $(\log X)^{\max\{3, 2r\}} \le Z \le X^{1/2}$, where $U \coloneqq {\log X}/{\log Z}$.
\end{lem}
\begin{proof}[Proof of Lemma \ref{lem:rfullsmooth}]
The lemma is a classical application of Rankin's trick. We start by letting $\eta \le \min\{1/3, 1/2r\}$ be a positive parameter to be chosen later, and observe that 
\begin{align}\label{eq:RankinApp1}\allowdisplaybreaks
\largesum_{\substack{n \le X: ~ P(n) \le Z\\n \text{ is }r\text{-full}}} 1 \le \largesum_{\substack{n\text{ is }r\text{-full}\\P(n) \le Z}} \left(\frac Xn\right)^{(1-\eta)/r} 
\ll X^{(1-\eta)/r} \exp\left(\largesum_{p \le Z} \frac1{p^{1-\eta}}\right), 
\end{align} 
where we have used the Euler product and noted that $\sum_p \sum_{v \ge r+1} p^{-v(1-\eta)/r} \ll \sum_p p^{-(1-\eta)(1+1/r)}$ $\ll_r 1$ since $(1-\eta)(1+1/r) \ge (1+1/r)(1-\min\{1/3, 1/2r\}) > 1$. 

Let $\eta \coloneqq \frac{\log U}{\log Z} \le \min\left\{\frac13, \frac1{2r}\right\}$. We write $\sum_{p \le Z} 1/{p^{1-\eta}} = \log_2 Z + \largesum_{p \le Z} (\exp(\eta \log p)-1)/p + O(1)$. Since $\eta\log p \le \log 2 \ll 1$ for all $p \le 2^{1/\eta}$, we find that the contribution of $p \le 2^{1/\eta}$ to the last sum above is 
$ \ll \eta\sum_{p \le 2^{1/\eta}} {\log p}/p \ll 1$, while that of $p \in (2^{1/\eta}, Z]$ is at most $(\exp(\eta \log Z)-1) \sum_{2^{1/\eta} < p \le Z} 1/p 
~ ~ \le ~ U(\log_2 U + O(1))$.  
Collecting estimates, we obtain $\sum_{p \le Z} 1/{p^{1-\eta}} = \log_2 Z + O(U \log_2 (3U)),$ which from \eqref{eq:RankinApp1} completes the proof of the lemma.
\end{proof}
Since $\alpha_v(\ell)>0$ for all $\ell>D_v+1$, it follows that for each $1 \le v < k$, the set $S_v \coloneqq \{\ell\text{ prime }: \alpha_v(\ell)=0\}$ consists only of primes of size $O(1)$, with the implied constant depending only on the polynomials $\Wivlist$. It is easy to show that if $q$ is $k$-admissible, then the $k$-free part of any positive integer $n$ satisfying $\gcd(f(n), q)=1$ must be supported on the primes in the set $\bigcup_{1 \le v<k} S_v$. As a consequence, we have the following important observation. 
\begin{lem}\label{lem:kfreepartbdd}
If $q$ is $k$-admissible, then the $k$-free part of any positive integer $n$ satisfying $\gcd(f(n), q)=1$ is bounded. More precisely, it is of size $O(1)$, where the implied constant depends only on the polynomials $\{\Wiv\}_{\substack{1 \le i \le K\\1\le v \le k}}$.
\end{lem}
The following estimate (see \cite[Lemma 2.4]{PSR23}) will be useful throughout the paper.
\begin{lem}\label{lem:primesum} Let $G\in \Z[T]$ be a fixed nonconstant polynomial. For each positive integer $q$, let $\alpha_G(q) \coloneqq \frac1{\phi(q)}\#\{u \in U_q: G(u) \in U_q\}$. We have, uniformly in $q$ and $x\ge 3q$, \[ \sum_{p \le x} \frac{\1_{(G(p),q)=1}}p = \alpha_G(q) \log_2 x  + O((\log_2{(3q)})^{O(1)}).  \]
\end{lem}
Coming to the proof of the upper bound implied in \eqref{eq:fnqcoprimeEstimate}, we define $y \coloneqq \exp(\sqrt{\log x})$ and start by removing those $n$ which are divisible by the $(k+1)$-th power of a prime exceeding $y$. Writing any such $n$ as $AB$ for some $k$-free $B$ and $k$-full $A$, Lemma \ref{lem:kfreepartbdd} shows that $B\ll 1$ so that the contribution of such $n$ to $\eqref{eq:fnqcoprimeEstimate}$ is 
\begin{align}\label{eq:Large power of prime exceeding y negligible}\allowdisplaybreaks
\largesum_{\substack{n \le x: ~ (f(n), q)=1\\\exists ~ p>y: ~ p^{k+1} \mid n}} 1 ~ \ll \largesum_{\substack{A \le x\\A \iskfull\\\exists ~ p>y: ~ p^{k+1} \mid n}} 1 &\le \largesum_{p>y} ~ ~ \largesum_{\substack{v \ge k+1\\ p^v \le x}} ~ ~ \largesum_{\substack{m \le x/p^v\\m \iskfull}} 1 \ll \largesum_{p>y} ~ ~ \largesum_{v \ge k+1} \left(\frac x{p^v}\right)^{1/k} 
\ll \left(\frac xy\right)^{1/k},   
\end{align}
where we have used the fact that the number of $k$-full integers up to $X$ is $O(X^{1/k})$ (see \cite{ES34}). 
The last expression above is negligible in comparison to the right hand side of \eqref{eq:fnqcoprimeEstimate}, hence, it remains to bound the number of $n$ not divisible by the $(k+1)$-th power of any prime greater than $y$ and satisfying $(f(n), q)=1$. 

We write any such $n$ in the form $BMN$, where $N$ is $y$-rough, $BM$ is $y$-smooth, $B$ is $k$-free, $M$ is $k$-full, and $B, M, N$ are pairwise coprime. By Lemma \ref{lem:kfreepartbdd}, we see that $B=O(1)$ and that $N$ is $k$-full. But since $n$ is not divisible by the $(k+1)$-th power of any prime exceeding $y$, $N$ must be the $k$-th power of a squarefree $y$-rough integer $A$. Consequently, 
\begin{equation}\label{eq:BMAkSplit1}
\largesum_{\substack{n \le x: ~ (f(n), q)=1\\ \nopgrystpkOnedivn}} 1 ~ ~ ~ \le \largesum_{\substack{B \le x\\ (f(B), q)=1\\B \iskfree}} ~ ~ ~ \largesum_{\substack{M \le x/B: ~ M\iskfull\\ P(M) \le y, ~ (f(M), q)=1}} ~ ~ \largesum_{\substack{A \le (x/BM)^{1/k}\\P^-(A)>y: ~ (f(A^k), q)=1\\A \issqfree}} 1.
\end{equation}
We now write the right hand side of the above inequality as $\Sigma_1+\Sigma_2$, where $\Sigma_1$ and $\Sigma_2$ count the contribution of $(B, M, A)$ with $M \le x^{1/2}$ and $M>x^{1/2}$, respectively. Any $A$ counted in $\Sigma_2$ satisfies $A \le (x/BM)^{1/k} \le x^{1/2k}/B^{1/k}$, so that 
$$\Sigma_2 \le \largesum_{\substack{B \le x\\ (f(B), q)=1\\B \iskfree}} ~ ~ ~ \largesum_{\substack{A \le x^{1/2k}/B^{1/k}\\P^-(A)>y: ~ (f(A^k), q)=1\\A \issqfree}} ~ ~ \largesum_{\substack{M \le x/BA^k: ~ P(M) \le y\\ M\iskfull, ~ (f(M), q)=1}} 1.$$
To bound the innermost sum, we invoke Lemma \ref{lem:rfullsmooth}; 
here $U = \frac{\log(x/BA^k)}{\log y} \ge \frac12\sqrt{\log x}$. This yields
$$\Sigma_2 \ll \largesum_{\substack{B \le x\\ (f(B), q)=1\\B \iskfree}} ~ ~ ~ \largesum_{\substack{A \le x^{1/2k}/B^{1/k}\\P^-(A)>y: ~ (f(A^k), q)=1\\A \issqfree}} ~ \frac{x^{1/k}}{B^{1/k} A} \exp\left(-\frac1{6k} \sqrt{\log x} \cdot \log_2 x\right).$$
Recalling that $B=O(1)$ and bounding the sum on $A$ trivially by $2\log x$, we deduce that $\Sigma_2 \ll x^{1/k}\exp\left(-\sqrt{\log x}\right)$, 
which is negligible compared to the right hand side of \eqref{eq:fnqcoprimeEstimate}. 

To estimate $\Sigma_1$, we invoke \cite[Theorem 01, p.\ 2]{HT88} 
on the multiplicative function $g(A) \coloneqq \mu(A)^2 \bbm_{P^-(A)>y} \bbm_{(f(A^k), q)=1}$, with $\mu$ denoting the M\"obius function. Since $M \le x^{1/2}$ and $B \ll 1$, 
\begin{align*}\allowdisplaybreaks
\Sigma_1 
\ll \frac{x^{1/k}}{\log x} \exp\left(\largesum_{y < p \le x} \frac{\bbm_{(W_k(p), q)=1}}p\right) \largesum_{\substack{M \le x^{1/2}: ~ M\iskfull\\ P(M) \le y, ~ (f(M), q)=1}} \frac1{M^{1/k}}. 
\end{align*}
But since the sum on $M$ above is no more than 
\begin{equation}\label{eq:Msmoothrecip}
\largesum_{\substack{M \iskfull\\ P(M) \le y, ~ (f(M), q)=1}} \frac1{M^{1/k}} \le \prod_{p \le y} \left(1+\frac{\bbm_{(f(p^k), q)=1}}p + O\left(\frac1{p^{1+1/k}}\right)\right) \ll \exp\left(\largesum_{p \le y} \frac{\bbm_{(W_k(p), q)=1}}p\right),
\end{equation}
it follows 
by an application of Lemma \ref{lem:primesum} to estimate the sum $\sum_{p \le x} {\bbm_{(W_k(p), q)=1}}/p$, that $\Sigma_1$ is absorbed in the right hand side of \eqref{eq:fnqcoprimeEstimate}. This establishes Proposition \ref{prop:fnqcoprimecount}. 

\section{The main term in Theorems \ref{thm:UnrestrictedInput_SqfrPoly} to \ref{thm:HigherPolyControl1}: Contribution of ``convenient" $n$}\label{sec:Convenientn}
In what follows, we define 
$$J \coloneqq \lfloor \log_3 x \rfloor \text{ and } y \coloneqq \exp((\log x)^{\epsilon/2}),$$
where $\epsilon$ is as in the statement of Theorem \ref{thm:UnrestrictedInput_SqfrPoly} and $\epsilon \coloneqq 1$ for Theorems \ref{thm:RestrictedInputGenMod1} to \ref{thm:HigherPolyControl1}. We call $n \le x$ \textsf{convenient} if the largest $J$ \textit{distinct} prime divisors of $n$ exceed $y$ and each appear to exactly the $k$-th power in $n$. In other words, $n$ is convenient iff it can be uniquely written in the form $n = m(P_J \cdots P_1)^k$ for $m \le x$ and primes $P_1, \dots, P_J$ satisfying 
\begin{equation}\label{eq:convnDefnIneq}
L_m \coloneqq \max\{y, P(m)\} < P_J < \dots < P_1.
\end{equation}
Note that any $n$ having $P_{Jk}(n) \le y$ must be inconvenient; on the other hand, if $n$ is inconvenient and satisfies $\gcd(f(n), q)=1$ then either $P_{Jk}(n) \le y$ or $n$ is divisible by the $(k+1)$-th power of a prime exceeding $y$. 
We start by showing that there are a negligible number of inconvenient $n \le x$ satisfying $\gcd(f(n),q)=1$.
\begin{prop}\label{prop:inconvnfnq=1}
We have as $x \to\infty$, 
\begin{equation}\label{eq:inconvfnq=1}
\largesum_{\substack{n \le x: ~(f(n), q)=1\\n\inconv}} 1 ~ ~ = ~ \osumnxfnq,
\end{equation}
uniformly in $k$-admissible $\qlelogxKZ$.
\end{prop}
\begin{proof} 
By \eqref{eq:Large power of prime exceeding y negligible} and \eqref{eq:fnqcoprimeEstimate}, 
the contribution of the $n$'s that are divisible by the $(k+1)$-th power of a prime exceeding $y$ is negligible. 
Letting $z \coloneqq x^{1/\log_2 x}$, we show that the same is true for the contribution of $z$-smooth $n$ to the left hand side of \eqref{eq:inconvfnq=1}. Indeed, writing any such $n$ in the form $AB$ for some $k$-free $B$ and $k$-full $A$, we have $P(A) \le z$ whereas (by Lemma \ref{lem:kfreepartbdd}) $B=O(1)$. 
Hence the contribution of $z$-smooth $n$ 
is, by Lemma \ref{lem:rfullsmooth}, 
\begin{equation}\label{eq:Smoothfnq=1}
\largesum_{\substack{n \le x: ~P(n) \le z\\(f(n), q)=1}} 1 ~  ~ \ll \largesum_{\substack{A \le x: ~ P(A) \le z\\A \iskfull}} 1 ~ ~ \ll x^{1/k} \exp\left(-\left(\frac1k+o(1)\right)\log_2 x \log_3 x\right),
\end{equation}
which is negligible compared to the right hand side of \eqref{eq:inconvfnq=1}. 

It remains to consider the contribution of those $n$ which are not $z$-smooth and are not divisible by the $(k+1)$-th power of a prime exceeding $y$. Since $n$ is inconvenient, we have $P_{Jk}(n) \le y$. Hence, $n$ can be written in the form $mP^k$ where $P \coloneqq P(n) > z$ and $m  = n/P^k$, so that $P_{Jk}(m) \le y$, $\gcd(m, P)=1$ and $f(n)= f(m) f(P^k)$. Given $m$, there are at most $\sum_{z < P \le (x/m)^{1/k}} 1 \ll x^{1/k}/m^{1/k}\log z$ many possibilities for $P$. Consequently, 
\begin{equation}\label{eq:Inconvnfnq1Main}
\largesum_{\substack{n \le x \inconv\\P(n)>z, ~ (f(n), q)=1\\ \nopgrystpkOnedivn}} 1 ~ ~ \le \largesum_{\substack{n \le x: ~P_{Jk}(n) \le y\\P(n)>z, ~ (f(n), q)=1\\ \nopgrystpkOnedivn}} 1 ~ ~ \ll ~ \frac{x^{1/k} \log_2 x}{\log x} \largesum_{\substack{m \le x\\P_{Jk}(m) \le y, ~ (f(m), q)=1\\ \nopgrystpkOnedivm}} \frac1{m^{1/k}}.
\end{equation}
As in the argument preceding \eqref{eq:BMAkSplit1}, we write any $m$ occurring in the above sum (uniquely) in the form $BMA^k$, where $B$ is $k$-free, $M$ is $k$-full, $A$ is squarefree, $P(BM) \le y < P^-(A)$, 
and $\Omega(A)\le J$ (since $P_{Jk}(n) \le y$). 
Since $B=O(1)$, we deduce that 
\begin{align*}\allowdisplaybreaks
\largesum_{\substack{m \le x\\P_{Jk}(m) \le y, ~ (f(m), q)=1\\ \nopgrystpkOnedivm}} \frac1{m^{1/k}} 
\ll \largesum_{\substack{M ~ k\text{-full}\\ P(M) \le y, ~ (f(M), q)=1}} \frac1{M^{1/k}} ~ ~ \largesum_{\substack{A \le x\\ ~ \Omega(A) \le J}} \frac1A. 
\end{align*}
The sum on $A$ is no more than $ (1+\sum_{p \le x} 1/p)^J \le (2 \log_2 x)^J \le \exp(O((\log_3 x)^2))$, while the sum on $M$ is $\ll \exp(\alpha_k \log_2 y + O((\log_2 (3q))^{O(1)}))$ by \eqref{eq:Msmoothrecip} and Lemma \ref{lem:primesum}. Altogether,
\begin{equation}\label{eq:mInconvRecip}
\largesum_{\substack{m \le x\\P_{Jk}(m) \le y, ~ (f(m), q)=1\\ \nopgrystpkOnedivm}} \frac1{m^{1/k}} \ll (\log x)^{\alpha_k \epsilon/2} \exp\big(O((\log_3 x)^2 + (\log_2 (3q))^{O(1)})\big),
\end{equation}
and inserting this into \eqref{eq:Inconvnfnq1Main} completes the proof via 
Proposition \ref{prop:fnqcoprimecount}. 
\end{proof} 
It is the convenient $n$ which give rise to the main term in the count of $n \le x$ satisfying the congruences $\finaimodq$. 
We shall spend the next few sections proving this.
\begin{thm}\label{thm:ConvenientMAINTerm}
Fix $K_0, B_0 > 0$ and assume that $\Wikset \subZT$ are nonconstant and multiplicatively independent. As $x \to\infty$, we have
\begin{equation*}
\largesum_{\substack{n \le x \conv\\\forallifinaimodq}} 1 ~ ~ \sim \frac1{\phi(q)^K} \sumnxfnq,
\end{equation*}   
uniformly in coprime residues $a_1, \dots, a_K$ to moduli $\qlelogxKZ$ lying in $\WUDkAdmSet$ and satisfying $\IFHWikBZ$.
\end{thm}
In this section and the next, we establish a weaker version of this result, where we reduce the congruences $\finaimodq$ to a bounded modulus. 
\begin{prop}\label{prop:ConvenientMainTermUptoQ0}
Fix $K_0, B_0 > 0$ and assume that $\Wikset \subZT$ are nonconstant and multiplicatively independent. There exists a constant $\lambda \coloneqq \lambda(\Wiklist; B_0)>0$ depending only on $\Wikset \subZT$ and $B_0$, such that as $x \to\infty$, we have
\begin{equation}\label{eq:CountReducedToQ0}
\largesum_{\substack{n \le x \conv\\\forallifinaimodq}} 1 ~ ~ = \phiQZphiqK \sumnxfnqfiQZ ~ ~ + ~ \ophiqKsumnxfq,
\end{equation}   
uniformly in coprime residues $a_1, \dots, a_K$ to $k$-admissible moduli $\qlelogxKZ$ satisfying $\IFHWikBZ$, where $Q_0$ is a divisor of $q$ satisfying $Q_0 \le \lambda$.
\end{prop}
\begin{proof} For any $N \ge 1$ and $\wifam \in U_q^K$, we define $$\VNKqwi \coloneqq \left\{(v_1, \dots, v_N) \in (U_q)^N: ~ (\forall i \in [K]) ~ \prod_{j=1}^N \Wik(v_j) \equiv w_i \pmod q\right\}.$$
We write each convenient $n$ uniquely in the form $m(P_J \cdots P_1)^k$, where $m, P_J, \dots, P_1$ satisfy \eqref{eq:convnDefnIneq}. Then $f_i(n) = f_i(m) \prod_{j=1}^J \Wik(P_j)$, so that the conditions $\finaimodq$ amount to $\gcd(f(m), q)=1$ and $(P_1, \dots, P_J)$ mod $q$ $\in \VqmPr \coloneqq \VJKqaifim$.  
Noting that the conditions $P_1 \cdots P_J \le (x/m)^{1/k}$ and $(P_1, \dots, P_J)$ mod $q$ $\in \VqmPr$ are both independent of the ordering of $P_1, \dots, P_J$, we obtain
\begin{align}\label{eq:ConvnInitialSplit}\allowdisplaybreaks
\sum_{\substack{n \le x\conv\\ \forallifinaimodq}} 1 ~ ~ 
= \largesum_{\substack{m \le x\\(f(m), q)=1}} ~\sum_{(v_1, \dots, v_J) \in \VqmPr} \frac1{J!} \sum_{\substack{P_1, \dots, P_J > L_m\\P_1 \cdots P_J \le (x/m)^{1/k}\\ P_1, \dots, P_J \text{ distinct }\\(\forall j) ~ P_j \equiv v_j \pmod q}} 1. \nonumber
\end{align}

Proceeding exactly as in \cite{PSR23} to remove the congruence conditions on $P_1,\dots,P_J$ by successive applications of the Siegel--Walfisz Theorem, we deduce that
\begin{equation}\label{eq:CongrueceRemoved}
\sum_{\substack{P_1, \dots, P_J > L_m\\P_1 \cdots P_J \le (x/m)^{1/k}\\ P_1, \dots, P_J \text{ distinct }\\(\forall j) ~ P_j \equiv v_j \pmod q}} 1 ~ = ~  \frac1{\phi(q)^J}\sum_{\substack{P_1, \dots, P_J > L_m\\P_1 \cdots P_J \le (x/m)^{1/k}\\ P_1, \dots, P_J \text{ distinct }}} 1 ~ + ~ O\left(\frac{x^{1/k}}{m^{1/k}} \exp\left(-{K_1} (\log{x})^{\epsilon/4}\right)\right)
\end{equation} 
for some constant $K_1 \coloneqq K_1(K_0)>0$. Collecting estimates and noting that $\#\VqmPr \le \phi(q)^J \le (\log x)^{K_0 J}$, we obtain
\begin{equation}\label{eq:convnSplitForm} 
\sum_{\substack{n \le x\conv\\\forallifinaimodq}} 1 = \largesum_{\substack{m \le x\\(f(m), q)=1}} \frac{\#\VqmPr}{\phi(q)^J} \Bigg(\frac1{J!}\sum_{\substack{P_1, \dots, P_J > L_m\\P_1 \cdots P_J \le (x/m)^{1/k}\\ P_1, \dots, P_J \text{ distinct }}} 1\Bigg) ~ + ~ O\left(x^{1/k} \exp\left(-\frac{K_1}2 (\log{x})^{\epsilon/4}\right)\right).    
\end{equation}
Here in the last step we have crudely bounded the sum $\sum_{\substack{m \le x\\(f(m), q)=1}} m^{-1/k}$ 
by writing each $m$ as $AB$ for some $k$-full $A$ and $k$-free $B$ satisfying $\gcd(A, B)=1$, and then noting that the sum $\sum 1/A$ is no more than $\prod_{p \le x} \left(1+1/p+O\left(1/{p^{1+1/k}}\right)\right)$.
The following proposition is a special case of the more general Proposition \ref{prop:VNKCountmegageneral} established in the next section, and will provide the needed estimate on the cardinalities of the sets $\VqmPr$.  
\begin{prop}\label{prop:Vqwi_ReducnToBddModulus}
Assume that $\Wikset$ are multiplicatively independent. 
There exists a constant $C_0 \coloneqq C_0(\Wiklist; B_0)> (8D)^{2D+2}$ depending only on $\Wikset$ and $B_0$, such that for \underline{any} constant $C>C_0$, the following estimates hold uniformly in coprime residues $\wifam$ to moduli $q$ satisfying $\alpha_k(q) \ne 0$ and $\IFHWikBZ$: We have
\begin{multline}\label{eq:VNKkCount_largeN}
\frac{\#\VNKqwi}{\phi(q)^N}\\ = \alphakratioN \phiQZphiqK \left\{\frac{\#\VNKQZwi}{\phi(Q_0)^N} + O\left(\frac1{C^N}\right)\right\} \mathlarger{\prod}_{\substack{\ell \mid q\\\ell>C_0}} \left(1+O\left(\frac{(4D)^N}{\ell^{N/D-K}}\right)\right), 
\end{multline}
uniformly for $N \ge KD+1$, where $Q_0$ is a $C_0$-smooth divisor of $q$ of size $O_C(1)$. Moreover 
\begin{equation}\label{eq:VNKkCount_smallN}
\frac{\#\VNKqwi}{\phi(q)^N} \le \frac{\big(\prod_{\ell^e \parallel q} e\big)^{\bbm_{N=KD}}}{q^{N/D}} ~ \exp\left(O(\omega(q)\right)), ~ \text{ for each }1 \le N \le KD. 
\end{equation}
\end{prop}
Applying \eqref{eq:VNKkCount_largeN} with $N \coloneqq J \ge KD+1$, and with $C$ chosen to be a constant exceeding $2 C_0^{C_0}$, we see that 
\begin{equation*}
\frac{\#\VqmPr}{\phi(q)^J}\\ = (1+o(1)) \frac{\alpha_k(q)^J}{\alpha_k(Q_0)^J} \phiQZphiqK \left\{\frac{\#\mathcal V_{Q_0, m}'}{\phi(Q_0)^J} + O\left(\frac1{C^J}\right)\right\},
\end{equation*}
where $\mathcal V_{Q_0, m}'$ $\coloneqq \mathcal V_{J, K}^{(k)}\big(Q_0; \aifimfam\big)$ and we have noted that $\sum_{\substack{\ell \mid q\\\ell>C_0}} {(4D)^J}/{\ell^{J/D-K}}$ $\le \left({4D}/{C_0^{1/(2D+2)}}\right)^J = o(1)$. We insert this into \eqref{eq:convnSplitForm}, and observe that since $\alpha_k(q) \ne 0$,  $Q_0 \mid q$ and $Q_0$ is $C_0$-smooth, we have $\alpha_k(Q_0) C \ge C\prod_{\substack{\ell \le C_0}} \left(1-\frac{\ell-2}{\ell-1}\right) \ge \frac C{C_0^{C_0}} \ge 2$. We obtain
\begin{multline}\label{eq:ConvnSplitFormSimplified}
\largesum_{\substack{n \le x\conv\\ \forallifinaimodq}} 1 \\= (1+o(1)) \phiQZphiqK \frac{\alpha_k(q)^J}{\alpha_k(Q_0)^J} \largesum_{\substack{m \le x\\(f(m), q)=1}}\frac{\#\VQZmPr}{\phi(Q_0)^J}  \Bigg(\frac1{J!}\sum_{\substack{P_1, \dots, P_J > L_m\\P_1 \cdots P_J \le (x/m)^{1/k}\\ P_1, \dots, P_J \text{ distinct }}} 1\Bigg) + \ophiqKsumnxfq,
\end{multline}
where by the arguments leading to \eqref{eq:convnSplitForm} and the observation $\#\{(v_1, \dots, v_J) \in U_q^J: ~ \prod_{j=1}^J W_k(v_j)$ $\in U_q\} = (\alpha_k(q) \phi(q))^J$, we have noted that
\begin{equation}\label{eq:convnsplitformfnq=1}
\sum_{\substack{n \le x\conv\\\gcd(f(n), q)=1}} 1 ~ = \alpha_k(q)^J \largesum_{\substack{m \le x\\(f(m), q)=1}} \Bigg(\frac1{J!}\sum_{\substack{P_1, \dots, P_J > L_m\\P_1 \cdots P_J \le (x/m)^{1/k}\\ P_1, \dots, P_J \text{ distinct }}} 1\Bigg) ~ + ~ O\left(x^{1/k} \exp\left(-\frac{K_1}2 (\log{x})^{\epsilon/4}\right)\right).   
\end{equation}
For each $\wifam \in U_q^K$, we define $\mathcal U_{J, K}\big(q, Q_0; \wifam\big)$ to be the set of tuples $(v_1, \dots, v_J) \in U_q^J$ satisfying $\prod_{j=1}^J \Wik(v_j) \in U_q$ and $\prod_{j=1}^J \Wik(v_j) \equiv w_i \pmod{Q_0}$ for each $i \in [K]$. 
Observe that any convenient $n$ satisfying $\gcd(f(n), q)=1$ and $\finaimodQZ$ for all $i \in [K]$, can be uniquely written in the form $n = m(P_J \cdots P_1)^k$, where $P_J, \dots, P_1$ are primes satisfying \eqref{eq:convnDefnIneq}, $\gcd(f(m), q)=1$ and $(P_1, \dots, P_J)$ mod $q$ $\in \mathcal U_m \coloneqq \mathcal U_{J, K}\big(q, Q_0;$ $\aifimfam\big)$. As such, by the arguments leading to \eqref{eq:convnSplitForm}, we obtain
\begin{multline}\label{eq:convnqQ0SplitForm}
\largesum_{\substack{n \le x \conv\\\gcd(f(n), q)=1\\(\forall i) ~ \finaimodQZ}} 1 ~ ~ = \largesum_{\substack{m \le x\\(f(m), q)=1}} \frac{\#\mathcal U_m}{\phi(q)^J} \Bigg(\frac1{J!}\sum_{\substack{P_1, \dots, P_J > L_m\\P_1 \cdots P_J \le (x/m)^{1/k}\\ P_1, \dots, P_J \text{ distinct }}} 1\Bigg) + ~ \ophiqKsumnxfq.    
\end{multline}
Now, a simple counting argument shows the following general observation: let $F \in \Z[T]$ be a nonconstant polynomial, and let $Q, d$ be positive integers such that $d \mid Q$ and $\alpha_F(Q) \coloneqq \frac1{\phi(Q)} \#\{u \in U_Q: F(u) \in U_Q\}$ is nonzero (hence so is $\alpha_F(d)$). Then for any $u \in U_d$ for which $F(u) \in U_d$, we have
\begin{equation}\label{eq:LiftPolyCoprimality}
\#\{U \in U_Q: ~ U \equiv u \pmod d, ~ F(U) \in U_Q\} = \frac{\alpha_F(Q)\phi(Q)}{\alpha_F(d) \phi(d)}.     
\end{equation} 
Using this for $F \coloneqq W_k = \prodik \Wik$ (so that $\alpha_F = \alpha_k$), we immediately obtain
\begin{align*}\allowdisplaybreaks
\#\mathcal U_{J, K}\big(q, Q_0; \wifam\big) 
= \left(\frac{\alpha_k(q)\phi(q)}{\alpha_k(Q_0)\phi(Q_0)}\right)^J \#\mathcal V_{J, K}^{(k)}\big(Q_0, \wifam\big)    
\end{align*} 
for all $\wifam \in U_q^K$. Applying this with $w_i \coloneqq a_i f_i(m)^{-1}$ and recalling that $\VQZmPr$ $= \mathcal V_{J, K}^{(k)}\big(Q_0; \aifimfam\big)$, we get from \eqref{eq:convnqQ0SplitForm},
\begin{multline*}
\largesum_{\substack{n \le x \conv\\\gcd(f(n), q)=1\\(\forall i) ~ \finaimodQZ}} 1 ~ ~ = \frac{\alpha_k(q)^J}{\alpha_k(Q_0)^J}\largesum_{\substack{m \le x\\(f(m), q)=1}} \frac{\VQZmPr}{\phi(Q_0)^J} \Bigg(\frac1{J!}\sum_{\substack{P_1, \dots, P_J > L_m\\P_1 \cdots P_J \le (x/m)^{1/k}\\ P_1, \dots, P_J \text{ distinct }}} 1\Bigg)  + ~ \ophiqKsumnxfq.
\end{multline*}
Comparing this with \eqref{eq:ConvnSplitFormSimplified}, we obtain
$$\largesum_{\substack{n \le x\conv\\ \forallifinaimodq}} 1 ~ = (1+o(1)) \phiQZphiqK \largesum_{\substack{n \le x \conv\\\gcd(f(n), q)=1\\(\forall i) ~ \finaimodQZ}} 1 + \ophiqKsumnxfq.$$
Finally, an application of Proposition \ref{prop:inconvnfnq=1} allows us to remove the condition of $n$ being convenient from the main term on the right hand side above. This completes the proof of Proposition \ref{prop:ConvenientMainTermUptoQ0}, up to the proof of Proposition \ref{prop:Vqwi_ReducnToBddModulus}, which we take up in the next section. \end{proof} 
\section{Counting solutions to congruences: Generalization of Proposition \ref{prop:Vqwi_ReducnToBddModulus}}\label{sec:VNKtilCount_MegaGeneral_Sec}
We devote this section to establishing a general version of Proposition \ref{prop:Vqwi_ReducnToBddModulus} which shall also be useful while dealing with the contribution of the inconvenient $n$ in the proof of Theorems \ref{thm:UnrestrictedInput_SqfrPoly} to \ref{thm:HigherPolyControl1}. To do this, we shall primarily make use of two bounds on character sums, which we state in the next two propositions. In this section, we deviate from the notation and hypotheses set up in the introduction, assuming \textit{only} what is introduced in the rest of the section.
\begin{prop}\label{prop:WeilBounds} Let $\ell$ be a prime, $\chi$ a Dirichlet character mod $\ell$, and $F \in \Z[T]$ a nonconstant polynomial which is not congruent mod $\ell$ to a polynomial of the form $c \cdot G(T)^{\ord(\chi)}$ for some $c \in \F_\ell$ and $G \in \FellT$, where $\ord(\chi)$ denotes the order of the character $\chi$. Then $$\left|\largesum_{u \bmod \ell} \chi(F(u))\right| \le (d-1)\sqrt\ell,$$ where $d$ is the degree of the largest squarefree divisor of $F$.
\end{prop}
This is a version of the Weil bounds and is a special case of \cite[Corollary 2.3]{wan97} (see also \cite{davenport39}, \cite{weil48} and \cite{schmidt76} for older results). We will also need an analogue of the above result for character sums to higher prime power moduli, and this input is provided by the following consequences of Theorems 1.2 and 7.1 and eqn. (1.15) in work of Cochrane \cite{cochrane02} (see \cite{CLZ03} for related results). 

In what follows, for a polynomial $H \in \Z[T]$, we denote by $H'$ or $H'(T)$ the formal derivative of $H$. 
Given a prime $\ell$, by the \textsf{$\ell$-critical polynomial associated to $H$} we shall mean the polynomial $\critpoly_H \coloneqq \ell^{-\ord_\ell(H')} H'$, which has integer coefficients and can be considered as a nonzero element of the ring $\FellT$. Moreover, if $H$ is not identically zero in $\F_\ell[T]$ (i.e., if $\ord_\ell(H) = 0$), then by the \textsf{$\ell$-critical points} of $H$, we shall mean the set $\critpts(H; \ell) \subset \F_\ell$ of zeros of the polynomial $\critpoly_H$ which are not zeros of $H$ (both polynomials considered mod $\ell$). Finally, for any $\theta \in \F_\ell$, we use $\multtheta(H)$ to denote the multiplicity of $\theta$ as a zero of $H$.   
\begin{prop}\label{prop:Cochrane}
Let $\ell$ be a prime, $g \in \Z[T]$ a nonconstant polynomial, and $t \coloneqq \ord_\ell(g')$. Consider an integer $e \ge t+2$ and a primitive character $\chi$ mod $\ell^e$. Let $M \coloneqq \max_{\theta \in \mathcal A(g; \ell)} \multtheta(\critpoly_g)$ be the maximum multiplicity of an $\ell$-critical point.
\begin{enumerate}
\item[(i)] For odd $\ell$, we have
$\left|\largesum_{u \bmod \ell^e} \chi(g(u))\right| \le \left(\sum_{\theta \in \mathcal A(g; \ell)} \multtheta(\critpoly_g) \right) \ell^{t/(M + 1)} ~ \ell^{e(1-1/(M+1))}$.
\item[(ii)] For $\ell=2$ and $e \ge t+3$, we have $\left|\largesum_{u \bmod 2^e} \chi(g(u))\right| \le (12.5) 2^{t/(M + 1)} ~ 2^{e(1-1/(M+1))}.$
In fact, the sum is zero if $g$ has no $2$-critical points.
\end{enumerate}
\end{prop}
In order to make use of the aforementioned bounds, we will need to understand the quantities that appear when we apply them. The following observations enable us to do this. 
\begin{prop}\label{prop:OrdDerivInfo}
Let $\Fiset \subZT$ be nonconstant and multiplicatively independent. There exists a constant $C_1 \coloneqq C_1(\Filist) \in \NatNos$ such that all of the following hold:
\begin{enumerate}
\item[\textbf{(a)}] For any prime $\ell$, there are $O(1)$ many tuples $(A_1, \dots, A_K) \in [\ell-1]^K$ for which $F_1^{A_1} \cdots F_K^{A_K}$ 
is of the form $c \cdot G^{\ell-1}$ in $\FellT$ for some $c \in \F_\ell$ and $G \in \FellT$; here, the implied constant depends at most on $\Fiset$. In fact, if $\ell>C_1$ and $\gcd(\ell-1, \beta(\Filist)) = 1$, then the only such tuple is $(A_1, \dots, A_K) = (\ell-1, \dots, \ell-1)$.
\item[\textbf{(b)}] For any prime $\ell$ and any $(A_1, \dots, A_K) \in \NatNos^K$ satisfying $\ord_\ell(\prod_{i=1}^K F_i)=0$ and $(A_1, \dots$ $, A_K) \not\equiv (0, \dots, 0) \pmod\ell$, we have in the two cases below, 
\begin{multline}\label{eq:prop_OrdDerivInfo_tau(ell)Bound}
\tau(\ell) \coloneqq \ord_\ell\left((T^{\phi(\ell^r)} F_1(T)^{A_1} \cdots F_K(T)^{A_K})'\right) = \ord_\ell(\FtilT)\\
\begin{cases}
    = 0, &\text {if }\ell>C_1, r \ge 2\\
    \le C_1, &\text {if } \ell \le C_1, \ord_\ell\left(\prodik F_i \right)= 0, r \ge C_1+2,
\end{cases}    
\end{multline}
where $\FtilT \coloneqq \sumik A_i \FiTderivFjT$. In either of the two cases above, any root $\theta \in \F_\ell$ of the polynomial $\mathcal C_\ell(T) \coloneqq \ell^{-\tau(\ell)} (T^{\phi(\ell^r)} F_1(T)^{A_1} \cdots F_K(T)^{A_K})'$ which is not a root of $T\prodik F_i(T)$, must be a root of the polynomial $\ell^{-\tau(\ell)}\FtilT$ of the same multiplicity.\footnote{Once again, the last three polynomials are being considered as nonzero elements of $\FellT$.} 
\end{enumerate} 
\end{prop}
\begin{proof}
We start by writing $F_i \eqqcolon r_i \prod_{j=1}^M G_j^{\mu_{ij}}$ as in the introduction, so that $r_i \in \Z$ and $G_1, \dots, G_M \in \Z[T]$ are irreducible, primitive and pairwise coprime, and $M= \omega(\Fiprod)$. 
Recall that $M \ge K$ and that the exponent matrix $E_0(\Filist)$ has $\Q$-linearly independent columns, making $\beta(\Filist)$ a nonzero integer. Further, since $G_j$ are pairwise coprime irreducibles,  the resultants $\Res(G_j, G_{j'})$ 
and discriminants $\disc(G_j)$ are nonzero integers for all $j \ne j' \in [M]$. Note that for any prime $\ell$ not dividing the leading coefficient of any $G_j$ and not dividing $\prod_{1 \le j \le M} \disc(G_j) \cdot \prod_{1 \le j \ne j' \le M} \Res(G_j, G_{j'})$, the  
product $\prod_{j=1}^M G_j$ is separable in $\FellT$. 

We also observe that since $(F_1^{c_1} \cdots F_K^{c_K})' = \left(\prod_{i=1}^K F_i^{c_i-1}\right) \sumik c_i \FiderivFj$, the multiplicative independence of the polynomials $\Fiset$ forces the polynomials $\Fiderivset \subZT$ to be $\Q$-linearly independent. 
Writing $D \coloneqq \deg(\Fiprod)$ and $\FiTderivFjT \eqqcolon$ $\sum_{j=0}^{D-1} c_{i, j} T^j$ for some $\{c_{i, j}\}_{0 \le j \le D-1} \subset \Z$, we find that the columns of the matrix 
\begin{equation}\label{eq:M0Def}
M_1 \coloneqq M_1(\Filist) \coloneqq
\begin{pmatrix}
c_{1, 0} & \cdots & c_{K, 0}\\
\cdots & \cdots & \cdots\\
\cdots & \cdots & \cdots\\
c_{1, D-1} & \cdots & c_{K, D-1}
\end{pmatrix} \in \mathbb M_{D\times K}(\Z)
\end{equation}
must be $\Q$-linearly independent. Consequently, 
the last diagonal entry $\betatil \coloneqq \betatil(\Filist) \in \Z\sm\{0\}$ 
is the largest invariant factor of $M_1$ (in size). 

We now let $C_1 \coloneqq C_1(\Filist)$ be any positive integer exceeding $\max\{|\betatil|, 2\}$ such that for any $\ell>C_1$, $\ell$ divides neither the product $\prod_{j=1}^M \disc(G_j) \cdot \prod_{\substack{1 \le j \ne j' \le M}} \Res(G_j, G_{j'}) \in \Z\sm\{0\}$ nor the leading coefficient of any of $F_1, \dots, F_K$  (hence also none of the leading coefficients of $G_1, \dots, G_M$), and we have $\ord_\ell(\Fiprod) = 0$.
We claim that any such $C_1$ satisfies the properties in the statement of the proposition. 

\textit{Proof of (a).} We may assume that $\ell>C_1$. 
Let $\beta \coloneqq \beta(\Filist)$. As mentioned before, the conditions defining $C_1$ force $G_1, \dots, G_M$ to be pairwise coprime in $\FellT$. Let $(A_1, \dots, A_K) \ne (0, \dots, 0)$ be any tuple of nonnegative integers for which $F_1^{A_1} \cdots F_K^{A_K}$ is of the form $c \cdot G^{\ell-1}$ in $\FellT$ for some $c \in \F_\ell$ and $G \in \FellT$. We claim that $A_1, \dots, A_K$ must all be divisible by $(\ell-1)/d_1$ where $d_1 \coloneqq \gcd(\ell-1, \beta)$. This will be enough to complete the proof of (a), since there are no more than $d_1^K \le |\beta|^K \ll 1$ many tuples $(A_1, \dots, A_K) \in [\ell-1]^K$ satisfying this latter property. 

To establish the above claim, 
we may assume without loss of generality that $G$ is monic, and note that $c \in \F_\ell^\times$ since $\ord_\ell(\Fiprod) = 0$ by definition of $C_1$. Write each $G_j$ as $\lambda_j H_j$ in the ring $\FellT$, for some $\lambda_j \in \FellUnits$ and nonconstant monic $H_j \in \FellT$ (which can be done since $\ell$ doesn't divide the leading coefficient of any $G_j$). Then $F_i = r_i \prod_{j=1}^M G_j^{\mu_{ij}} = \rho_i \prod_{j=1}^M H_j^{\mu_{ij}}$ for some $\rho_i \in \FellUnits$. Since $c \cdot G^{\ell-1} = \prodik F_i^{A_i} = \left(\prodik \rho_i^{A_i}\right) \cdot \prod_{1 \le j \le M} H_j^{\sumik \mu_{ij} A_i}$
in $\FellT$, and $G, H_1, \dots, H_M$ are all monic, we find that $G^{\ell-1} = \prod_{1 \le j \le M} H_j^{\sumik \mu_{ij} A_i}$. But now since $\prod_{1 \le j \le M} G_j$ is separable in $\FellT$, so is $\prod_{1 \le j \le M} H_j$, 
and we deduce that 
$\sumik \mu_{ij} A_i \equiv 0 \pmod{\ell-1}$ for each $j \in [M]$. This can be rewritten as the matrix congruence $(0 \cdots \cdots 0)^\top \equiv E_0 (A_1 \cdots A_K)^\top \pmod{\ell-1}$, where each side is an $M \times 1$ matrix and $Y^\top$ denotes the transpose of a matrix $Y$. 

Now since $M \ge K$ and $E_0$ has full rank, there exist $P_0 \in GL_{M \times M}(\Z)$ and $R_0 \in GL_{K \times K}(\Z)$ for which $P_0 E_0 R_0$ is the Smith Normal Form $\diag(\beta_1, \dots, \beta_K)$ of $E_0$, with $\beta_1, \dots, \beta_K \in \Z\sm\{0\}$ being the invariant factors of $E_0$, so that $\beta_i \mid \beta_{i+1}$ for all $1 \le i < K$ and $\beta=\beta(\Filist) = \beta_K$. Thus $P_0 E_0 = \diag(\beta_1, \dots, \beta_K) R_0^{-1}$ and writing $(q_{ij})_{1 \le i, j \le K} \coloneqq R_0^{-1}$, we find that 
\begin{align*}\allowdisplaybreaks
\begin{pmatrix}
0\\ \cdots \\ \cdots \\ 0   
\end{pmatrix}_{M \times 1} \equiv P_0 E_0 \begin{pmatrix}
A_1\\ \cdots \\ A_K \end{pmatrix}_{K \times 1} 
\equiv \begin{pmatrix}\beta_1(q_{11} A_1 + \cdots + q_{1K}A_K)\\ \cdots \\ \beta_K(q_{K1} A_1 + \cdots + q_{KK}A_K)\\ 0 \\ \cdots\\ 0\end{pmatrix}_{M \times 1} \pmod{\ell-1}.
\end{align*}
Hence for each $i \in [K]$, $\beta_i(q_{i1} A_1 + \cdots + q_{iK} A_K) \equiv 0 \pmod{\ell-1}$, so that $(\ell-1)/\gcd(\ell-1, \beta_i)$ divides $q_{i1} A_1 + \cdots + q_{iK} A_K$. But since $\beta_i \mid \beta_K$, it follows that $(\ell-1)/\gcd(\ell-1, \beta_K) = (\ell-1)/d_1$ also divides $q_{i1} A_1 + \cdots + q_{iK} A_K$ for each $i \in [K]$. We obtain 
\begin{equation}\label{eq:matrixcomputn}
\begin{pmatrix}
0\\ \cdots \\ 0   
\end{pmatrix}_{K \times 1} \equiv \begin{pmatrix}q_{11} A_1 + \cdots + q_{1K}A_K\\ \cdots \\ q_{K1} A_1 + \cdots + q_{KK}A_K\end{pmatrix}_{K \times 1} \equiv R_0^{-1} \begin{pmatrix}
A_1\\ \cdots \\ A_K \end{pmatrix}_{K \times 1} \left(\bmod ~ \frac{\ell-1}{d_1}\right),
\end{equation}
establishing the desired claim that $(A_1, \dots, A_K) \equiv (0, \dots, 0) \left(\bmod ~ \frac{\ell-1}{d_1}\right)$.

\textit{Proof of (b).} 
We start by noting that 
\begin{equation}\label{eq:TF1...FKDeriv}
(T^{\phi(\ell^r)} F_1(T)^{A_1} \cdots F_K(T)^{A_K})' = \phi(\ell^r) T^{\phi(\ell^r)-1} \prodik F_i(T)^{A_i} + T^{\phi(\ell^r)} \left(\prodik F_i(T)^{A_i-1} \right) \FtilT,
\end{equation}
where $\FtilT$ is as in the statement of the proposition. We claim that $\ord_\ell(\Ftil) \le \bbm_{\ell \le C_1} C_1$ for all primes $\ell$ satisfying $\ord_\ell(\Fiprod)=0$ and for all nonnegative integers $A_1, \dots, A_K$ satisfying $(A_1, \dots, A_K) \not\equiv (0, \dots, 0) \bmod \ell$. To show this, we proceed as in the proof of (a), but working with the matrix $M_1$ defined in \eqref{eq:M0Def} in place of the exponent matrix $E_0$. Observe that $\FtilT = \sum_{j=0}^{D-1} \left(\sumik c_{i, j} A_i\right) T^j$, hence if $\kappaell \coloneqq \ord_\ell(\Ftil)$, then $\ell^\kappaell$ divides all the entries of the matrix $M_1 (A_1 \cdots A_K)^\top$. Since $M_1$ has full rank and $D = \sumik \deg F_i \ge K$ many rows, and since $(A_1, \dots, A_K) \not\equiv (0, \dots, 0) \bmod \ell$, an argument entirely analogous to the one leading to \eqref{eq:matrixcomputn} shows that $\ell^{\kappa(\ell)}$ divides the last invariant factor $\betatil$ of $M_1$. Hence $\ord_\ell(\Ftil) = \kappaell \le v_\ell(\betatil)$ and our claim follows as $|\betatil|<C_1$.

As a consequence, we find that $\ord_\ell\left(T^{\phi(\ell^r)} \left(\prodik F_i(T)^{A_i-1} \right) \FtilT\right) = \ord_\ell(\Ftil) \le \bbm_{\ell \le C_1} C_1$ for all primes $\ell \le C_1$ satisfying $\ord_\ell(\Fiprod) = 0$, and also for all primes $\ell>C_1$ (for which the condition $\ord_\ell(\Fiprod) = 0$ is automatic by definition of $C_1$). But now since $\ord_\ell(\phi(\ell^r)) \ge 1$ for $r \ge 2$ and $\ord_\ell(\phi(\ell^r)) \ge C_1+1$ for $r \ge C_1+2$, \eqref{eq:TF1...FKDeriv} shows that $\tau(\ell) = \ord_\ell\left(T^{\phi(\ell^r)} \left(\prodik F_i(T)^{A_i-1} \right) \FtilT\right)$, 
establishing subpart (b) of the proposition. 

Finally, since in both the cases of \eqref{eq:prop_OrdDerivInfo_tau(ell)Bound}, we have $\tau(\ell) < r-1$, the identity \eqref{eq:TF1...FKDeriv} reveals that 
$$\mathcal C_\ell(T) \equiv \ell^{-\tau(\ell)} \left(T^{\phi(\ell^r)} \prodik F_i(T)^{A_i}\right)' \equiv T^{\phi(\ell^r)} \left(\prodik F_i(T)^{A_i-1}\right) \left(\ell^{-\tau(\ell)} \FtilT \right) \text{  in the ring }\FellT.$$
As such, any root of the polynomial $\theta \in \F_\ell$ of $\mathcal C_\ell(T)$ (considered as a nonzero element of $\FellT$) which is not a root of $T\prodik F_i(T)$, must be a root of $\ell^{-\tau(\ell)} \FtilT$, and $\theta$ must have the same multiplicity in $\mathcal C_\ell(T)$ and $\ell^{-\tau(\ell)}\FtilT$. This completes the proof of Proposition \ref{prop:OrdDerivInfo}. 
\end{proof}
We now come to the main result of this section: the promised generalization of Proposition \ref{prop:Vqwi_ReducnToBddModulus}. The following notation and conventions will be relevant only in the rest of the section. 

Let $\GirsetOKL \subZT$ be a fixed collection of nonconstant polynomials such that for each $r \in [L]$, the polynomials $\GirsetiOK \subZT$ are multiplicatively independent. Define $D_0 \coloneqq \max_{1 \le r \le L} ~ \sumik \deg \Gir$. Let $N \ge 1$ and $\FijsetOKN \subZT$ be a family of polynomials such that for each $j \in [N]$, the vector $\FijVectiOK$ coincides with one of the vectors $\GijPrVectiOK$ for some $j' \in [L]$ (possibly depending on $j$). In this case we define, for any integer $q$,
$$\alphajtilq \coloneqq \frac1{\phi(q)} \#\{u \in U_q: \prodik \Fij(u) \in U_q\} = \frac1{\phi(q)} \#\{u \in U_q: \prodik G_{i, j'}(u) \in U_q\},$$
and let $\alphaStNq \coloneqq \prod_{j=1}^N \alphajtilq$. For any $\wifam \in U_q^K$, define
$$\VNKtilqwi \coloneqq \left\{(v_1, \dots, v_N) \in U_q^N: ~ (\forall i \in [K]) ~ \prod_{j=1}^N \Fij(v_j) \equiv w_i \pmod q\right\}.$$
Fix $B_0>0$. In the next three results, the implied constants may depend only on $B_0$ and on the fixed collection of polynomials $\GirsetOKL$ (besides other parameters declared explicitly).
\begin{prop}\label{prop:VNKCountmegageneral}
There exists a constant $C_0 \coloneqq C_0\left(\GirsetOKL; B_0\right)> (8D_0)^{2D_0+2}$ depending only on $\GirsetOKL$ and $B_0$, such that for \underline{any} constant $C>C_0$, the following hold.
\begin{enumerate}
\item[\textbf{(a)}] Uniformly for $N \ge KD_0 + 1$ and in coprime residues $\wilist$ to moduli $q$ satisfying $\alphaStNq \ne 0$ and $\IFHGirBZ$ for each $r \in [L]$, we have\begin{multline}\label{eq:VNKkCount_largeN_Gen} \allowdisplaybreaks
\frac{\#\VNKtilqwi}{\phi(q)^N}\\ = \alphaNStqRatio \phiQZphiqK \left\{\frac{\#\VNKtilQZwi}{\phi(Q_0)^N} + O\left(\frac1{C^N}\right)\right\} \mathlarger{\prod}_{\substack{\ell \mid q\\\ell>C_0}} \left(1+O\left(\frac{(4D_0)^N}{\ell^{N/D_0 -K}}\right)\right), 
\end{multline} 
where $Q_0$ is a $C_0$-smooth divisor of $q$ of size $O_C(1)$. 
\item[\textbf{(b)}] 
\text{For any \textbf{fixed} }$N \ge 1$ and uniformly in coprime residues $\wilist$ mod $q$, we have
\begin{align}\label{eq:VNKkCount_smallN_Gen} \allowdisplaybreaks
\frac{\#\VNKtilqwi}{\phi(q)^N} \le \frac{\big(\prod_{\ell^e \parallel q} e\big)^{\bbm_{N=KD_0}}}{q^{\min\{K, N/D_0\}}} ~ \expOomegaq. 
\end{align}
\end{enumerate}
\end{prop}
\begin{proof}
In what follows, $q$ is an arbitrary positive integer (unless stated otherwise). We may assume that $\alphaStNq \ne 0$, for both the assertions of the proposition are vacuous or tautological otherwise. In particular, this means that $\ord_\ell(\prodik \prod_{j=1}^N F_{i, j}) = 0$ for each prime $\ell \mid q$. Fix $C_0 \coloneqq C_0\left(\GirsetOKL; B_0\right)$ to be any constant exceeding $B_0$, $(32 D_0)^{2D_0+2}$, the sizes of the leading and constant coefficients of $\GirsetOKL$, as well as the constants $C_1(\GirlistiOK)$ coming from applications of Proposition \ref{prop:OrdDerivInfo} to each of the families $\GirsetiOK$ of multiplicatively independent polynomials. We will show that any such choice of $C_0$ suffices. 

We first consider the case $D_0>1$; we will deal with the $D_0=1$ case towards the end of this argument. For an arbitrary positive integer $Q$ and coprime residues $\wilist$ mod $Q$, we apply 
the orthogonality of Dirichlet characters to detect the congruences defining $\VNKtilQwi$. This yields
\begin{align}\label{eq:OrthogGen}\allowdisplaybreaks
\#\VNKtilQwi 
&= \frac1{\phi(Q)^K} \largesum_{\chi_1, \dots, \chi_K \bmod Q} \chiibarwiprod \prodjN Z_{Q; ~ \chi_1, \dots, \chi_K} (\FijlistiOK),
\end{align}
where $Z_{Q; ~ \chi_1, \dots, \chi_K}(\FijlistiOK) \coloneqq \largesum_{v \bmod Q} \chi_{0, Q}(v) \prodik \chi_i(\Fij(v))$ and $\chi_{0, Q}$ denotes (as usual) the trivial character mod $Q$. 

We show the following estimates, both uniform in residues $\wilist \in U_{\ell^e}$ for primes $\ell>C_0$:   
\begin{enumerate} 
\item[\textbf{(i)}] If $\alpha_N^*(\ell) \ne 0$ and $\gcd(\ell-1, \beta(\GirlistiOK))=1\text{ for each }r\in [L]$, 
then
\begin{equation}\label{eq:VNKtil_LargePrimes_LargeN}
\frac{\#\widetilde{\mathcal V}_{N, K}(\ell^e; \wifam)}{\phi(\ell^e)^N} 
= \displaystyle{\frac{\alpha_N^*(\ell)}{\phi(\ell^e)^K} \left(1 + O\left(\frac{(4D_0)^N}{\ell^{N/D_0-K}}\right)\right)},
\end{equation}
uniformly in $N \ge KD_0+1$.
\item[\textbf{(ii)}] For each fixed $N \ge 1$, there is a constant $K'$ depending at most on $N$ and $\GirsetOKL$ such that 
\begin{equation}\label{eq:VNKtil_LargePrimes_SmallN}
\frac{\#\widetilde{\mathcal V}_{N, K}(\ell^e; \wifam)}{\phi(\ell^e)^N} \le \displaystyle{K' ~ ~\frac{e^{\bbm_{N = KD_0}}}{(\ell^e)^{\min\{K, N/D_0\}}}}.
\end{equation}
\end{enumerate}
To show these, we start by applying \eqref{eq:OrthogGen} with $Q \coloneqq \ell^e$ to get
\begin{multline}\label{eq:OrthogPrimePower_Bound}
\frac{\#\widetilde{\mathcal V}_{N, K}(\ell^e; \wifam)}{\phi(\ell^e)^N}\\ \le \frac1{\phi(\ell^e)^K} \Bigg\{1+ \frac1{\phi(\ell^e)^N} \largesum_{\chiituplist \ne \chiZelltuplist \bmod \ell^e} ~ \prod_{j=1}^N |Z_{\ell^e; ~ \chi_1, \dots, \chi_K} (\FijlistiOK)| \Bigg\};      
\end{multline}
in addition, if $\alpha_N^*(\ell) \ne 0$, then from $Z_{\ell^e; ~ \chi_{0, \ell}, \dots, \chi_{0, \ell}} (F_{1, j}, \dots, $ $F_{K, j}) =$ $= \alphajtil(\ell) \phi(\ell^e)$, we have
\begin{multline}\label{eq:OrthogPrimePower}
\frac{\#\widetilde{\mathcal V}_{N, K}(\ell^e; \wifam)}{\phi(\ell^e)^N} = \frac{\alpha_N^*(\ell)}{\phi(\ell^e)^K} \Bigg\{1+ \\\frac1{\alpha_N^*(\ell) \phi(\ell^e)^N} \largesum_{\chiituplist \ne \chiZelltuplist \bmod \ell^e} ~ \chiibarwiprod \prod_{j=1}^N Z_{\ell^e; ~ \chi_1, \dots, \chi_K} (\FijlistiOK) \Bigg\}.    
\end{multline}
Now consider any tuple $\chiituplist \ne \chiZelltuplist$ mod $\ell^e$ and any $j \in [N]$. Let $\ell^{e_0} \coloneqq \lcm[\cond(\chi_1), \dots, \cond(\chi_K)]$ $\in \{\ell, \dots, \ell^e\}$. Using $\chi_1, \dots, \chi_K$ to also denote the characters mod $\ell^{e_0}$ inducing $\chi_1, \dots, \chi_K$ respectively, we get 
\begin{equation}\label{eq:ReducntoPrimIdentity}
Z_{\ell^e; ~ \chi_1, \dots, \chi_K} (\FijlistiOK) = \ell^{e-e_0} ~ Z_{\ell^{e_0}; ~ \chi_1, \dots, \chi_K} (\FijlistiOK)
\end{equation} 
Since $\ell>C_0>2$, 
the character group mod $\ell^{e_0}$ is generated by the character $\psieZ$ given by $\psieZ(\gamma) \coloneqq \exp(2 \pi i/\phi(\ell^{e_0}))$, 
for some generator $\gamma$ of $U_{\ell^{e_0}}$. As such, there exists a tuple $(A_1, \dots, A_K) \in [\phi(\ell^{e_0})]$ satisfying $\chi_i = \psieZ^{A_i}$ for each $i$, and 
\begin{equation}\label{eq:AiTupRestric}
(A_1, \dots, A_K) \not\equiv  
\begin{cases}
(0, \dots, 0) \pmod{\ell}, &\text{if } e_0>1,\\
(0, \dots, 0) \pmod{\ell-1}, &\text{if } e_0=1,
\end{cases}    
\end{equation} 
since at least one of $\chi_1, \dots, \chi_K$ is primitive mod $\ell^{e_0}$. This gives
\begin{equation}\label{eq:Reducn to Single Charac}
Z_{\ell^{e_0}; ~ \chi_1, \dots, \chi_K} (\FijlistiOK) = \largesum_{v \bmod \ell^{e_0}} \psieZ \left(v^{\phi(\ell^{e_0})} \prodik F_{i, j}(v)^{A_i}\right).    
\end{equation}
We now consider two possibilities, namely when $e_0 = 1$ or $e_0 \ge 2$. 

\textit{Case $1$:} Suppose $e_0 = 1$. For each $j \in [N]$, consider $j' \in [L]$ satisfying $(G_{i, j'})_{i=1}^K = (\Fij)_{i=1}^K$. By Proposition \ref{prop:OrdDerivInfo}(a), we see 
there are $O_L(1)$ 
many possible tuples $\chiituplist$ of characters mod $\ell^e$ having $\lcmcondchilist = \ell$, 
for which $T^{\phi(\ell)}\prodik F_{i, j}(T)^{A_i} = T^{\phi(\ell)}\prodik G_{i, j'}(T)^{A_i}$ is of the form $c \cdot G(T)^{\ell-1}$ in $\FellT$ for some $j \in [N]$ (here $A_i$ are as above). Moreover if $\gcd(\ell-1, \beta(\GirlistiOK))=1$ for all $r \in [L]$, then there is no such tuple $\chiituplist$. For all the remaining tuples $\chiituplist$ with $\lcmcondchilist = \ell$, we may invoke Proposition \ref{prop:WeilBounds} to obtain, \textbf{for all $j \in [N]$}, 
\begin{multline*}
|Z_{\ell; ~ \chi_1, \dots, \chi_K} (\FijlistiOK)| = \Bigg|\largesum_{v \bmod \ell} \psi_1 \left(v^{\phi(\ell)} \prodik F_{i, j}(v)^{A_i}\right)\Bigg| \le \left(\sumik \deg F_{i, j}\right) \ell^{1/2} 
\le D_0 \ell^{1/2}.
\end{multline*}
By \eqref{eq:ReducntoPrimIdentity}, we deduce that for all but $O_L(1)$ many tuples $(\chi_1, \dots, \chi_K)$ of characters mod $\ell^e$ satisfying $\lcmcondchilist = \ell$, we have
\begin{equation}\label{eq:Zell^eBoundFor e_0=1} 
|Z_{\ell^e; ~ \chi_1, \dots, \chi_K} (\FijlistiOK)|  \le D_0 \ell^{e-1/2} ~ \text{ for every }j \in [N],
\end{equation} 
and when $\gcd(\ell-1, \beta(\GirlistiOK))=1$ for all $r \in [L]$, this inequality is true for all $(\chi_1, \dots, \chi_K)$ mod $\ell^e$ satisfying $\lcmcondchilist = \ell$. 

\textit{Case 2:} Now assume that $e_0 \ge 2$. Consider an arbitrary $j \in [N]$ and let $\GijPrVectiOK = (\Fij)_{i=1}^K$ for some $j' \in [L]$. Since $\ell>C_0 > C_1(\GijPrlistiOK)$ and $e_0 \ge 2$, Proposition \ref{prop:OrdDerivInfo}(b) and \eqref{eq:AiTupRestric} show that $\tau(\ell) \coloneqq \ord_\ell \left((T^{\phi(\ell^{e_0})} \prodik F_{i, j}(T)^{A_i})'\right)$ $= 0$. Consequently, \eqref{eq:Reducn to Single Charac} and Proposition \ref{prop:Cochrane}(i) yield
$|Z_{\ell^{e_0}; ~ \chi_1, \dots, \chi_K} (\FijlistiOK)| \le \left(\sum_{\theta \in \mathcal A_\ell} \mu_\theta(\mathcal C_\ell) \right) \ell^{e_0(1-1/(M_\ell+1))},$
where $\mathcal A_\ell \subset \F_\ell$ denotes the set of 
$\ell$-critical points of the polynomial $T^{\phi(\ell^{e_0})} \prodik F_{i, j}(T)^{A_i}$, 
$\mathcal C_\ell(T) \coloneqq (T^{\phi(\ell^{e_0})} \prodik F_{i, j}(T)^{A_i})'$ and 
$M_\ell \coloneqq \max_{\theta \in \mathcal A_\ell} \mu_\theta(\mathcal C_\ell)$. Moreover, by the last assertion in Propositon \ref{prop:OrdDerivInfo}, any $\theta \in \mathcal A_\ell$ is a root of the polynomial $\FtilT \coloneqq \sumik A_i F_{i, j}'(T) \prod_{\substack{1 \le r \le K\\ r \ne i}} F_{r, j}(T)$ (a nonzero element of $\FellT$), and $\mu_\theta(\critpoly_\ell) = \mu_\theta(\Ftil)$. As such, $M_\ell \le \sum_{\theta \in \mathcal A_\ell} \mu_\theta(\critpoly_\ell) \le D_0-1$, yielding $|Z_{\ell^{e_0}; ~ \chi_1, \dots, \chi_K} (\FijlistiOK)| \le (D_0-1) \ell^{e_0(1-1/D_0)}$. Thus, by \eqref{eq:ReducntoPrimIdentity}, 
\begin{equation}\label{eq:Zell^eBoundFor e_0>1}
|Z_{\ell^e; ~ \chi_1, \dots, \chi_K} (\FijlistiOK)| \le (D_0-1) \ell^{e-e_0/D_0}  ~ ~ \text {if }~ ~ \ell^{e_0} \coloneqq \lcmcondchilist \in \{\ell^2, \dots, \ell^e\}.
\end{equation} 
For any $e_0 \in \{1, \dots, e\}$ there are at most $\ell^{e_0K}$ tuples $(\chi_1, \dots, \chi_K)$ of characters mod $\ell^e$ having $\lcmcondchilist 
 = \ell^{e_0}$. Combining \eqref{eq:Zell^eBoundFor e_0>1} with the respective assertion in \eqref{eq:Zell^eBoundFor e_0=1}, 
 we get
 \begin{multline}\label{eq:SimplificnStep1}
\largesum_{\chiituplist \ne \chiZelltuplist \bmod \ell^e} ~ \left|\prod_{j=1}^N Z_{\ell^e; ~ \chi_1, \dots, \chi_K} (\FijlistiOK)\right| 
\le D_0^N \ell^{eN} \largesum_{1 \le e_0 \le e} \ell^{e_0(K-N/D_0)},
\end{multline}
for any prime power $\ell^e$ with $\ell>C_0$ satisfying $\alpha_N^*(\ell) \ne 0$ and $\gcd(\ell-1, \beta(\GirlistiOK))=1$ for each $r \in [L]$. (In the last inequality above, we have used the fact that $D_0 \ge 2$.) Now for each $j \in [N]$, $\alphajtil(\ell) \ge 1-D_0/(\ell-1) > 1-D_0/(C_0-1) > 1/2$, so that $\alpha_N^*(\ell) \ge 1/2^N$. If $N \ge KD_0+1$, then $\ell^{K-N/D_0} \le \ell^{-1/D_0} \le C_0^{-1/D_0} \le 1/2$, and \eqref{eq:SimplificnStep1} yields 
\begin{multline}\label{eq:CharacSum_LargePrime_Explicit}
\frac1{\alpha_N^*(\ell) \phi(\ell^e)^N}\largesum_{\chiituplist \ne \chiZelltuplist \bmod \ell^e} ~ \left|\prod_{j=1}^N Z_{\ell^e; ~ \chi_1, \dots, \chi_K} (\FijlistiOK)\right|
\le \frac{2 (4D_0)^N}{\ell^{N/D_0-K}}.
\end{multline} 
Inserting this bound into \eqref{eq:OrthogPrimePower} shows the assertion \eqref{eq:VNKtil_LargePrimes_LargeN}.  
On the other hand for \textit{any} prime power $\ell^e$ with $\ell>C_0$, \eqref{eq:Zell^eBoundFor e_0=1} and \eqref{eq:Zell^eBoundFor e_0>1} show that for each fixed $N \ge 1$,
\begin{align}\label{eq:FixedN_LargePrime_CharacSum_Intermed} \allowdisplaybreaks
\largesum_{\chiituplist \ne \chiZelltuplist \bmod \ell^e} ~ \left|\prod_{j=1}^N Z_{\ell^e; ~ \chi_1, \dots, \chi_K} (\FijlistiOK)\right| 
\ll ~ \phi(\ell^e)^N + D_0^N \ell^{eN} \largesum_{1 \le e_0 \le e} \ell^{e_0(K-N/D_0)}. 
\end{align}
Thus for a fixed $N \ge KD_0+1$, a calculation analogous to \eqref{eq:CharacSum_LargePrime_Explicit} yields 
$$\frac1{\phi(\ell^e)^N}\largesum_{\chiituplist \ne \chiZelltuplist \bmod \ell^e} ~ \left|\prod_{j=1}^N Z_{\ell^e; ~ \chi_1, \dots, \chi_K} (\FijlistiOK)\right| \ll ~ 1.$$
On the other hand if $N \in \{1, \dots, KD_0\}$, the expression in \eqref{eq:FixedN_LargePrime_CharacSum_Intermed} leads to
\begin{multline*}
\frac1{\phi(\ell^e)^N}\largesum_{\chiituplist \ne \chiZelltuplist \bmod \ell^e} ~ \left|\prod_{j=1}^N Z_{\ell^e; ~ \chi_1, \dots, \chi_K} (\FijlistiOK)\right|  
\ll e^{\bbm_{N=KD_0}} \ell^{e(K-N/D_0)}.
\end{multline*}
Inserting the last two bounds displays into \eqref{eq:OrthogPrimePower_Bound} 
yields \eqref{eq:VNKtil_LargePrimes_SmallN}. 

Now for an arbitrary $q$, 
we let $\qtil \coloneqq \prod_{\substack{\ell^e \parallel q\\ \ell \le C_0}} \ell^e$ denote the $C_0$-smooth part of $q$. By \eqref{eq:OrthogGen}, 
\begin{equation}\label{eq:VNKqtil_Orthog}
\#\VNKtilqtilwi = \frac1{\phi(\qtil)^K} \largesumchiiOKmodqtil \chiibarwiprod \prodjN Z_{\qtil; ~ \chi_1, \dots, \chi_K} (\FijlistiOK).
\end{equation} 
Given a constant $C>C_0$, we fix $\kappa$ to be any integer constant exceeding $C \cdot (30 D_0 C_0^{C_0})^{2C_0}$, and let $Q_0 \coloneqq \prod_{\ell^e \parallel \qtil} ~ \ell^{\min\{e, \kappa\}} = \prod_{\ell \le C_0} \ell^{\min\{v_\ell(q), \kappa\}}$ denote the largest $(\kappa+1)$-free divisor of $\qtil$. Write the expression on the right hand side of \eqref{eq:VNKqtil_Orthog} as $\SPr + \SDPr$, where $\SPr$ denotes the contribution of those tuples $\chiituplist$ mod $\qtil$ for which $\lcmcondchilist$ is $(\kappa+1)$-free, or equivalently, those $\chiituplist$ for which $\lcmcondchilist$ divides $Q_0$.

For each tuple $\chiituplist$ counted in $\SPr$, there exists a unique tuple $(\psi_1, \dots, \psi_K)$ of characters mod $Q_0$ inducing $\chiituplist$ mod $\qtil$, respectively. Noting that $\alphajtil(\qtil) = \alphajtil(Q_0)$, a straightforward calculation using \eqref{eq:LiftPolyCoprimality} shows that   
\begin{multline*}
Z_{\qtil; ~ \chi_1, \dots, \chi_K} (\FijlistiOK) = \largesum_{v \bmod \qtil} \chi_{0, \qtil} ~(v) \prodik \chi_i(\Fij(v)) 
= \frac{\phi(\qtil)}{\phi(Q_0)} Z_{Q_0; ~ \psi_1, \dots, \psi_K} (\FijlistiOK)
\end{multline*}
for each $j \in [N]$. Consequently,  
\begin{align*}\allowdisplaybreaks
\SPr 
= \frac1{\phi(\qtil)^K} \left(\frac{\phi(\qtil)}{\phi(Q_0)}\right)^N \largesum_{\substack{\psi_1, \dots, \psi_K \bmod Q_0}} \overline\psi_1(w_1) \cdots \overline\psi_K(w_K) \prodjN Z_{Q_0; ~ \psi_1, \dots, \psi_K} (\FijlistiOK),
\end{align*}
so that invoking \eqref{eq:OrthogGen} with $Q \coloneqq Q_0$, we obtain
\begin{equation}\label{eq:SPrimeExpr}
\frac{\SPr}{\phi(\qtil)^N} = \left(\frac{\phi(Q_0)}{\phi(\qtil)}\right)^K \frac{\#\VNKtilQZwi}{\phi(Q_0)^N}.    
\end{equation}
We now deal with the remaining sum 
$$\SDPr = \frac1{\phi(\qtil)^K} \largesum_{\substack{\chi_1, \dots, \chi_K \bmod \qtil\\ \lcmcondchilist \text{ is not }(\kappa+1)\text{-free}}} \chiibarwiprod \prodjN Z_{\qtil; ~ \chi_1, \dots, \chi_K} (\FijlistiOK).$$
For each tuple $(\chi_1, \dots, \chi_K)$ of characters mod $\qtil$ considered in the sum above, we factor $\chi_i \eqqcolon \prod_{\ell^e \parallel \qtil} \chi_{i, \ell}$, where $\chi_{i, \ell}$ is a character mod $\ell^e$. With $e_\ell \coloneqq v_\ell\left(\lcmcondchilist\right)$, we observe that since $\cond(\chi_i) = \prod_{\ell^e \parallel q} \cond(\chi_{i, \ell})$ and each $\cond(\chi_{i, \ell})$ is a power of $\ell$, 
we must have $\lcm[\cond(\chi_{1, \ell}), \dots, \cond(\chi_{K, \ell})] = \ell^{e_\ell}$. Letting $\chi_{1, \ell}, \dots, \chi_{K, \ell}$ also denote the characters mod $\ell^{e_\ell}$ inducing $\chi_{1, \ell}, \dots, \chi_{K, \ell}$ mod $\ell^e$ respectively (for each $\ell^e \parallel \qtil$), we see that at least one of $\chi_{1, \ell}, \dots, \chi_{K, \ell}$ must be primitive mod $\ell^{e_\ell}$. Furthermore for each $j \in [N]$, 
$Z_{\qtil; ~ \chi_1, \dots, \chi_K} (\FijlistiOK)$ $= \prod_{\ell^e \parallel \qtil} Z_{\ell^e; ~ \chi_{1, \ell}, \dots, \chi_{K, \ell}} (\FijlistiOK),$
so that 
\begin{equation}\label{eq:S''PrepForCochrane}
|Z_{\qtil; ~ \chi_1, \dots, \chi_K} (\FijlistiOK)| \le \left(\prod_{\substack{\ell^e \parallel \qtil\\e_\ell \le \kappa}} \phi(\ell^e)\right) \prod_{\substack{\ell^e \parallel \qtil\\ e_\ell \ge \kappa+1}} \left(\ell^{e-e_\ell} |Z_{\ell^{e_\ell}; ~ \chi_{1, \ell}, \dots, \chi_{K, \ell}} (\FijlistiOK)|\right).
\end{equation}
We will show that 
prime powers for all $\ell^e \parallel \qtil$ with $e_\ell \ge \kappa+1$, we have
\begin{equation}\label{eq:Cochrane to Each Component}
|Z_{\ell^{e_\ell}; ~ \chi_{1, \ell}, \dots, \chi_{K, \ell}} (\FijlistiOK)| \le (D_0 C_0^{C_0}) ~ \ell^{e_\ell(1-1/D_0)}.   
\end{equation}
This follows for odd $\ell$, by essentially the same argument as that given for \eqref{eq:Zell^eBoundFor e_0>1}, the only difference is that this time we use \textit{both} the assertions in 
\eqref{eq:prop_OrdDerivInfo_tau(ell)Bound} since $e_\ell \ge \kappa+1 > (30 D_0 C_0)^{2 C_0} + 1 > C_0+2$. So assume that $\ell=2$, i.e. $e_2 = v_2(\lcmcondchilist) \ge \kappa+1 \ge 31$. 

We shall use Proposition \ref{prop:Cochrane}(ii) to bound the sum $Z_{2^{e_2}; ~ \chi_{1, 2}, \dots, \chi_{K, 2}} (\FijlistiOK)$. To do this, we observe that since $e_2 \ge 4$, the characters $\psi, \eta$ mod $2^{e_2}$ defined by 
$$\psi(5) \coloneqq \exp(2\pi i/2^{e_2-2}), ~ \psi(-1) \coloneqq 1 ~ ~ \text{  and  } ~ ~ \eta(5) \coloneqq 1, \eta(-1) \coloneqq -1$$
generate the character group mod $2^{e_2}$. Hence for each $i \in [K]$, there exist $A_{i, 2} \in [2^{e_2-2}]$ and $B_{i, 2} \in [2]$ satisfying $\chi_{i, 2} = \psi^{A_{i, 2}} \eta^{B_{i, 2}}$ and $(A_{1, 2}, \dots, A_{K, 2}) \not\equiv (0, \dots, 0) \bmod 2$ (since $e_2 \ge 4$ and at least one of $\chi_{1, 2}, \dots, \chi_{K, 2}$ is primitive mod $2^{e_2}$). This allows us to write 
\begin{equation*}
Z_{2^{e_2}; ~ \chi_{1, 2}, \dots, \chi_{K, 2}} (\FijlistiOK) = \largesum_{v \bmod 2^{e_2}} \psi\left(g(v)\right) ~ \eta\left(v^2 \prodik F_{i, j}(v)^{B_{i, 2}}\right),
\end{equation*}
where $g(T) \coloneqq \prodik F_{i, j}(T)^{A_{i, 2}}$. Now $\eta$ is induced by the nontrivial character mod $4$ and $v^2 \prodik F_{i, j}(v)^{B_{i, 2}} \equiv \prodik F_{i, j}(\iota)^{B_{i, 2}}$ if $v \equiv \iota \pmod 4$ ($\iota=\pm 1$).  
Thus, writing $v \coloneqq 4u+\iota$ gives
\begin{multline*}
Z_{2^{e_2}; ~ \chi_{1, 2}, \dots, \chi_{K, 2}} (\FijlistiOK)\\ = \eta\left(\prodik F_{i, j}(1)^{B_{i, 2}}\right)\largesum_{\substack{u \bmod 2^{e_2-2}}} \psi\left(h_1(u)\right) + \eta\left(\prodik F_{i, j}(-1)^{B_{i, 2}}\right)\largesum_{\substack{u \bmod 2^{e_2-2}}} \psi\left(h_{-1}(u)\right), 
\end{multline*} 
where $h_\iota(T) \coloneqq g(4T + \iota)$. Note that 
\begin{multline}\label{eq:Z2^e2 InTermsOf psi eta}
Z_{2^{e_2}; ~ \chi_{1, 2}, \dots, \chi_{K, 2}} (\FijlistiOK)\\ = \frac14 \cdot \eta\left(\prodik F_{i, j}(1)^{B_{i, 2}}\right)\largesum_{\substack{u \bmod 2^{e_2}}} \psi\left(h_1(u)\right) + \frac14 \cdot \eta\left(\prodik F_{i, j}(-1)^{B_{i, 2}}\right)\largesum_{\substack{u \bmod 2^{e_2}}} \psi\left(h_{-1}(u)\right).
\end{multline}
We will now show that the first of the two terms must have size no more than $(12.5) \cdot 2^{2D_0+C_0} \cdot 2^{e_2(1-1/D_0)}$; by an analogous argument, so does the second. If the first term is nonzero, then $\prodik F_{i, j}(1)^{B_{i, 2}} \equiv 1 \pmod 2$, so that $\ord_2\left(\prodik F_{i, j}(4T+1)^{A_{i, 2}-1}\right) = 0$. On the other hand, \eqref{eq:prop_OrdDerivInfo_tau(ell)Bound} shows that with $\widetilde G(T) \coloneqq \sumik A_{i, 2} F_{i, j}'(T) \prod_{\substack{1 \le r \le K\\ r \ne i}} F_{r, j}(T)$, we have $\ord_2\left(\widetilde G(T)\right) \le C_0$, so that $\ord_2(\widetilde G(4T+1)) \le C_0 + 2\deg \widetilde G(T) \le C_0 + 2(D_0-1)$. Combining these observations, we find that $h_1'(T) = 4 \widetilde G(4T+1) \prodik F_{i, j}(4T+1)^{A_{i, 2}-1}$ has $t \coloneqq \ord_2(h_1') \le 2D_0 + C_0$. Since any $2$-critical point of $h_1(T) = \prodik F_{i, j}(4T+1)^{A_{i, 2}}$ is a root of the polynomial $2^{-\ord_2(\widetilde G(4T+1))} \widetilde G(4T+1)$ mod $2$, it follows that the maximum multiplicity of such a $2$-critical point is no more than $\deg \widetilde G(4T+1) - 1 \le D_0 - 1$. As such, an application of Proposition \ref{prop:Cochrane}(ii)   
with $m \coloneqq e_2-2 \ge (2D_0 + C_0) +3 \ge t+3$ and $\chi \coloneqq \psi$, 
shows that 
the first term in \eqref{eq:Z2^e2 InTermsOf psi eta} has size at most $(12.5) \cdot 2^{2D_0+C_0} \cdot 2^{e_2(1-1/D_0)}$, proving our claim. This shows that $|Z_{2^{e_2}; ~ \chi_{1, 2}, \dots, \chi_{K, 2}} (\FijlistiOK)| \le 25 \cdot 2^{2D_0+C_0} \cdot 2^{e_2(1-1/D_0)}$, completing the proof of \eqref{eq:Cochrane to Each Component}.

Setting $C_1 \coloneqq D_0 C_0^{C_0}$ and combining \eqref{eq:S''PrepForCochrane} with \eqref{eq:Cochrane to Each Component}, we find that for each $j \in [N]$,
$$|Z_{\qtil; ~ \chi_1, \dots, \chi_K} (\FijlistiOK)| \le \left(\prod_{\substack{\ell^e \parallel \qtil\\e_\ell \le \kappa}} \phi(\ell^e)\right) \prod_{\substack{\ell^e \parallel \qtil\\ e_\ell \ge \kappa+1}} \left(\ell^{e-e_\ell} \cdot C_1 \ell^{e_\ell(1-1/D_0)}\right) \le \frac{(2 C_1)^{C_0} \phi(\qtil)}{A^{1/D_0}},$$
where $A \coloneqq \prod_{\ell^e \parallel \qtil: ~ e_\ell \ge \kappa+1} ~ \ell^{e_\ell}$ denotes the $(\kappa+1)$-full part of $\lcmcondchilist$, i.e, the largest $(\kappa+1)$-full divisor of $\lcmcondchilist$. For a divisor $d$ of $\qtil$, there are $\le d^K$ tuples $(\chi_1, \dots, \chi_K)$ of characters mod $\qtil$ for which $\lcmcondchilist = d$. Hence, summing the bound in the above display over all possible $\chiituplist$ occurring in the sum $\SDPr$, we obtain 
\begin{align*}\allowdisplaybreaks
|\SDPr| &\le \frac1{\phi(\qtil)^K} \largesum_{\substack{A \mid \qtil: ~A>1\\A\iskappaPlusOneFull}} ~ ~ ~ ~ ~ ~ ~ ~ ~ ~ \largesum_{\substack{d \mid \qtil\\(\kappa+1)\text{-full part of }d\text{ is }A}} ~ ~ d^K \cdot \frac{(2 C_1)^{C_0 N} \phi(\qtil)^N}{A^{N/D_0}}\\
&\ll \frac{\phi(\qtil)^N}{\phi(\qtil)^K} \cdot (2 C_1)^{C_0 N} \largesum_{\substack{A \mid \qtil: ~A>1\\A\iskappaPlusOneFull}} \frac1{A^{N/D_0-K}}.
\end{align*}
In the last step above, we have noted that for any $d$ dividing $\qtil$ whose $(\kappa+1)$-full part is $A$, we have $d \ll A$. Continuing, 
\begin{equation}\label{eq:SDPr Simplificn Step1}
\frac{|\SDPr|}{\phi(\qtil)^N} \ll \frac{(2 C_1)^{C_0 N}}{\phi(\qtil)^K} \left\{\prod_{\ell^e\parallel \qtil} ~ \left(1 + \sum_{\kappa+1 \le \nu \le e} \frac1{\ell^{\nu(N/D_0-K)}} \right)-1\right\}.
\end{equation}
Now if $N \ge KD_0 + 1$, then since $\kappa > C \cdot (30 D_0 C_0^{C_0})^{2C_0} \ge D_0(D_0+3)$, we see that the sum on $\nu$ above is no more than $2^{-\kappa (N/D_0-K)} \left(1-{2^{-1/D_0}}\right)^{-1} \le \frac{2^{D_0+2}}{2^{\kappa/D_0}} \le \frac12$. Hence $\log(1+ \sum_{\kappa+1 \le \nu \le e} {\ell^{-\nu(N/D_0-K)}}) \ll 2^{-\kappa(N/D_0-K)} \ll 2^{-\kappa N/D_0}$.
In addition, 
since $P(\qtil) \le C_0$, \eqref{eq:SDPr Simplificn Step1} gives 
\begin{equation}\label{eq:SDPrLargeN}
\frac{|\SDPr|}{\phi(\qtil)^N} \ll \frac{(2 C_1)^{C_0 N}}{\phi(\qtil)^K} \left\{\exp\left(O\left(\frac1{2^{\kappa N/D_0}}\right)\right)-1\right\} \ll \frac1{\phi(\qtil)^K} \cdot \left(\frac{(2 C_1)^{C_0}}{2^{\kappa/D_0}}\right)^N \ll \frac{C^{-N}}{\phi(\qtil)^K}, 
\end{equation}
where in the last step, we have recalled that $\kappa/D_0 > D_0^{-1} \cdot C \cdot (30 D_0 C_0^{C_0})^{2C_0} > C \cdot (2 C_1)^{C_0}$. Combining \eqref{eq:SDPrLargeN} with \eqref{eq:SPrimeExpr}, we deduce that 
\begin{equation}\label{eq:VNKtil(qtil;wi)_LargeN}
\frac{\#\VNKtilqtilwi}{\phi(\qtil)^N} = \frac{\SPr+\SDPr}{\phi(\qtil)^N} = \left(\frac{\phi(Q_0)}{\phi(\qtil)}\right)^K \left\{ \frac{\#\VNKtilQZwi}{\phi(Q_0)^N} + O\left(\frac1{C^{N}}\right)\right\},
\end{equation} 
uniformly for $N \ge KD_0+1$ and in coprime residues $\wilist$ to any modulus $q$. In particular, since $\qtil/\phi(\qtil) \le 2^{\omega(\qtil)} \ll 1$, we have 
\begin{equation}\label{eq:VNKtil(qtil;wi)_LargeN_Bound}
\frac{\#\VNKtilqtilwi}{\phi(\qtil)^N} \ll ~ \frac1{\phi(\qtil)^K} ~ \ll ~ \frac1{\qtil^K}.
\end{equation} 
Note that we did not make use of any invariant factor hypothesis to derive \eqref{eq:VNKtil(qtil;wi)_LargeN} or \eqref{eq:VNKtil(qtil;wi)_LargeN_Bound}.

On the other hand, for each $N \in [KD_0]$, we have $1 + \sum_{\kappa+1 \le \nu \le e} {\ell^{-\nu(N/D_0-K)}} \le \sum_{0 \le \nu \le e} {\ell^{\nu(K-N/D_0)}}$ $\ll e^{\bbm_{N=KD_0}} ~ \ell^{e(K-N/D_0)}$. Multiplying this over the $O(1)$ primes $\ell$ dividing $\qtil$ yields, from \eqref{eq:SDPr Simplificn Step1}, ${|\SDPr|}/{\phi(\qtil)^N}$ 
$\ll {\left(\prod_{\ell^e\parallel \qtil} ~ e\right)^{\bbm_{N=KD_0}}}\Big/{\qtil^{N/D_0}}.$
Combining this with the trivial bound $|\SPr|/\phi(\qtil)^N \ll \phi(\qtil)^{-K} \ll \qtil^{-K} \ll \qtil^{-N/D_0}$ coming from \eqref{eq:SPrimeExpr}, we find that for each $N \in [KD_0]$, we have
\begin{equation}\label{eq:VNKtil(qtil;wi)_SmallN}
\frac{\#\VNKtilqtilwi}{\phi(\qtil)^N} \ll \frac{\left(\prod_{\ell^e\parallel \qtil} ~ e\right)^{\bbm_{N=KD_0}}}{\qtil^{N/D_0}}, ~ ~ \text{ uniformly in }q\text{ and }(w_i)_{i=1}^K \in U_q^K.
\end{equation}

To establish Proposition \ref{prop:VNKCountmegageneral}(a), we multiply \eqref{eq:VNKtil(qtil;wi)_LargeN} with the relations \eqref{eq:VNKtil_LargePrimes_LargeN} for all $\ell^e \parallel q$ with $\ell>C_0$, noting that any such $\ell$ exceeds $B_0$ and hence automatically satisfies $\gcd(\ell-1, \beta(\GirlistiOK))=1$ for each $r \in [L]$ as $q$ satisfies $IFH(\GirlistiOK; B_0)$. (Note here that $\prod_{\ell \mid q: ~ \ell>C_0} \alpha_N^*(\ell) = \alphaStNq/\alphaNStQZ$.) 
On the other hand, (b) follows by multiplying the relations \eqref{eq:VNKtil_LargePrimes_SmallN} over all $\ell>C_0$, with \eqref{eq:VNKtil(qtil;wi)_SmallN} (resp. \eqref{eq:VNKtil(qtil;wi)_LargeN_Bound}) for each $N \in [KD_0]$ (resp. each $N \ge KD_0+1$). 
This establishes Proposition \ref{prop:VNKCountmegageneral} for $D_0>1$. 

Now we consider the case $D_0=1$, so that $K = 1$ and $G_{1, r}(T) \coloneqq G_r(T) \coloneqq R_r T + S_r$ for some integers $R_r$ and $S_r$ with $R_r \ne 0$. We set $F_j \coloneqq F_{1, j}$ for each $j \in [N]$ and show that  
\begin{equation}\label{eq:VN1Count_LargePrime_Linear}
\frac{{\#\widetilde{\mathcal V}_{N, 1}(\ell^e; w)}}{\phi(\ell^e)^N} = \frac{\alpha_N^*(\ell)}{\phi(\ell^e)} \left(1+O\left(\Big(\frac2{\ell-1}\Big)^{N-1}\right)\right),
\end{equation}
uniformly for $N \ge 1$ and for $\ell^e \parallel q$ with $\ell>C_0$. To this end, we start by using \eqref{eq:OrthogPrimePower} to write 
\begin{equation}\label{eq:Orthog_LargePrime_Linear}
\frac{{\#\widetilde{\mathcal V}_{N, 1}(\ell^e; w)}}{\phi(\ell^e)^N} = \frac{\alpha_N^*(\ell)}{\phi(\ell^e)} \left(1 + \frac1{\alpha_N^*(\ell)\phi(\ell^e)^N} \largesum_{\chi \ne \chi_{0, \ell}\bmod{\ell^e}} \overline\chi(w) \prodjN Z_{\ell^e; \chi} (F_j)\right).
\end{equation}
If $\condofchi = \ell^{e_0}$ for some $e_0 \in \{2, \dots, e\}$, then it is easy to see that $Z_{\ell^e; \chi} (F_j) = 0$ for any $j \in [N]$: this is immediate by orthogonality if $S_{j'}=0$, and follows from Proposition \ref{prop:Cochrane} otherwise, 
since the polynomial $T^{\phi(\ell^{e_0})} (R_{j'} T + S_{j'})$ has no $\ell$-critical points. 
On the other hand, if $\condofchi = \ell$, then $|Z_{\ell^e; \chi}(F_j)| = \ell^{e-1} \left|\sum_{v \bmod \ell} \chi(R_{j'} v + S_{j'}) - \chi(S_{j'})\right| = \ell^{e-1} \left|\sum_{u \bmod \ell} \chi(u) - \chi(S_{j'})\right| \le \ell^{e-1}$; here we have recalled that $\ell \nmid R_{j'}$ (by choice of $C_0$). 
Now \eqref{eq:VN1Count_LargePrime_Linear} follows by combining these observations with the fact that there are $\ell-2$  many characters mod $\ell^e$ with conductor $\ell$.

Letting $\qtil \coloneqq \prod_{\substack{\ell^e \parallel q\\\ell \le C_0}} \ell^e$ as before, we fix an integer $\kappa > C_0+3$, and write $\#\VNOnetilqtilw = \frac1{\phi(\qtil)} \largesum_{\chi \bmod \qtil} \overline\chi(w) \prodjN Z_{\qtil; ~ \chi} (F_j) = \SPr + \SDPr,$
where 
$\SPr$ again denotes the sum over those $\chi$ mod $\qtil$ for which $\cond(\chi)$ is ($\kappa+1$)-free. Then \eqref{eq:SPrimeExpr} continues to hold, and $\SDPr=0$, once again 
since there are no critical points. This yields
${\VNOnetilqtilw}/{\phi(\qtil)^N} = ({\phi(Q_0)}/{\phi(\qtil)}) \cdot ({\#\VNOnetilQZw}/\phi(Q_0)^N)$, which along with \eqref{eq:VN1Count_LargePrime_Linear} establishes Proposition \ref{prop:VNKCountmegageneral} in the remaining case $D_0=1$. \end{proof}
While proving Theorem \ref{thm:RestrictedInputSqfreeMod1}, we will also need the following variant of the Proposition \ref{prop:VNKCountmegageneral}, whose argument in a simpler version of that given for \eqref{eq:VNKtil_LargePrimes_SmallN}. Indeed applying \eqref{eq:OrthogPrimePower_Bound} with $e \coloneqq 1$, and recalling the two assertions around  \eqref{eq:Zell^eBoundFor e_0=1}, we obtain the following corollary.
\begin{cor}\label{cor:VNKCountSqfrMegaGen}
In the setting preceding Proposition \ref{prop:VNKCountmegageneral}, the following estimates hold uniformly in coprime residues $\wilist$ to squarefree moduli $q$. 
\begin{enumerate}
\item[\textbf{(a)}] For each fixed $N \ge 2K+1$, 
\begin{equation}\label{eq:VNKkCountGen_Sqfreeq_LargeN}
\frac{\#\VNKtilqwi}{\phi(q)^N} \ll 
\displaystyle{\frac1{\phi(q)^K} \exp\Big(O(\sqrt{\log q})\Big)},
\end{equation}
if $q$ satisfies $\IFHGirBZ$ for each $r \in [L]$.
\item[\textbf{(b)}] For each fixed $N \ge 1$,  
\begin{equation}\label{eq:VNKkCountGen_Sqfreeq_SmallN}
\frac{\#\VNKtilqwi}{\phi(q)^N} 
\ll \displaystyle{\frac1{q^{\min\{K, N/2\}}} \expOomegaq}.
\end{equation} 
\end{enumerate}
\end{cor}
Proposition \ref{prop:Vqwi_ReducnToBddModulus} is a special case of Proposition \ref{prop:VNKCountmegageneral}, with $L \coloneqq 1$ and $(G_{i, 1})_{i=1}^K \coloneqq (W_{i, k})_{i=1}^K$, so that $D_0 = D = \sumik \deg \Wik$, $(F_{i, j})_{i=1}^K = (W_{i, k})_{i=1}^K$, $\alphajtil(q) = \alpha_k(q)$, $\alphaStNq = \alpha_k(q)^N$, and $\VNKtilqwi = \VNKqwi$. (Here of course, all the quantities on the right hand side are in accordance with notation until the previous section.) This also completes the proof of Proposition \ref{prop:ConvenientMainTermUptoQ0}. 

\textbf{Remark.} Taking $K = L = N = 1$ and $G_{1, 1} = H \in \Z[T]$ with $\deg H = d > 1$ in \eqref{eq:VNKkCount_smallN_Gen}, we get 
\begin{equation}\label{eq:KonyaSub}
\frac1{\phi(q)}\#\{v \in U_q: H(v) \equiv w \pmod q\} \ll \frac{\expOomegaq}{q^{1/d}} \ll_\delta \frac1{q^{1/d-\delta}}
\end{equation}
for any fixed $\delta>0$. This is only slightly weaker than the results of Konyagin in \cite{konyagin79b, konyagin79a}.

In order to deduce Theorem \ref{thm:ConvenientMAINTerm} from Proposition \ref{prop:ConvenientMainTermUptoQ0}, we apply the orthogonality of Dirichlet characters to see that the main term in the right hand side of \eqref{eq:CountReducedToQ0} is equal to 
$$\frac{1}{\phi(q)^K} \sumnxfnq + \frac1{\phi(q)^K} \largesum_{\chiituplist \ne (\chi_{0, Q_0}, \dots, \chi_{0, Q_0}) \bmod Q_0} \Bigg(\prodik \overline\chi_i(a_i) \Bigg)\largesum_{n \le x} \bbmfnqOne \prodik \chi_i(f_i(n)).$$
Let $Q \coloneqq \prod_{\ell \mid q} \ell$ denote the \textsf{radical} of $q$. To obtain Theorem \ref{thm:ConvenientMAINTerm}, it remains to prove that each inner sum above is $\osumnxfnq$.  
For $Q \ll 1$, this follows by applying Theorem N to the divisor $Q^* \coloneqq \lcm[Q, Q_0] \ll 1$ of $q$. (Note that as $q$ lies in $\WUDkAdmSet$, so does $Q^*$, since $q$ and $Q^*$ have the same prime factors). So we may assume that $Q$ is sufficiently large.  
Theorem \ref{thm:ConvenientMAINTerm} would follow once we show the result below. Here $\lambda$ and $Q_0$ are as in Proposition \ref{prop:ConvenientMainTermUptoQ0}.
\begin{thm}\label{thm:RemCharacSums}
There exists a constant $\delta_0 \coloneqq \delta_0(\lambda)>0$ such that, uniformly in moduli $\qlelogxKZ$ lying in $\WUDkAdmSet$ and having sufficiently large radical, we have 
$$\largesum_{n \le x} \chi_1(f_1(n)) \cdots \chi_K(f_K(n)) \bbmfnqOne ~ ~ \ll ~ \frac{x^{1/k}}{(\log x)^{1-(1-\delta_0) \alpha_k(Q)}}$$
for all tuples of characters $\chiituplist \ne (\chi_{0, Q_0}, \dots, \chi_{0, Q_0}) \bmod Q_0$.
\end{thm}
Let $\CkQZ$ denote the set of tuples of characters $(\chi_1, \dots, \chi_K)$ mod $Q_0$, not all trivial, such that $\prodik \chi_i(\Wik(u))$ is constant on the set $R_k(Q_0) = \{u \in U_{Q_0}: W_k(u) \in U_{Q_0}\}$. 
To prove Theorem \ref{thm:RemCharacSums}, we separately consider the two cases when a tuple of characters mod $Q_0$ lies in $\CkQZ$ or not.
\section{Proof of Theorem \ref{thm:RemCharacSums} for nontrivial tuples of characters not in $\CkQZ$}\label{sec:HalaszReducnProof}
For any integer $d \ge 1$ and any nontrivial tuple $\psiituplist$ of characters mod $d$ not lying in $\mathcal C_k(d)$, we have $\left|\sum_{u \bmod d} ~ \chi_{0, d}(u) \psi_1(W_{1, k}(u)) \cdots \psi_K(W_{K, k}(u)) \right| ~ < ~ \alpha_k(d) \phi(d)$. With $\lambda$ as in Proposition \ref{prop:ConvenientMainTermUptoQ0}, we define the constant $\delta_1 \coloneqq \delta_1(\Wiklist; B_0) \in (0, 1)$ to be
\begin{multline*}
\max_{\substack{d \le \lambda\\\alpha_k(d) \ne 0}} ~ ~ \max_{\substack{\psiituplist \ne (\chi_{0, d}, \dots, \chi_{0, d}) \bmod d\\\psiituplist ~ \not\in ~ \mathcal C_k(d)}} ~ ~ \frac1{\alpha_k(d)\phi(d)} \left|\sum_{u \bmod d} ~ \chi_{0, d}(u) \psi_1(W_{1, k}(u)) \cdots \psi_K(W_{K, k}(u)) \right|.
\end{multline*}
Then since $Q_0 \le \lambda$, we have for any nontrivial tuple $\chiituplist \not\in \CkQZ$, 
\begin{equation}\label{eq:Q0CharacSum_Bound}
\left|\sum_{u \bmod Q_0} ~ \chi_{0, Q_0}(u) \chi_1(W_{1, k}(u)) \cdots \chi_K(W_{K, k}(u)) \right| \le \delta_1 \alpha_k(Q_0) \phi(Q_0).  
\end{equation}
We set $\delta \coloneqq (1-\delta_1)/2$ and $Y \coloneqq \exp((\log x)^{\delta/3})$. 
To establish Theorem \ref{thm:RemCharacSums} for all  
$\chiituplist \not\in \CkQZ$, it suffices to show that 
\begin{equation}\label{eq:HalaszReduc_MainRemaining}
\largesum_{\substack{n \le x\\ \nopgrYstpkOnedivn}} \chi_1(f_1(n)) \cdots \chi_K(f_K(n)) ~ \bbm_{(f(n), q)=1} ~ \ll ~ \frac{x^{1/k}}{(\log x)^{1-(\delta_1+\delta)\alpha_k}},    
\end{equation}
since by the arguments before \eqref{eq:Large power of prime exceeding y negligible}, the contribution of $n$'s 
not counted above is negligible. 
Writing any $n$ counted in \eqref{eq:HalaszReduc_MainRemaining} uniquely as $BMA^k$ (as in \eqref{eq:BMAkSplit1}), we see that the sum equals 
\begin{multline}\label{eq:BMAkSplit_forHalaszReducn}
\largesum_{\substack{B \le x\\ P(B) \le Y\\B \iskfree}} \bbm_{(f(B), q)=1} \left(\prodik \chi_i(f_i(B))\right) \largesum_{\substack{M \le x/B\\M\iskfull\\ P(M) \le Y}} \bbm_{(f(M), q)=1} \left(\prodik \chi_i(f_i(M))\right) \\ \largesum_{\substack{A \le (x/BM)^{1/k}}} \bbm_{P^-(A)>Y} ~ \bbm_{(f(A^k), q)=1} ~ \mu(A)^2 \prodik \chi_i(f_i(A^k))
\end{multline} 
Moreover, the arguments leading to the bound for $\Sigma_2$ in section \ref{sec:TechnicalPrep} show that the tuples $(B, M, A)$ having $M>x^{1/2}$ give negligible contribution 
to the above sum. It thus remains to consider the contribution of tuples $(B, M, A)$ with $M \le x^{1/2}$. To deal with such tuples, we will establish the following general upper bound uniformly for $X \ge \exp((\log Y)^2)$: 
\begin{equation}\label{eq:ASumBound}
\largesum_{\substack{A \le X}} \bbm_{P^-(A)>Y} ~ \bbm_{(f(A^k), q)=1} ~ \mu(A)^2 \prodik \chi_i(f_i(A^k)) \ll \frac X{(\log X)^{1-\alpha_k(\delta_1 + \delta/2)}}.
\end{equation}
We apply a quantitative version of Hal\'asz's Theorem \cite[Corollary III.4.12]{tenenbaum15} on the multiplicative function 
$F(A) \coloneqq \bbm_{P^-(A)>Y} ~ \bbm_{(f(A^k), q)=1} ~ \mu(A)^2 \prodik \chi_i(f_i(A^k))$, taking $T \coloneqq \log X$. 
This requires us to put, for each $t \in [-T, T]$, a lower bound on the sum below (which is the square of a certain ``pretentious distance"):
\begin{align}\label{eq:DXDef} \allowdisplaybreaks
\nonumber
\mathcal D(X; t) &\coloneqq \largesum_{p \le X} ~ \frac1p \left(1- \Ree\Bigg(\bbm_{p>Y} ~ \bbm_{(f(p^k), q)=1} ~ \mu(p)^2 ~ p^{-\ii t} \prodik \chi_i(f_i(p^k))\Bigg)\right)\\
&= \begin{multlined}[t]
(1-\alpha_k) \log_2 X + \alpha_k \log_2 Y + \largesum_{\substack{Y< p \le X\\(W_k(p), q)=1}} ~ \frac1p \left(1- \Ree\Bigg(p^{-\ii t} \prodik \chi_i(\Wik(p))\Bigg)\right)\\ + \OloglogqBddPower; 
\end{multlined}
\end{align}
here the second line uses Lemma \ref{lem:primesum}. 
To get this lower bound, we proceed analogously to the proof of \cite[Lemma 3.3]{PSR}. The key idea is to split the range of the above sum into blocks of small multiplicative width, so that the complex number $p^{-\ii t}$ is essentially constant for all $p$ lying in a given block. More precisely, we cover the interval $(Y, X]$ with finitely many disjoint intervals $\mcI \coloneqq \left(\eta, \eta(1+1/\log^2 X)\right]$ for certain choices of $\eta \in (Y, X]$, choosing the smallest $\eta$ to be $Y$ and allowing the rightmost endpoint of such an interval to jut out slightly past $X$ but no more than $X(1+1/\log^2 X)$. Then the last sum in \eqref{eq:DXDef} equals
\begin{align}\label{eq:CoveringWithIntervals}
\largesum_{\mcI} \largesum_{\substack{p \in \mcI\\(W_k(p), q)=1}} ~ \frac1p \left(1- \Ree\Bigg(p^{-\ii t} \prodik \chi_i(\Wik(p))\Bigg)\right) + O\left(\frac1{\log^3 X}\right)
\end{align}
Consider any $\mcI$ occurring in the sum above. For each $p \in \mcI$, we have
$$|p^{-\ii t} - \eta^{-\ii t}| \le \left|\int_{t \log \eta}^{t\log p} \exp(-i\varrho)\,d\varrho\right| \le |t \log p - t \log \eta| \le \frac{|t|}{\log^2 X} \le \frac1{\log X}.$$
This shows that each inner sum in \eqref{eq:CoveringWithIntervals} is equal to
\begin{equation}\label{eq:Extracting r^(-it)}
\largesum_{\substack{u \in U_q\\(W_k(u), q)=1}} ~ \left(1- \Ree\Bigg(\eta^{-\ii t} \prodik \chi_i(\Wik(u))\Bigg)\right) \largesum_{\substack{p \in \mcI\\p \equiv u \pmod q}} ~ \frac1p + O\left(\frac1{\log X}\largesum_{p \in \mcI} \frac1p\right)
\end{equation}
Note that $p = (1+o(1))\eta$ for all $p \in \mcI$. (Here and in what follows, the asymptotic notation refers to the behavior as $x \rightarrow \infty$, and is uniform in the choice of $\mcI$.) For parameters $Z, W$ depending on $X$, we write $Z \gtrsim W$ to mean $Z \ge (1+o(1))W$. By the Siegel Walfisz Theorem, 
$$\largesum_{\substack{p \in \mcI\\p \equiv u \pmod q}} ~ \frac1p ~ ~ \gtrsim ~ ~ \frac1\eta ~ \largesum_{\substack{p \in \mcI\\p \equiv u \pmod q}} ~1 ~\gtrsim ~~ ~\frac1{\phi(q)} \cdot \frac1\eta ~ \largesum_{p \in \mcI} 1 ~ ~ \gtrsim ~ ~ \frac1{\phi(q)} ~ \largesum_{p \in \mcI} ~\frac1p.$$
Hence the main term in \eqref{eq:Extracting r^(-it)} is
\begin{align*}\allowdisplaybreaks
\gtrsim \frac1{\phi(q)} ~ \largesum_{p \in \mcI} ~\frac1p ~ \largesum_{\substack{u \in U_q\\(W_k(u), q)=1}} ~ \left(1- \Ree\Bigg(\eta^{-\ii t} \prodik \chi_i(\Wik(u))\Bigg)\right)
\gtrsim ~ (\alpha_k - \alpha_k \delta_1) \left(\largesum_{p \in \mcI} ~\frac1p\right),
\end{align*}
where in the last step, 
we have used \eqref{eq:LiftPolyCoprimality} and \eqref{eq:Q0CharacSum_Bound} to see that
\begin{multline*}
\frac1{\phi(q)} \Bigg|\largesum_{\substack{u \in U_q\\(W_k(u), q)=1}} ~ \prodik \chi_i(\Wik(u)) \Bigg| 
= \frac{\alpha_k(q)}{\alpha_k(Q_0) \phi(Q_0)} \Bigg|\largesum_{\substack{r \bmod Q_0}} ~ \chi_{0, Q_0}(r) \prodik\chi_i(\Wik(r))\Bigg| \le \alpha_k \delta_1.
\end{multline*}
Inserting the bound obtained in the previous display into \eqref{eq:Extracting r^(-it)}, we find that each inner sum in \eqref{eq:CoveringWithIntervals} is 
$\gtrsim ~ \alpha_k (1-\delta_1) \sum_{p \in \mcI} ~1/p + O\left((\log X)^{-1} \sum_{p \in \mcI} 1/p\right)$. 
The $O$-term when summed over all $\mcI$ is $\ll (\log X)^{-1} \sum_{p \le 2X} p^{-1} \ll \log_2 X/\log X$. Thus, the main term in \eqref{eq:CoveringWithIntervals} is at least $\alpha_k \left(1-\delta_1 - \frac\delta2\right) (\log_2 X - \log_2 Y)$. Inserting this into \eqref{eq:DXDef} yields
$$\mathcal D(X; t) \ge \left(1-\alpha_k\left(\delta_1 + \frac\delta2\right)\right) \log_2 X + \alpha_k\left(\delta_1 + \frac\delta2\right)\log_2 Y + \OloglogqBddPower,$$
uniformly for $t \in [-T, T]$. As such,  
Corollary \cite[III.4.12]{tenenbaum15} establishes the claimed bound \eqref{eq:ASumBound}.

Now for each $M \le x^{1/2}$, we have $(x/BM)^{1/k} \gg x^{1/2k}$. 
Applying \eqref{eq:ASumBound} to each of the innermost sums in \eqref{eq:BMAkSplit_forHalaszReducn}, we see that the total contribution of all tuples $(B, M, A)$ with $M \le x^{1/2}$ is
\begin{align*}\allowdisplaybreaks
\ll \largesum_{\substack{B \ll 1}} ~ ~ \largesum_{\substack{M \le x^{1/2}: ~ M\iskfull\\ P(M) \le Y, ~ (f(M), q)=1}} ~ \frac{(x/BM)^{1/k}}{(\log x)^{1-\alpha_k(\delta_1 + \delta/2)}} ~ \ll ~ \frac{x^{1/k}}{(\log x)^{1-\alpha_k(\delta_1 + \delta)}},
\end{align*}
where we have 
used \eqref{eq:Msmoothrecip} (with $Y$ in place of $y$) and Lemma \ref{lem:primesum}. 
This proves \eqref{eq:HalaszReduc_MainRemaining}, and hence also Theorem \ref{thm:RemCharacSums} for all nontrivial tuples of characters $\chiituplist$ mod $Q_0$ not in $\CkQZ$. \hfill \qedsymbol 
\section{Proof of Theorem \ref{thm:RemCharacSums} for tuples of characters in $\CkQZ$}\label{sec:scourfield}
It suffices to consider the case when $x$ is an integer, and we will do so in the rest of the section. Our argument consists of suitably modifying the Landau--Selberg--Delange method for mean values of multiplicative functions (see for instance \cite[Chapter II.5]{tenenbaum15}), and to study the behavior of a product of $L$-functions raised to complex powers by accounting for the presence of Siegel zeros modulo $q$. This is partly inspired from work of Scourfield \cite{scourfield85} and will also need some results from her paper. We will denote complex numbers in the standard notation $s = \sigma + it$. \footnote{The parameters $\sigma$ and $\sigma_k$ (to be defined later) in this section have nothing to do with the divisor functions $\sigma_r(n) = \sum_{d \mid n} d^r$ mentioned in the introduction. We are not working with the divisor functions in this section.} To begin with, we consider the Dirichlet series
$$\Fchis \coloneqq \sum_{n \ge 1} \frac{\bbmfnqOne}{n^s} \prod_{i=1}^K \chi_i(f_i(n)) = \sum_{n \ge 1} \frac{\bbmfnQOne}{n^s} \prodchiifin$$
which is absolutely convergent in the half-plane $\sigma>1$. Let $\cchihat$ denote the constant value of $\prodik \chi_i(W_{i, k}(u))$ on the set $R_k(Q_0) = W_k^{-1}(U_{Q_0}) \cap U_{Q_0}$. In the rest of the section, we assume that the complex plane has been cut along the line $\sigma \le 1/k$ if $\alpha_k(Q)$ and $\cchihat$ are not both $1$, while if $\alpha_k(Q) = \cchihat = 1$, then the complex plane is cut along the line $\sigma \le \beta_e/k$. (In the last case, if there is no Siegel zero, then there is no cut.) Fix $\mu_0$ satisfying $\max\{0.7, k/(k+1)\} < \mu_0 < 1$.

\subsection{Analysis of the Dirichlet series.} 
We start by giving a meromorphic continuation of $\Fchis$ to a larger region. To do this, set $\ElQt \coloneqq \log(Q(|tk|+1))$ and recall that there exists an absolute constant $c_1>0$ such that the product $\prod_{\psi \bmod Q} L(s, \psi)$ has at most one zero $\beta_e$ (counted with multiplicity) in the region $\sigma>1-c_1/\log(Q(|t|+1))$, called the ``Siegel zero", which is necessarily real and simple. If $\beta_e$ exists, then it is a root of $L(s, \psi_e)$ for some real character $\psi_e$ mod $Q$, which we will be referring to as the ``exceptional character". By reducing $c_1$ if necessary, we may assume that $c_1<1-\mu_0$, and that the conductor of $\psi_e$ (which is squarefree) is large enough that it is not $(D+2)$-smooth.   
\begin{lem}\label{lem:FchisMeroCont}
The Dirichlet series $\Fchis$ is absolutely convergent on the half-plane $\sigma> \frac1k$, where it satisfies
\begin{equation}\label{eq:FchiMeroContProd}
\Fchis = \FOneskcchihat ~ \gskcchihat ~ \GchiOnes ~ \GchiTwos
\end{equation}
with
\begin{align*}\allowdisplaybreaks
&\FOnesk = \left(\prodQOnemidQ \prodPrimChar \LskpsiGammapsi \right)^{\alpha_k(Q)}\\
&\gsk = \left(\prodQOnemidQ \prodPrimChar \prod_{\ell \mid \frac{Q}{Q_1}} \OneminPsiellksGamPsi \right)^{\alpha_k(Q)}, ~ \gamma(\psi) = \frac1{\alpha_k(Q) \phi(Q)}\largesum_{\substack{v \in U_Q\\W_k(v) \in U_Q}} \chibar(v).
\end{align*}
Here, the functions $\FOnesk$, $\gsk$, $\GchiOnes$ and $\GchiTwos$ satisfy the following properties:
\begin{enumerate}
\item[(i)] $\FOnesk$ is holomorphic and nonvanishing in the region $\left\{s:\sigma>\OnekOneMinuscOneLQtk\right\}\sm\left\{\OneByk, \betaeByk\right\}$. 
\item[(ii)] $\gsk$ and $\GchiOnes$ are holomorphic and nonvanishing in the half-plane $\sigma> \mu_0/k$, and we have, uniformly for all $s$ in this region,
\begin{equation}\label{eq:LogerivgGchiOne}
\max\left\{\AbsLDgsk, \AbsLDGchiOnes \right\} \ll \max\{1, (\log Q)^{1-\sigma k}\} ~ \log\log Q.    
\end{equation}
\item[(iii)] $\GchiTwos$ is holomorphic in the half-plane $\sigma> \mu_0/k$, wherein $|\GchiTwos|, |\GchiTwosDer| \ll 1$.
\end{enumerate}
\end{lem}
\begin{proof}
For all $s$ in the region $\sigma>1$, we can use the Euler product of $\Fchis$ to write 
\begin{multline}\label{eq:FchiEulerProdConeq1}
\Fchis 
= \left(\prodbInUQWkbInUQ \prodpbmodQ \OneminOnepksMincchi\right) \cdot \left(\prodpMidQWkpInUQ \OneminOnepksMincchi\right)\\
\cdot \prod_p \left(1+\sum_{v \ge 1} \frac{\bbmfpvQOne}{\pvs} \prodchiifipv\right) \OneminOneWkpQOnepkscchihat
\end{multline}
Since $q$ and $Q$ are supported on the same primes, $Q$ is also $k$-admissible. By Lemma \ref{lem:kfreepartbdd} and the fact that $(\chi_1, \cdots, \chi_K) \in \mathcal C_k(Q_0)$, we thus find that
\begin{multline}\label{eq:ProdConvergenceCalculn}
\prod_{p \gg 1} \left(1+\sum_{v \ge 1} \frac{\bbmfpvQOne}{\pvs} \prodchiifipv\right) \OneminOneWkpQOnepkscchihat\\ = \prod_{p \gg 1} \left(1+\frac{\cchihat \bbmWkpQOne}{\pks} + O\left(\frac1{p^{(k+1)\sigma}}\right)\right) \OneminOneWkpQOnepkscchihat = \prod_{p \gg 1} \left(1 + O\left(\frac1{p^{(k+1)\sigma}}\right)\right), 
\end{multline}
which is an absolutely convergent product in the half plane $\sigma>1/k$, showing the absolute convergence of the Dirichlet series $\Fchis$ in the same half plane. 

Now for $\sigma>1/k$, the orthogonality of Dirichlet characters mod $Q$ and the fact that $\log L(sk, \psi) = \sum_{p, v} \psi(p^v)/{p^{vsk}}$ show that the logarithm of the first double product in \eqref{eq:FchiEulerProdConeq1} is equal to 
\begin{multline*}\allowdisplaybreaks
\cchihat\sumbInUQWkbInUQ \left\{\sum_p \frac1{\phi(Q)} \sum_{\psi \bmod Q} \psibarb \frac{\psi(p)}{p^{ks}} + \sum_{v \ge 2}\sum_{p \equiv b \pmod Q} \frac1{v p^{vks}} \right\} \\ 
= \alpha_k(Q)\cchihat \sumPsimodQ \gamma(\psi) \log \Lskpsi + \cchihat\sumbInUQWkbInUQ \sum_{v \ge 2} \left(\sumpbmodQ \frac1{v\pvks} - \sum_{p: ~ p^v \equiv b \pmod q} \frac1{v\pvks}\right).
\end{multline*}
We insert this into \eqref{eq:FchiEulerProdConeq1}, noting that $L(sk, \psi) = L(sk, \psi^*) \prod_{\ell \mid \frac Q{Q_1}} (1-\psi^*(\ell)/\ell^{sk})$ and that $\gamma(\psi) = \gamma(\psi^*)$ if the primitive character $\psi^*$ mod $Q_1$ induces $\psi$ mod $Q$. This yields \eqref{eq:FchiMeroContProd}, with
$$\GchiTwos \coloneqq \prod_{p \le C_k} \left(1+\sum_{v \ge 1} \frac{\bbmfpvQOne}{\pvs} \prodchiifipv\right) \OneminOneWkpQOnepkscchihat$$ 
and
\begin{multline}\label{eq:GchiOnesDef}
\GchiOnes \coloneqq \prod_{p > C_k} \left(1+\sum_{v \ge 1} \frac{\bbmfpvQOne}{\pvs} \prodchiifipv\right) \OneminOneWkpQOnepkscchihat\\
\cdot \prodpMidQWkpInUQ \OneminOnepksMincchi \cdot \exp\left(\cchihat \sumbInUQWkbInUQ \sum_{v \ge 2} \left(\sumpbmodQ \frac1{v\pvks} - \sum_{p: ~ p^v \equiv b \pmod q} \frac1{v\pvks}\right)\right),
\end{multline}
where $C_k>2^{k/\mu_0}$ is a constant exceeding any $k$-free integer $n$ satisfying $\gcd(f(n), q)=1$; recall that by Lemma \ref{lem:kfreepartbdd}, $C_k$ can be chosen to depend only on $\{\Wiv\}_{\substack{1 \le i \le K\\1 \le v \le k}}$ (and $\mu_0$). 

Now (i) follows by the result quoted before the statement of the Lemma and (iii) is immediate by a mechanical calculation. It is also clear that $\gsk$ is holomorphic and nonvanishing in the half-plane $\sigma>0$ and the assertion of \eqref{eq:LogerivgGchiOne} relevant to $g(sk)$ is an immediate consequence of \cite[Lemma 9(ii)]{scourfield85}. To show the assertions for $\GchiOnes$, we 
recall that for each prime $p>C_k$, the first local factor defining $\GchiOnes$ in \eqref{eq:GchiOnesDef} is $1+\cchihat\bbmWkpQOne/\pks + O(p^{-(k+1)\sigma})$, 
whereupon a computation analogous to \eqref{eq:ProdConvergenceCalculn} shows that the first product (over primes $p>C_k$) in \eqref{eq:GchiOnesDef} is absolutely convergent and defines a holomorphic function in the half plane $\sigma>\mu_0/k$. (Here is it important that $\mu_0/k>1/(k+1)$.) Likewise the exponential factor in \eqref{eq:GchiOnesDef} defines a holomorphic function in the same half plane, hence so does $\GchiOnes$. To see that $\GchiOnes$ is also nonvanishing in this region, we need only see that 
the condition $p>C_k>2^{k/\mu_0}$ guarantees the nonvanishing of each of the factors in the (absolutely convergent) product over $p>C_k$. Finally, a straightforward computation using \eqref{eq:GchiOnesDef} shows that for $\sigma>\mu_0/k$, we have
$$\LDGchiOnes = -\cchihat k \sumpMidQWkpInUQ \frac{\log p}{\pks} + O(1) \ll \sum_{p \mid Q} \frac{\log p}{p^{k\sigma}},$$
completing the proof of \eqref{eq:LogerivgGchiOne} via \cite[Lemma 3(i)(a)]{scourfield85}.
\end{proof}
Our objective is to relate the sum in Theorem \ref{thm:RemCharacSums} to the Dirichlet series $\Fchis$ by an effective version of Perron's formula, and shift the contour to the left of the line $\sigma=1/k$. As such, we will need the following proposition in order to estimate the resulting integrals. 

To set up, we choose $\epsilon_1 \coloneqq \epsilon_1(\lambda)$ to be a constant (depending only on $\lambda$) satisfying $0< \epsilon_1 < 1-\cos(2\pi/d)$ for any positive integer $d \le \lambda$. Consider the functions 
\begin{align*}\allowdisplaybreaks
&\FchiTils \coloneqq \FOneskcchihat ~ \gskcchihat ~ \GchiOnes\\
&\HchiTils \coloneqq \FchiTils \sMinOnekAlphcchi \sMinBetaekAlphcchiGam, ~ \Hchis \coloneqq \frac{\FchiTils}s \sMinOnekAlphcchi,
\end{align*}
where here and in what follows, any term or factor involving $\beta_e$ is to be understood as omitted if the Siegel zero doesn't exist. By assertions (i) and (ii) of the previous lemma, we see that:
\begin{itemize}
\item $\FchiTils$ is holomorphic and nonvanishing in the region $\left\{s:\sigma>\OnekOneMinuscOneLQtk, ~ s \ne \OneByk, \betaeByk\right\}$, 
\item $\Hchis$ is holomorphic and nonvanishing in the region $\left\{s:\sigma>\OnekOneMinuscOneLQtk, ~ s \ne \betaeByk\right\}$, 
\item $\HchiTils$ is holomorphic and nonvanishing in the region $\left\{s:\sigma>\OnekOneMinuscOneLQtk\right\}$ 
\end{itemize} (Recall our branch cut conventions elucidated at the start of the section.) Let $T \coloneqq \exp(\sqrt{\log x})$.
\begin{prop}\label{prop:ResidualIntgeralsBoundingTools} We have the following bounds:
\begin{enumerate}
\item[(i)] $\AbsHchiOnek \ll \logxalphakepsFiv$.
\item[(ii)] $|\HchiTils| \ll \logxalphakepsFour$ uniformly for real $s$ satisfying $\frac1k\left(1-\frac{c_1}{4\log Q}\right) \le s \le \frac1k$.
\item[(iii)] $|\Fchis| \ll (\log x)^{(1/2 + \epsilon_1)\alpha_k(Q)}$ uniformly for complex numbers $s$ satisfying $\sigma \ge \ScaledOneMincOneTwoElqt$, $|t| \le T$ and $|s-\theta/k| \gg 1/\ElQt$ for $\theta \in \{1, \beta_e\}$.
\item[(iv)] Uniformly in real $s \le 1/k$ satisfying $s \ge \frac1k\left(\frac23 + \frac{\beta_e}3\right)$ (if the Siegel zero exists) or $s \ge \frac1k\left(1 - \frac{c_1}{4 \log Q}\right)$ (otherwise), we have $$\left|\HchiFracOnek G_{\chi, 2}\left(\frac1k\right) - \Hchis \GchiTwos\right| \ll (\log x)^{(1/20+\alpha_k(Q)/5)\epsilon_1} \left(\frac1k-s\right).$$ 
\end{enumerate}    
\end{prop}
\begin{proof} The following general observation will play an important role in our arguments: We have $|\HchiTils| \asymp |\widetilde{H}_\chi(w)|$ uniformly in complex numbers $s$ and $w$ satisfying $\TextIm(s) = \TextIm(w) \eqqcolon t$, $|s-w| \ll \ElQt^{-1}$ and $\Ree(w) \ge \Ree(s) \ge \ScaledOneMincOneTwoElqt$. 

Indeed by the definitions of $\HchiTils$ and $\FchiTils$, we have
\begin{equation}\label{eq:HchiTilLD}
\left|\frac{\widetilde{H}_\chi'(z)}{\widetilde{H}_\chi(z)}\right| = \left|\cchihat k \left(\frac{F_1'(kz)}{F_1(kz)} + \frac{\alpha_k(Q)}{kz-1} - \frac{\alpha_k(Q) \gamma(\psi_e)}{kz-\beta_e}\right) + \cchihat k \frac{g'(kz)}{g(kz)} + \frac{G_{\chi, 1}'(z)}{G_{\chi, 1}(z)} \right| \ll \ElQt
\end{equation}
uniformly for complex numbers $z = u+it$ satisfying $u \ge \frac1k\left(1-\frac{c_1}{2\ElQt}\right)$. Here in the last step, we have applied \eqref{eq:LogerivgGchiOne} and \cite[Lemma 15(i)]{scourfield85}, the latter with $\xi(t) \coloneqq \exp(6 \ElQt)$. The general observation now follows by writing $\log\left(\widetilde{H}_\chi(w)/\HchiTils\right) = \int_{\Ree(s)}^{\Ree(w)} \widetilde{H}'_\chi(u+it)/{\widetilde{H}_\chi(u+it)} \, \mathrm{d}u$.

(i) Let $b_k(t) \coloneqq \frac1k\left(1+\frac{c_3}{\ElQt}\right)$ for some absolute constant $c_3>0$. By the above observation and the definitions of $\FchiTils$, $\HchiTils$ and $\Hchis$, it follows that 
\begin{align}\label{eq:AbsHchiFracOnek Beginnings}\allowdisplaybreaks
&\AbsHchiFracOnek \ll \AbsHchiTilFracOnek \BetaFactor \ll |\HchiTil(b_k(0))| ~ \BetaFactor\\
&\ll |\FchiTil(b_k(0))| ~ (\log Q) \BetaFactorTwice \ll |F_1(k b_k(0)) g(k b_k(0))|^{\Ree(\cchihat)} ~ (\log Q)^2 \BetaFactorTwice. \nonumber
\end{align}
Here in the last bound, we have noted that $|\GchiOne(b_k(0))| \ll \log_2 Q$, as is evident from the fact that $\prodpMidQWkpInUQ (1-p^{-k b_k(0)})^{-1} \ll \exp(\sum_{p \mid Q} 1/p) \ll \exp(\sum_{p \le \omega(Q)} 1/p) \ll \log \omega(Q) \ll \log_2 Q$. 

Now proceeding as in \cite[Lemma 8]{scourfield85}, we see that for all $s$ with $\sigma>1/k$, we have
\begin{equation}\label{eq:f(n^k), Q Dir Ser}
\sum_{n \ge 1} \frac{\bbm_{(f(n^k), Q)=1}}{n^{ks}} = F_1(ks) ~ g(ks) ~ \widetilde G(s),
\end{equation}  
where 
\begin{multline*}
\widetilde G(s) = \prod_p \left(1+\sum_{v \ge 2} \frac1{\pvks} (\bbm_{(f(p^{kv}), Q)=1} - \bbm_{(W_k(p), Q)=1} ~ \bbm_{(f(p^{k(v-1)}), Q)=1})\right)\cdot \\\prodpMidQWkpInUQ \left(\OneminOnepks\right)^{-1}
\cdot \exp\left(\sumbInUQWkbInUQ \sum_{v \ge 2} \left(\sumpbmodQ \frac1{v\pvks} - \sum_{p: ~ p^v \equiv b \pmod Q} \frac1{v\pvks}\right)\right).
\end{multline*}
Uniformly for $s$ with $\sigma \ge 1/k$, we observe that the infinite product above has size at least $1-\sum_{p, v \ge 2} 1/{p^v} \gg 1$ and at most $\exp(\sum_{p, v \ge 2} 1/p^v) \ll 1$. Likewise, the exponential factor has size 
$\asymp 1$ in the same region. Moreover, for $\sigma \ge 1/k$, the product over $p \mid Q$ is $\asymp \exp(\sum_{p \mid Q: ~ (W_k(p), Q)=1} ~ p^{-k \sigma})$, which is $\gg 1$ and $\ll \exp(\sum_{p \mid Q} ~ p^{-1}) \ll \log_2 Q$. Putting these observations together, we find that $1 \ll \widetilde G(s) \ll \log_2 Q$ for $\sigma \ge 1/k$. 
Applying the lower bound on $s \coloneqq b_k(0)$, the equality \eqref{eq:f(n^k), Q Dir Ser} yields
$$|F_1(k b_k(0)) ~ g(k b_k(0))| \ll \sum_{n \ge 1} \frac{\bbm_{(f(n^k), Q)=1}}{n^{k b_k(0)}} \le \zeta(k b_k(0)) = \frac1{kb_k(0)-1} + O(1) \ll \log Q,$$
so that from \eqref{eq:AbsHchiFracOnek Beginnings}, we obtain $\AbsHchiOnek \ll (\log Q)^3 \BetaFactorTwice$. Subpart (i) now follows as $Q \le (\log x)^{K_0}$ and $1-\beta_e \gg_{\epsilon_1} Q^{-\epsilon_1/20{K_0}} \gg_{\epsilon_1} (\log x)^{-\epsilon_1/20}$ by Siegel's Theorem.  

(ii) By the observation made at the start of the proof, we have $|\HchiTils| \ll |\HchiTil(1/k)| \ll \AbsHchiOnek \BetaFactor \ll \AbsHchiOnek (\log x)^{\alpha_k(Q)\epsilon_1/20}$. The result now follows from (i).

(iii) By the aforementioned observation, we have $|\HchiTils| \ll |\HchiTil(b_k(t)+it)|$, and since $|s-\theta/k| \gg 1/\ElQt$, we have ${b_k(t)+it-\theta/k} \asymp {s-\theta/k}$ for $\theta \in \{1, \beta_e\}$. Thus 
$|\FchiTils| \ll |\FchiTil(b_k(t)+it)|$. Using  \eqref{eq:GchiOnesDef}   
and replicating the arguments that led to the bounds on $\widetilde G(s)$ above, we also obtain $(\log_2 Q)^{-1} \ll G_{\chi, 1}(s) $ $\ll \log_2 Q$ for $\sigma \ge 1/k$, so that $|\FchiTils| \ll (\log_2 Q) \cdot |F_1(k(b_k(t)+it)) g(k(b_k(t)+it))|^{\Ree(\cchihat)}$. From \eqref{eq:f(n^k), Q Dir Ser} and the bounds on $\widetilde G(s)$, we thus get $|\FchiTils| \ll (\log_2 Q) \left|\sum_{n \ge 1} {\bbm_{(f(n^k), Q)=1}}/{n^{k(b_k(t) + it)}}\right|^{\Ree(\cchihat)} \ll (\log_2 Q) \left(\sum_{n \ge 1} {\bbm_{(f(n^k), Q)=1}}/{n^{kb_k(t)}}\right)^{\Ree(\cchihat)}$, whence $|\FchiTils| \ll (\log_2 Q)^2 |F_1(k b_k(t)) g(k b_k(t))|^{\Ree(\cchihat)} \ll (\log_2 Q)^3 |\FchiTil(b_k(t))|$. By definitions of $b_k(t)$ and $\HchiTil(b_k(t))$, we have $|\FchiTils| \ll (\log_3 x)^3 |\HchiTil(b_k(t))| \ElQt^{\alpha_k(Q)} (1-\beta_e)^{-\alpha_k(Q)}$. Finally, recall that $|t| \le T = \exp(\sqrt{\log x})$, that $1-\beta_e \gg_{\epsilon_1} (\log x)^{-\epsilon_1/20}$, and that $|\HchiTil(b_k(t))| \ll |\HchiTil(1/k)| \ll  (\log x)^{\alpha_k(Q) \epsilon_1/4}$ (by subpart (ii) the general observation at the start of the proof). This yields $|\FchiTils| \ll (\log x)^{\alpha_k(Q)(1/2+\epsilon_1)}$, and Lemma \ref{lem:FchisMeroCont}(iii) applies. 

(iv) It suffices to show that uniformly for $s$ satisfying the same conditions as in this subpart,
\begin{equation}\label{eq:Propn(iv)ETS}
|\Hchis| + |H_\chi'(s)| \ll (\log x)^{\alpha_k(Q)\epsilon_1/5} \left(\log Q + \frac1{1-\beta_e}\right).    
\end{equation}
(Here as usual, the second term on the right is omitted if there is no Siegel zero, otherwise it dominates.) Indeed once we establish \eqref{eq:Propn(iv)ETS}, then from the bound $1-\beta_e \gg_{\epsilon_1} (\log x)^{-\epsilon_1/20}$, it follows that $|\Hchis| + |H_\chi'(s)| \ll (\log x)^{(1/20+\alpha_k(Q)/5)\epsilon_1}$, which combined with Lemma \ref{lem:FchisMeroCont}(iii) and the observation $\left|\HchiOnek G_{\chi, 2}\left(1/k\right) - \Hchis \GchiTwos\right| = \left|\int_s^{1/k} (H_\chi(u) \GchiTwo(u))' \, \mathrm{d}u\right|$ completes the proof of the subpart. 
To show \eqref{eq:Propn(iv)ETS}, we recall that $H_\chi(s)$ is non-vanishing for $s$ as in the subpart. Further \eqref{eq:HchiTilLD} applies with $z=s$ for all $s$ considered in this subpart, yielding 
$$\left|\frac{{H}_\chi'(s)}{{H}_\chi(s)}\right| = \left|\frac{\widetilde{H}_\chi'(s)}{\widetilde{H}_\chi(s)} - \frac1s + \frac{\alpha_k(Q) \cchihat \gamma(\psi_e)}{s-\beta_e/k} \right| \ll \ElQ(0) + 1 + \frac1{1-\beta_e} \ll \log Q + \frac1{1-\beta_e}.$$
As a consequence,
$$\left|\log\frac{H_\chi(1/k)}{H_\chi(s)}\right| = \left|\int_s^{1/k} \frac{H_\chi'(u)}{H_\chi(u)} \, \mathrm{d}u\right| \ll \left(\frac1k-s\right) \left(\log Q + \frac1{1-\beta_e}\right) \ll 1,$$
showing that $|H_\chi(s)| \asymp \AbsHchiOnek$ uniformly for all $s$ in the statement. Collecting these bounds, we obtain for all such $s$,
$$|\Hchis| + |H_\chi'(s)| \ll \AbsHchiFracOnek + \left|\frac{{H}_\chi'(s)}{{H}_\chi(s)}\right| \cdot \left|\frac{H_\chi(s)}{H_\chi(1/k)}\right| \cdot \AbsHchiFracOnek \ll \AbsHchiFracOnek \left(\log Q + \frac1{1-\beta_e}\right),$$
so that the desired bound \eqref{eq:Propn(iv)ETS} now follows from subpart (i). This concludes the proof.
\end{proof}

\subsection{Perron's formula and the contour shifts} 
We first show that there is some $X$ sufficiently close to $x$ for which the error term arising from an effective Perron's formula is small. 
\begin{lem}\label{lem:PerronMinimizer}
Let $h \coloneqq x/\log^2 x$. There exists a positive integer $X \in (x, x+h]$ satisfying $$\sum_{\substack{3X/4 < n < 5X/4\\n \ne X}} \frac{\bbmfnQOne}{|\log(X/n)|} \ll X^{1/k} \log X.$$
\end{lem}
\begin{proof}
This would follow once we show that 
\begin{equation}\label{eq:averaged}
\sum_{x < X \le x+h} ~ ~ \sum_{\substack{3X/4 < n < 5X/4\\n \ne X}} \frac{\bbmfnQOne}{|\log(X/n)|} \ll x^{1/k} h \log x,
\end{equation}
with the outer sum being over integers $X \in (x, x+h]$. (Recall that $x \in \Z^+$ in this entire section.) To show this, we write the sum on the left hand side as $S_1 + S_2$, where $S_1$ denotes the contribution of the case $3X/4 < n \le X-1$. Writing any $n$ contributing to $S_1$ as $X-v$ for some integer $v \in [1, X/4)$, we see that $|\log(X/n)| = -\log(1-v/X) \gg v/X \gg v/x$. Recalling that $n=Bm$ for some $k$-free $B$ of size $O(1)$ and some $k$-full $m$, we thus have 
\begin{align*}\allowdisplaybreaks
S_1 &\le \sum_{3x/4 < n < x+h} ~ \sum_{\substack{x < X \le x+h\\n+1 \le X < 4n/3}} \frac{\bbmfnQOne}{|\log(X/n)|} \ll x\sum_{B \ll 1} ~ \sum_{\substack{\frac{3x}{4B} < m < \frac{x+h}{B}\\m \text{ is }k\text{-full}}} \sum_{\substack{1 \le v < \frac{x+h}4\\x<v+Bm \le x+h}} \frac1v\\
&\ll x \sum_{1 \le v \le \frac{x+h}4} \frac1v \sum_{B \ll 1} ~ \sum_{\substack{\frac{x-v}{B} < m \le \frac{x-v+h}{B}\\m \text{ is }k\text{-full}}} 1 \ll x \log x \left(x^{1/k}\frac hx + x^{1/(k+1)}\right) \ll x^{1/k} h \log x,
\end{align*}
where we have bounded the last inner sum on $m$ using the Erd\"os-Szekeres estimate on the count of $k$-full integers (see \cite{ES34}). This shows that the sum $S_1$ is bounded by the right hand expression in \eqref{eq:averaged}, and similarly so is the sum $S_2$, establishing \eqref{eq:averaged}.
\end{proof}
To complete the proof of Theorem \ref{thm:RemCharacSums}, it suffices to establish the bound therein for $X$ in place of $x$, for once we do so, we may simply note that 
$$\Big| \sum_{x < n \le X} \chi_1(f_1(n)) \cdots \chi_K(f_K(n)) \bbmfnqOne \Big| \le \sum_{x < n \le X} \bbmfnQOne \le \sum_{B \ll 1} \sum_{\substack{\frac{x}{B} < m \le \frac{X}{B}\\m \text{ is }k\text{-full}}} 1 \ll \frac{x^{1/k}}{\log^2 x}.$$
To show the bound in Theorem \ref{thm:RemCharacSums} for $X$, we start by applying an effective version of Perron's formula \cite[Theorem II.2.3]{tenenbaum15}. To bound the resulting error, we use Lemma \ref{lem:PerronMinimizer} and note that 
\begin{align*}\allowdisplaybreaks
&X^{\ScaledOnekOnelogX} \left(\sum_{n \le 3X/4} + \sum_{n \ge 5X/4}\right) \frac{\bbmfnQOne}{T|\log(X/n)| n^{\ScaledOnekOnelogX}} \ll \frac{X^{1/k}}T \sum_{B \ll 1} ~ \sum_{\substack{m \ge 1\\m\iskfull}} \frac1{m^{\ScaledOnekOnelogX}}\\
&\ll \frac{X^{1/k}}T \prod_p \left(1 + \frac1{p^{1+1/\log X}} + O\left(\frac1{p^{1+1/k}}\right)\right) \ll \frac{X^{1/k}}T \exp\left(\sum_p \frac1{p^{1+1/\log X}}\right) \ll \frac{X^{1/k} \log X}T,
\end{align*}
with the last bound above being a consequence of Mertens' Theorem along with the fact that 
$$\sum_{p>X} \frac1{p^{1+1/\log X}} \le \sum_{j \ge 0} \sum_{X^{2^j}< p \le X^{2^{j+1}}} \frac1{p^{1+1/\log X}} \le \sum_{j \ge 0} \exp(-2^j) \sum_{X^{2^j}< p \le X^{2^{j+1}}} \frac1p \ll 1.$$
(Recall that $T = \exp(\sqrt{\log x}) \ge \exp\left(\frac12 \sqrt{\log X}\right)$.) As such, \cite[Theorem II.2.3]{tenenbaum15} yields 
\begin{equation}\label{eq:PerronApp}
\sum_{n \le X} \chi_1(f_1(n)) \cdots \chi_K(f_K(n)) \bbmfnQOne = \frac1{2\pi i} \int_{\ScaledOnekOnelogX - iT}^{\ScaledOnekOnelogX + iT} \frac{\Fchis X^s}s \, \ds + O\left(\frac{X^{1/k} \log X}T\right).  
\end{equation}
Our arguments will be divided into three possibilities:\\ 
\textsf{Case 1:} When $(\alpha_k(Q), \cchihat) \ne (1, 1)$ and there is a Seigel zero $\beta_e$ mod $Q$.\\
\textsf{Case 2:} When $(\alpha_k(Q), \cchihat) \ne (1, 1)$ and there is no Seigel zero mod $Q$.\\
\textsf{Case 3:} When $(\alpha_k(Q), \cchihat) = (1, 1)$.\\
In Case 1, we will be assuming henceforth that $\beta_e > 1-\frac{5c_1}{24 \log Q}$; otherwise decreasing $c_1$ reduces to Case 2. Let $\beta^* \coloneqq \frac23 + \frac{\beta_e}3$ and $\sigma_k(t) \coloneqq \ScaledOneMincOneFourElqt$, so that $\frac{\beta_e}k>\sigma_k(0)$. Let $\delta, \delta_1 \in (0, \beta_e/10k)$ satisfy $\sigma_k(0)<\frac{\beta_e}k-2\delta_1 < \frac{\beta_e}k+2\delta_1<\frac{\beta^*}k< \frac1k-2\delta$. Consider the contours
\begin{itemize}
\item $\Gamma_2$, the horizontal segment traversed from $\ScaledOnekOnelogX + iT$ to $\sigma_k(T) + iT$.
\item $\Gamma_3$, the part of the curve $\sigma_k(t) + it$ traversed from $t=T$ to $t=0$.
\item $\Gamma_4 \coloneqq \Gamma_4(\delta_1)$, the segment traversed from $\sigma_k(0)$ to $\beta_e/k-\delta_1$ \textbf{above} the branch cut. 
\item $\Gamma_5 \coloneqq \Gamma_5(\delta_1)$, the semicircle of radius $\delta_1$ centered at $\beta_e/k$, lying in the upper half plane and traversed clockwise. 
\item $\Gamma_6 \coloneqq \Gamma_6(\delta_1)$, the segment traversed from $\beta_e/k+\delta_1$ to $\beta^*/k$ \textbf{above} the branch cut. 
\item $\Gamma_7 \coloneqq \Gamma_7(\delta)$, the segment traversed from $\beta^*/k$ to $1/k-\delta$ \textbf{above} the branch cut. 
\item $\Gamma_8 \coloneqq \Gamma_8(\delta)$, the circle of radius $\delta$ centered at $1/k$, traversed clockwise from the point  $1/k-\delta$ above the branch cut to its reflection below the branch cut.
\item $\Gamma_4^* \coloneqq \Gamma_4^*(\delta)$, the segment traversed from $\sigma_k(0)$ to $1/k-\delta$ \textbf{above} the branch cut.
\item $\Gamma_5^* \coloneqq \Gamma_5^*(\delta_1)$, the circle of radius $\delta_1$ centered at $\beta_e/k$, traversed clockwise from the point $\beta_e/k-\delta_1$ above the branch cut to its reflection below the branch cut.
\end{itemize} 
Here $\Gamma_5^*(\delta_1)$ is relevant only when our branch cut is along $\sigma \le \beta_e/k$ (i.e., when $\alpha_k(Q) = \cchihat = 1$ and $\beta_e$ exists), while the rest of the contours are defined irrespective of the branch cut. For a contour $\Omega$, let $-\overline\Omega$ denote the contour given by the complex conjugate of $\Omega$ traversed in the opposite direction and \textbf{below} the respective branch cuts. (Note that $-\overline{\Gamma_5}$ is still traversed \textbf{clockwise} but below the branch cut.) 
We define the contour $\Gamma_1$ by
$$\Gamma_1 \coloneqq \begin{cases}
\sum_{j=2}^8 \Gamma_j + \sum_{j=2}^7 (\MinGammajBar), & \text{ under Case 1}\\
\Gamma_2 + \Gamma_3 + \Gamma_4^* + \Gamma_8 + (\MinGammaFourStBar) + (\MinGammaThrBar) + (\MinGammaTwoBar), & \text{ under Case 2}\\
\sum_{j=2}^4 \Gamma_j + \Gamma_5^* + \sum_{j=2}^4 (\MinGammajBar), & \text{ under Case 3.}
\end{cases}$$
In Case 3, if $\beta_e$ doesn't exist, then there is no branch cut and $\Gamma_4$, $\overline\Gamma_4$ and $\Gamma_5^*$ are excluded from $\Gamma_1$. In all three cases, the integrand in \eqref{eq:PerronApp} is analytic in the region enclosed by $\Gamma_1$ and the segment joining $\ScaledOnekOnelogX - iT$ and $\ScaledOnekOnelogX + iT$. 
(Note that if $\cchihat=1$, the definitions of $\WUDkAdmSet$ and $\GchiOne$, $\GchiTwo$ in Lemma \ref{lem:FchisMeroCont} give $G_{\chi, 2}(1/k) = 0$, canceling the simple pole of $F_1(sk)$ at $s=1/k$. In particular, this happens in Case 3.) 
So   
\begin{equation}\label{eq:PerronApp2}
\sum_{n \le X} \chi_1(f_1(n)) \cdots \chi_K(f_K(n)) \bbmfnQOne = -\frac1{2\pi i} \int_{\Gamma_1} \frac{\Fchis X^s}s \, \ds + O\left(\frac{X^{1/k} \log X}T\right).  
\end{equation} 
We now proceed to estimate the integrals occurring on the right hand side above. In the following proposition, any result about an integral is valid whenever the corresponding contour is a part of $\Gamma_1$: so for instance, the assertion on $\Gamma_8$ (resp. $\Gamma_5^*$) holds under Cases 1 or 2 (resp. Case 3), those on $\Gamma_5$ and $\Gamma_6$  hold under Case 1, and the bound involving $\Gamma_4$ holds under Cases 1 and 3. Let $I_j$ (resp. $\overline{I_j}$, $I_j^*$) denote the corresponding integral along $\Gamma_j$ (resp. $\GammajBar$, $\Gamma_j^*$). 
\begin{prop}\label{prop:ResIntBounds} We have the following bounds:
\begin{itemize}
\item[(i)] $|I_2| + |\overline{I_2}| + |I_3| + |\overline{I_3}| \ll X^{1/k} \exp(-\kappa_0 \sqrt{\log X})$ for some constant $\kappa_0 \coloneqq \kappa_0(c_1, k)>0$.
\item[(ii)] $\max\{|I_4 + \overline{I_4}|, |I_6 + \overline{I_6}|\} \ll X^{1/k} \exp(-\sqrt{\log X})$ uniformly in $\delta, \delta_1$ as above.
\item[(iii)] $\lim_{\delta_1 \rightarrow 0+} |I_5| = \lim_{\delta_1 \rightarrow 0+} |\overline{I_5}| = \lim_{\delta_1 \rightarrow 0+} |I_5^*| = \lim_{\delta \rightarrow 0+} |I_8| = 0$.
\end{itemize}
\end{prop}
\begin{proof} To show subpart (i), we use the fact that since $\beta_e>1-5c_1/24\log Q$, any $s$ lying on $\Gamma_2$, $\Gamma_3$ or their conjugates satisfies the requirements of Proposition \ref{prop:ResidualIntgeralsBoundingTools}(iii). As such, (i) follows immediately from Proposition \ref{prop:ResidualIntgeralsBoundingTools}(iii) and the fact that $|s| \gg |t|+1$ for all $s$. 

For subpart (ii), we note that for all $s \in \Gamma_4$, we have $(s-1/k)^{-\alpha_k(Q) \cchihat} = (1/k-s)^{-\alpha_k(Q) \cchihat} ~ e^{-i \pi \alpha_k(Q) \cchihat}$ and $(s-\beta_e/k)^{\alpha_k(Q) \cchihat \gamma(\psi_e)} = (\beta_e/k-s)^{\alpha_k(Q) \cchihat \gamma(\psi_e)} ~ e^{i \pi \alpha_k(Q) \cchihat \gamma(\psi_e)}$. (This is clear if the branch cut is along $\sigma \le 1/k$, and also if the branch cut is along $\sigma \le \beta_e/k$ which is when $(\alpha_k(Q), \cchihat) = (1, 1)$.) Likewise, for all $s \in \overline{\Gamma_4}$, we have $(s-1/k)^{-\alpha_k(Q) \cchihat} = (1/k-s)^{-\alpha_k(Q) \cchihat} ~ e^{i \pi \alpha_k(Q) \cchihat}$ and $(s-\beta_e/k)^{\alpha_k(Q) \cchihat \gamma(\psi_e)} = (\beta_e/k-s)^{\alpha_k(Q) \cchihat \gamma(\psi_e)} ~ e^{-i \pi \alpha_k(Q) \cchihat \gamma(\psi_e)}$. Since $e^{\pm i \pi \alpha_k(Q) \cchihat (\gamma(\psi_e)-1)} \ll 1$, the definitions of $\FchiTils$ and $\HchiTils$ show that
\begin{align*}\allowdisplaybreaks
|I_4 + \overline{I_4}| \ll \left|\int_{\sigma_k(0)}^{\beta_e/k-\delta_1} \frac{\HchiTils \GchiTwos X^s}s \left(\frac1k-s\right)^{-\alpha_k(Q) \cchihat} \left(\frac{\beta_e}k-s\right)^{\alpha_k(Q) \cchihat \gamma(\psi_e)} \, \ds\right|.
\end{align*}
But now by Lemma \ref{lem:FchisMeroCont}(iii) and Proposition \ref{prop:ResidualIntgeralsBoundingTools}(ii), we see that 
\begin{align*}\allowdisplaybreaks
|I_4 + \overline{I_4}| &\ll X^{\beta_e/k} (\log X)^{\alpha
_k(Q) \epsilon_1/4} \BetaFactor \int_{\sigma_k(0)}^{\beta_e/k-\delta_1} \left(\frac{\beta_e}k-s\right)^{\alpha_k(Q) \Ree(\cchihat \gamma(\psi_e))} \, \ds\\
&\ll X^{\beta_e/k} (\log X)^{3\alpha
_k(Q) \epsilon_1/10} \cdot \left(\frac{\beta_e}k-\sigma_k(0)\right)^{1+\alpha_k(Q) \Ree(\cchihat \gamma(\psi_e))} \ll X^{1/k} \exp(-\sqrt{\log X}).
\end{align*}
Here we have recalled that $\beta_e \le 1-c(\epsilon_1)/Q^{\epsilon_1/20{K_0}} \le 1-c(\epsilon_1)/(\log X)^{\epsilon_1/20}$ for some constant $c(\epsilon_1)>0$, and (as argued before Lemma \ref{lem:FchisMeroCont}) that $Q_e \coloneqq \cond(\psi_e)$ has a prime factor $\ell_e>D+2$, which upon factoring $\psi_e = \prod_{\ell \mid Q} \psi_{e, \ell}$ with $\psi_{e, \ell}$ being a character mod $\ell$, led to
\begin{equation}\label{eq:alphak(Q)Gamma(psie)Bound}
\alpha_k(Q) |\gamma(\psi_e)| \le \alpha_k(Q) \prod_{\ell \mid Q_e} \left|\frac{\sum_{v:~vW_k(v) \in U_\ell} \overline\psi_{e, \ell}(v)}{\alpha_k(\ell)(\ell-1)}\right| \le \frac1{\ell_e-1} \Bigg|\sum_{\substack{v \bmod {\ell_e}\\ W_k(v) \equiv 0 \pmod{\ell_e}}} \overline\psi_{e, \ell}(v)\Bigg| \le \frac D{D+1}.  
\end{equation}
This shows the desired bound on $I_4$ in (ii), and the assertion for $I_6$ is entirely analogous. 

Coming to subpart (iii), we parametrize the points of $\Gamma_5$ by $s=\beta_e/k + \delta_1 e^{i\theta}$ where $\pi \ge \theta \ge 0$. Since $\widetilde M \coloneqq \sup_{\left|s-\frac{\beta_e}k\right| \le \frac12\left(\frac{\beta_e}k-\sigma_k(0)\right)} ~ |\HchiTils|$ is finite, we have for all sufficiently small $\delta_1>0$, 
$$|I_5| \ll \widetilde M \int_0^\pi X^{\beta_e/k+\delta_1} \left(\frac{1-\beta_e}k-\delta_1\right)^{-\alpha_k(Q) \Ree(\cchihat)} \delta_1^{1+\alpha_k(Q)\Ree(\cchihat \gamma(\psi_e))} \, \mathrm d\theta \ll   \frac{\widetilde M X^{\beta_e/k+\delta_1} \delta_1^{1/(D+1)}}{\left(\frac{1-\beta_e}k-\delta_1\right)^{\alpha_k(Q)}},$$ 
where we have again seen that $1+\alpha_k(Q)\Ree(\cchihat \gamma(\psi_e)) \ge 1/(D+1)$ by \eqref{eq:alphak(Q)Gamma(psie)Bound}. The last expression shows that $\lim_{\delta_1 \rightarrow 0+} |I_5| = 0$, and the assertions on $|\overline{I_5}|$ and $|I_5^*|$ are proved similarly. The same argument also shows that $|I_8| \ll M^* X^{1/k+\delta} \delta^{1-\alpha_k(Q) \Ree(\cchihat)} \left(\frac{1-\beta_e}k-\delta\right)^{-\alpha_k(Q)}$ for all sufficiently small $\delta>0$, where $M^* = \sup_{\left|s-\frac1k\right| \le \frac{1-\beta^*}k} ~ |\HchiTils|$.  This yields $\lim_{\delta \rightarrow 0+} |I_8|=0$, because $\alpha_k(Q) \Ree(\cchihat) < 1$ whenever $(\alpha_k(Q), \cchihat) \ne (1, 1)$.
\end{proof}
Now in case 3, we let $\delta_1 \downarrow 0$ in \eqref{eq:PerronApp2} and invoke the relevant assertions of Proposition \ref{prop:ResIntBounds} to obtain $\sum_{n \le X} \chi_1(f_1(n)) \cdots \chi_K(f_K(n)) \bbmfnQOne \ll X^{1/k} \exp(-\kappa_1 \sqrt{\log X})$ for some constant $\kappa_1>0$. Hence to complete the proof of Theorem \ref{thm:RemCharacSums}, it suffices to assume that $(\alpha_k(Q), \cchihat) \ne (1, 1)$. In case 1, we obtain, by letting $\delta \downarrow 0$ and $\delta_1 \downarrow 0$ in \eqref{eq:PerronApp2},
\begin{equation}\label{eq:PostVanishing}
\sum_{n \le X} \chi_1(f_1(n)) \cdots \chi_K(f_K(n)) \bbmfnQOne = -
\lim_{\delta \rightarrow 0+} \frac{I_7 + \overline{I_7}}{2\pi i} + O(X^{1/k} \exp(-\kappa_1 \sqrt{\log X})).
\end{equation}
By an argument analogous to that given for Proposition \ref{prop:ResIntBounds}(ii), it is easy to see that the above limit exists. Furthermore, writing $(s-1/k)^{-\alpha_k(Q) \cchihat} = (1/k-s)^{-\alpha_k(Q) \cchihat} ~ e^{\pm i \pi \alpha_k(Q) \cchihat}$ as before, we see that the limit in \eqref{eq:PostVanishing} is equal to 
\begin{equation*}\label{eq:ScourMainTerm}
\frac{\sin(\pi \alpha_k(Q) \cchihat)}\pi \int_{\beta^*/k}^{1/k} \Hchis \GchiTwos X^s \left(\frac1k-s\right)^{-\alpha_k(Q) \cchihat} \ds,
\end{equation*}
We write the above integral as $\HchiOnek G_{\chi, 2}(1/k) I_1 - I_2$, where $I_1 \coloneqq \int_{\beta^*/k}^{1/k} X^s (1/k-s)^{-\alpha_k(Q) \cchihat} \, \ds$. Letting $s=1/k-u/\log X$, and using $\beta^* = 2/3+\beta_e/3 \le 1-c(\epsilon_1)/3(\log X)^{\epsilon_1/20}$ along with a standard bound on the tail of the integral defining a Gamma function \cite[Lemma 7]{scourfield85}, we get
$$I_1 = \frac{X^{1/k}}{(\log X)^{1-\alpha_k(Q) \cchihat}} \left\{\Gamma(1-\alpha_k(Q) \cchihat) + O(\exp(-\sqrt{\log X}))\right\}.$$
Now using Proposition \ref{prop:ResidualIntgeralsBoundingTools}(iv) and making the same change of variable, we find that 
$$I_2 \ll (\log X)^{\left(\frac1{20}+\frac{\alpha_k(Q)}5\right)\epsilon_1} \int_{\beta^*/k}^{1/k} X^s \left(\frac1k-s\right)^{1-\alpha_k(Q) \Ree(\cchihat)} \ds \ll \frac{X^{1/k}}{(\log X)^{2-\alpha_k(Q) \Ree(\cchihat) - (1/20+\alpha_k(Q)/5)\epsilon_1}}$$
as $\Gamma(2-\alpha_k(Q)\Ree(\cchihat)) \ll 1$. Collecting estimates, we obtain from \eqref{eq:PostVanishing},
\begin{multline}\label{eq:MasterEstim}
\sum_{n \le X} \bbmfnQOne \prodik \chi_1(f_1(n)) = \frac{\HchiOnek \GchiTwo\left(1/k\right)}{\Gamma(\alpha_k(Q) \cchihat)} \cdot \frac{X^{1/k}}{(\log X)^{1-\alpha_k(Q)\cchihat}} \Big(1+O(\exp(-\sqrt{\log X}))\Big)\\ + O\left(\frac{X^{1/k}}{(\log X)^{2-\alpha_k(Q)\Ree(\cchihat) - (1/20 + \alpha_k(Q)/5)\epsilon_1}}\right),
\end{multline}
by  
the reflection formula for the Gamma function and as $\Gamma(z) \gg 1$ for all $z$ with $|z| \le 2$. 

If $\cchihat \ne 1$, then $\Ree(\cchihat) \le \cos(2\pi/\phi(Q_0)) < 1-\epsilon_1$. Lemma \ref{lem:FchisMeroCont}(iii) and Proposition \ref{prop:ResidualIntgeralsBoundingTools}(i) yield
\begin{equation*}
\sum_{n \le X} \bbmfnQOne \prodik \chi_1(f_1(n)) 
\ll \frac{X^{1/k}}{(\log X)^{1-\alpha_k(Q)(\Ree(\cchihat) + \epsilon_1/5)}} \ll \frac{X^{1/k}}{(\log X)^{1-\alpha_k(Q)(1-\delta_0)}},
\end{equation*}
with $\delta_0 \coloneqq \delta_0(\lambda) \coloneqq \min\{3\epsilon_1/4, 1-\epsilon_1/2\}$. On the other hand, if $\cchihat = 1$, then since $q \in \WUDkAdmSet$, we must have $\GchiTwo(1/k)=0$ (as observed before \eqref{eq:PerronApp2}). Hence, \eqref{eq:MasterEstim} yields  
$$\sum_{n \le X} \chi_1(f_1(n)) \cdots \chi_K(f_K(n)) \bbmfnQOne
\ll \frac{X^{1/k}}{(\log X)^{2-\alpha_k(Q) - (1/20 + \alpha_k(Q)/5)\epsilon_1}} \ll \frac{X^{1/k}}{(\log X)^{1-\alpha_k(Q)(1-\delta_0)}},$$
completing the proof of Theorem \ref{thm:RemCharacSums} in case 1.

Finally in case 2, \eqref{eq:PerronApp2} and Proposition \ref{prop:ResIntBounds} lead to the following analogue of \eqref{eq:PostVanishing}:
\begin{equation}
\sum_{n \le X} \chi_1(f_1(n)) \cdots \chi_K(f_K(n)) \bbmfnQOne = -
\lim_{\delta \rightarrow 0+} \frac{I_4^* + \overline{I_4^*}}{2\pi i} + O(X^{1/k} \exp(-\kappa_0 \sqrt{\log X})).
\end{equation}
An argument entirely analogous to the one given above leads to the sharper variant of \eqref{eq:MasterEstim} with the $\exp(-\sqrt{\log X})$ replaced by $\exp\left(-\frac{c_1\log X}{8k{K_0}\log_2 X}\right)$, completing the proof of Theorem \ref{thm:RemCharacSums}.

This finally concludes the proof of Theorem \ref{thm:ConvenientMAINTerm}. 
In order to establish Theorems \ref{thm:UnrestrictedInput_SqfrPoly} to \ref{thm:HigherPolyControl1}, we thus need to appropriately bound the contributions of inconvenient $n$'s considered in the respective theorems. We take this up in the next several sections. 

\section{Equidistribution to restricted moduli: Proof of Theorem \ref{thm:UnrestrictedInput_SqfrPoly}}\label{sec:UnrestrictedInputGen_Proof}
By Theorem \ref{thm:ConvenientMAINTerm}, it remains to show that 
\begin{equation}\label{eq:UnrestrictedTheorem Inconvenient Contribn}
\largesum_{\substack{n \le x \inconv\\\forallifinaimodq}} 1 ~ ~ = ~ \ophiqKsumnxfq \quadxtoinfty,  
\end{equation}
uniformly in coprime residues $\aifam$ to $k$-admissible moduli $\qlelogxKZ$, under any one of the conditions (i)-(iii) of Theorem \ref{thm:UnrestrictedInput_SqfrPoly}. 

To show this, we set $z \coloneqq x^{1/\log_2 x}$ and recall that, by \eqref{eq:Smoothfnq=1}, \eqref{eq:Large power of prime exceeding y negligible} and \eqref{eq:fnqcoprimeEstimate}, the $n$'s that are either $z$-smooth or divisible by the $(k+1)$-th power of a prime exceeding $y$ give negligible contribution to the left hand side of \eqref{eq:UnrestrictedTheorem Inconvenient Contribn} in comparison to the right hand side. 
The remaining $n$ can be written in the form $m P^k$, where $P \coloneqq P(n) > z$, $P_{Jk}(m) \le y$, $m$ is not divisible by the $(k+1)$-th power of a prime exceeding $y$, and $\gcd(m, P)=1$, so that $f_i(n) = f_i(m) \Wik(P)$. 
Given $m$, the number of possible $P$ is, by the Brun-Titchmarsh inequality, 
$$\ll \frac{V_{1, q}''}{\phi(q)} \cdot \frac{(x/m)^{1/k}}{\log (z/q)} \ll \frac{V_{1, q}''}{\phi(q)} \cdot \frac{x^{1/k} \log_2 x}{m^{1/k}\log x},$$
where $V_{1, q}'' \coloneqq \max \left\{\#\mathcal V_{1, K}^{(k)}\big(q; \wifam\big): ~ \wifam \in U_q^K\right\}$. Summing this over possible $m$, we get 
$$\largesum_{\substack{n \le x\inconv\\P(n)>z; ~ \nopgrystpkOnedivn\\\forallifinaimodq}} 1 ~ ~ \ll ~ \frac{V_{1, q}''}{\phi(q)} \cdot \frac{x^{1/k}}{(\log x)^{1-\alpha_k \epsilon/2}}\exp\big(O((\log_3 x)^2 + (\log_2 (3q))^{O(1)})\big)$$
via \eqref{eq:mInconvRecip}. By Proposition \ref{prop:fnqcoprimecount}, the quantity on the right hand side above is negligible compared to the right hand side of \eqref{eq:UnrestrictedTheorem Inconvenient Contribn} whenever $q^{K-1} V_{1, q}'' \ll (\log x)^{(1-2\epsilon/3) \alpha_k}$. But this does hold under any one of conditions (i)-(iii) in the statement of Theorem \ref{thm:UnrestrictedInput_SqfrPoly}, because:
\begin{enumerate}
\item[(i)] $V_{1, q}'' \ll 1$ if at least of one of $\Wikset$ is linear. 
\item[(ii)] $V_{1, q}'' \ll \Dmin^{\omega(q)}$ if $q$ is squarefree, since $\#\mathcal V_{1, K}^{(k)}(\ell; \wifam) \le \Dmin$ for all $\ell \gg 1$.
\item[(iii)] $V_{1, q}'' \ll_\delta q^{1-1/\Dmin+\delta}$ by \eqref{eq:KonyaSub}. With $\delta \coloneqq \frac{\epsilon}{4(1-\epsilon)}\left(K-\frac1{\Dmin}\right)$, this yields $q^{K-1} V_{1, q}'' \ll (\log x)^{(1-2\epsilon/3)\alpha_k}$ under condition (iii) of the theorem.
\end{enumerate}
This establishes \eqref{eq:UnrestrictedTheorem Inconvenient Contribn}, completing the proof of Theorem \ref{thm:UnrestrictedInput_SqfrPoly}. \hfill \qedsymbol 
\subsection{Optimality in the ranges of $q$ in Theorem \ref{thm:UnrestrictedInput_SqfrPoly}.} 
\label{subsec:UnrestrictedInput_Optimality} 
In all our examples below, $\{\Wik\}_{i=1}^K$ $\subZT$ will be nonconstant 
with $\prodik \Wik$ separable over $\Q$. Then $\beta(\Wiklist) = 1$, guaranteeing that any integer satisfies $IFH(W_{1, k}, \dots, W_{K, k}; 1)$. We claim that there exists a constant $\widetilde C$ $\coloneqq \widetilde C(W_{1, k}, \dots, W_{K, k})$ such that any $\widetilde C$-rough $k$-admissible integer $q$ lies in $\WUDkAdmSet$. Indeed, viewing a character of $U_q^K$ as a tuple of characters mod $q$,\footnote{Here $U_q^K$ is the direct product of $U_q$ taken $K$ times.} the condition \eqref{eq:JtWUD_Crit} becomes vacuously true whenever $\mathcal T_k(q) \coloneqq \{(W_{1, k}(u), \cdots, W_{K, k}(u)) \in U_q^K: u \in U_q\}$ generates the group $U_q^K$. Now under the canonical isomorphism $U_q^K \to \prod_{\ell^e \parallel q} U_{\ell^e}^K$, the set $\mathcal T_k(q)$ maps to $\prod_{\ell^e \parallel q} \mathcal T_k(\ell^e)$. Thus by \cite[Lemma 5.13]{narkiewicz84}, if $\mathcal T_k(q)$ does not generate $U_q^K$, then there is some $\ell^e \parallel q$ and some tuple of characters $(\psi_1, \cdots, \psi_K) \ne \chiZelltuplist$ mod $\ell^e$ for which $\prodik \psi_i(\Wik(u))$ is constant on the set $R_k(\ell^e)$. Our claim now follows from \cite[Lemma 5]{narkiewicz82}. 

Fix any $k \in \NatNos$. Let $\CZeroTilde > \max\{\widetilde C, 4KD\}$ be any constant depending only on the polynomials $\Wikset$, which also 
exceeds the size of the leading coefficient and (nonzero) discriminant of $\prodik \Wik$. 
Then by Theorem N, $f_1, \dots, f_K$ are jointly weakly equidistributed modulo any (fixed) $\CZeroTilde$-rough $k$-admissible integer. 
Fix a prime $\ell_0 > \CZeroTilde$, and consider any nonconstant polynomials $\{\Wiv\}_{\substack{1 \le i \le K\\1 \le v \le k-1}}$ $\subZT$ all of whose coefficients are divisible by $\ell_0$, so that $\alpha_v(\ell_0) = 0$ for each $v < k$. Our moduli $q$ will have $P^-(q) = \ell_0$, so that 
$\alpha_v(q) = 0$ for all $v<k$. In each example below, we will show that $\alpha_k(q) \ne 0$, so that $q$ is $k$-admissible and lies in $\WUDkAdmSet$ by definition of $\CZeroTilde$. 
The constant $K_0$ (in the assumption $\qlelogxKZ$) is taken large enough in terms of $\{W_{i, k}\}_{i=1}^K$.  

\textbf{Optimality under condition (i)}. We show that for any $K \ge 2$, the range of $q$ in Theorem \ref{thm:UnrestrictedInput_SqfrPoly}(i) is optimal, -- even if \textit{all} of $\Wiklist$ are assumed to be linear, for \textit{any} choice of (pairwise coprime) linear functions. Indeed, consider $\Wik(T) \coloneqq c_i T+b_i \in \Z[T]$ for nonzero integers $c_i$ and integers $b_i$ satisfying $b_i/c_i \ne b_j/c_j$ for all $i \ne j$. Then $\prodik \Wik$ is clearly separable in $\Q[T]$. Choose a nonzero integer $b$ such that $\prod_{i=1}^K (c_i b + b_i) \ne 0$. Let $\CZeroTilde > \max\{|b|, |c_i b+b_i|: 1 \le i \le K\}$ be any constant satisfying the aforementioned requirements, so that any $q$ with $P^-(q)=\ell_0>\CZeroTilde$ is coprime to $b$ and to $\prodik \Wik(b) = \prodik (c_i b+b_i)$. Thus  
$\alpha_k(q) \ne 0$ and $q \in \WUDkAdmSet$. 
Now any prime $P \le x^{1/k}$ satisfying $P \equiv b \pmod q$ also satisfies $f_i(P^k) = \Wik(P) \equiv c_i b + b_i \pmod q$ for all $i \in [K]$. The Siegel--Walfisz Theorem thus shows that there are $\gg {x^{1/k}}\big/{\phi(q) \log x}$ 
many $n \le x$ satisfying $f_i(n) \equiv c_i b + b_i \pmod q$ for all $i \in [K]$. 
By Proposition \ref{prop:fnqcoprimecount}, this last expression grows strictly faster than $\phi(q)^{-K} \#\{n \le x: (f(n), q)=1\}$ as soon as $q \ge (\log x)^{(1+\epsilon)\alpha_k/(K-1)}$ for any fixed $\epsilon \in (0, 1)$, showing that the range of $q$ in Theorem \ref{thm:UnrestrictedInput_SqfrPoly} under condition (i) is essentially optimal.  
Note that with $Y \in [2(1+\epsilon)\log_2 x/(K-1), (K_0/2)\log_2 x]$, the squarefree integer $q \coloneqq \prod_{\ell_0 \le \ell \le Y} \ell$ satisfies all desired conditions; in particular $(\log x)^{(1+\epsilon)/(K-1)} \le q \le (\log x)^{K_0}$ and $P^-(q) = \ell_0$. 

\textbf{Optimality under condition (ii).} 
To show that the range of squarefree $q$ in Theorem \ref{thm:UnrestrictedInput_SqfrPoly}(ii) is optimal, we define $\Wik(T) \coloneqq \prod_{1 \le j \le d} (T-2j) + 2(2i-1) \in \Z[T]$ for some fixed $d>1$. Eisenstein's criterion at the prime $2$ shows that each $\Wik$ is irreducible in $\Q[T]$, and the distinct $\Wik$'s differ by a constant, making $\prodik \Wik$ separable over $\Q$. Now $2 \in U_q$, and $\Wik(2) = 2(2i-1) \le 2(2K-1) < 4KD < \CZeroTilde < P^-(q)$ for each $i \in [K]$. Thus, $q \in \WUDkAdmSet$ and $(2(2i-1))_{i=1}^K \in U_q^K$. Further, any prime $P$ satisfying $\prod_{1 \le j \le d} (P-2j) \equiv 0 \pmod q$ also satisfies $f_i(P^k) = \Wik(P) \equiv 2(2i-1) \pmod q$ for each $i$. Since $2d = 2\deg \Wik < 4KD< P^-(q)$, we see that $2, 4, \dots, 2d$ are all distinct coprime residues modulo each prime dividing $q$, whereupon it follows that the congruence $\prod_{1 \le j \le d} (v-2j) \equiv 0 \pmod q$ has exactly $d^{\omega(q)}$ distinct solutions $v \in U_q$ for squarefree $q$. Hence, there are $\gg \frac{d^{\omega(q)}}{\phi(q)} \cdot \frac{x^{1/k}}{\log x}$ many primes $P \le x^{1/k}$ satisfying $f_i(P^k) \equiv 2(2i-1) \pmod q$ for all $i$, so there are also at least as many $n \le x$ for which all $f_i(n) \equiv 2(2i-1) \pmod q$.   
The last expression grows strictly faster than $\phi(q)^{-K} \#\{n \le x: (f(n), q)=1\}$ as soon as $q^{K-1} \Dmin^{\omega(q)} = q^{K-1} d^{\omega(q)} > (\log x)^{(1+\epsilon)\alpha_k}$ for any fixed $\epsilon>0$, showing that the range of $q$ in Theorem \ref{thm:UnrestrictedInput_SqfrPoly}(ii) is essentially optimal. 

Note that it is possible to construct squarefree $\qlelogxKZ$ satisfying the much stronger requirement 
that $d^{\omega(q)} > (\log x)^{(1+\epsilon)\alpha_k}$ (and $P^-(q) = \ell_0$). Indeed, let $q \coloneqq \prod_{\ell_0 \le \ell \le Y} \ell$ for some $Y \le (K_0/2)\log_2 x$. Then $\omega(q) = \sum_{\ell_0 \le \ell \le Y} 1 \ge Y/2\log Y$, while by the Chinese Remainder Theorem and the Prime Ideal Theorem, $\alpha_k(q) \le \kappa'/\log Y$ for some constant $\kappa' \coloneqq \kappa'(\Wiklist; \ell_0)$. 
So we need only choose $Y \in (4 \kappa' \log_2 x/\log d, (K_0/2) \log_2 x)$ to have $\qlelogxKZ$ and $d^{\omega(q)} > (\log x)^{(1+\epsilon)\alpha_k}$.  

For future reference, we observe that any $n$ of the form $P^k$ with $P$ a prime exceeding $q$ satisfies $P(n_k)>q$. Hence in the above setting, we have shown the stronger lower bound
\begin{equation}\label{eq:PknGreaterThanqInsuff}
\largesum_{\substack{n \le x: ~ P(n_k)>q\\(\forall i)~ f_i(n) \equiv 2(2i-1) \pmod q}} 1 ~ \ge \largesum_{\substack{q< P \le x^{1/k}\\\prod_{1 \le j \le d} (P-2j) \equiv 0 \pmod q}} 1 ~ \gg ~ \frac{d^{\omega(q)}}{\phi(q)} \cdot \frac{x^{1/k}}{\log x}.    
\end{equation}
\textbf{Optimality under condition (iii).} 
Fix $d>1$ 
and define $\Wik(T) \coloneqq (T-1)^d + i \in \Z[T]$, so that $\prodik W_{i, k}(T+1) = \prodik (T^d+i)$ is clearly separable in $\Q[T]$, hence so is $\prodik \Wik(T)$. Let $q \coloneqq Q^d$ for some $Q \le (\log x)^{K_0/d}$ satisfying $P^-(Q) = \ell_0$. Then $1 \in R_k(q)$, showing that $q \in \WUDkAdmSet$. Moreover, $i \in U_q$ for each $i \in [K]$, and any prime $P \equiv 1 \pmod Q$ satisfies $f_i(P^k) = \Wik(P) = (P-1)^d + i \equiv i \pmod{q}$. Consequently, there are $\gg {x^{1/k}}\big/{q^{1/d} \log x}$ many $n \le x$ satisfying $f_i(n) \equiv i \pmod q$ for all $i$, 
and this last expression grows strictly faster than $\phi(q)^{-K} \#\{n \le x: (f(n), q)=1\}$ as soon as $q^{K-1/\Dmin} = q^{K-1/d} \ge (\log x)^{(1+\epsilon) \alpha_k}$ for some fixed $\epsilon \in (0, 1)$. This establishes that the range of $q$ in condition (iii) of Theorem \ref{thm:UnrestrictedInput_SqfrPoly} is optimal, and concrete examples of moduli $q$ satisfying the conditions imposed so far, are those of the form $Q^d$, with $Q$ lying in $[(\log x)^{(1+\epsilon)(K-1/d)^{-1}/d}, (\log x)^{K_0/d}]$ and having least prime factor $\ell_0$.  
\section{Restricted inputs to general moduli: Proof of Theorem \ref{thm:RestrictedInputGenMod1}} \label{sec:Thm_RestrictedInputGenModGen_Proof} 
Fix $T \in \NatNos_{>1}$. We first show that 
as $x \to\infty$ and uniformly in $k$-admissible $\qlelogxKZ$,
\begin{equation}\label{eq:PTnleq Negligible}
\largesum_{\substack{n \le x: ~P_T(n) \le q\\\gcd(f(n), q)=1}} 1 ~ ~ = o\Bigg(\sumnxfnqgcd\Bigg), ~ ~ \largesum_{\substack{n \le x: ~P_T(n_k) \le q\\\gcd(f(n), q)=1}} 1 ~ ~ = o\Bigg(\sumnxfnqgcd\Bigg).
\end{equation}
The first asymptotic is immediate by Proposition \ref{prop:inconvnfnq=1} as $P_{Jk}(n) \le P_T(n)$. 
To show the second, 
we write any $n$ counted in the left side uniquely in the form $n=BN^k A$, where $B$ is $k$-free, $A$ is $(k+1)$-full and the exponent of any prime in $A$ is not a multiple of $k$. Then $n_k = N$, and $B, N, A$ are pairwise coprime, so that $f(n) = f(B) f(N^k) f(A)$, and 
\begin{equation}\label{eq:fnq=1 PTnk le q}
\largesum_{\substack{n \le x: ~P_T(n_k) \le q\\\gcd(f(n), q)=1}} 1 ~ ~ 
\le \largesum_{\substack{B \le x\\B \iskfree\\(f(B), q)=1}} ~ ~ ~\largesum_{\substack{N, A: ~ N^k A \le x/B\\P_T(N) \le q; ~ A\text{ is }(k+1)\text{-full}\\\gcd(f(N^k)f(A), q)=1}} 1.
\end{equation} 
If $A>x^{1/2}$, then $N \le (x/AB)^{1/k} \le x^{1/2k}$. Since $A$ is $(k+1)$-full, the contribution of the tuples $(B, N, A)$ with $A>x^{1/2}$ is $\ll \sum_{B \ll 1} \sum_{N \le x^{1/2k}} (x/BN^k)^{1/(k+1)} \ll x^{1/k - 1/2k(k+1)}$, which is negligible.
On the other hand, if $A \le x^{1/2}$, then given $B$ and $A$, \cite[Lemma 2.3]{PSR22} shows there are $\ll {x^{1/k} (\log_2 x)^T}\big/{B^{1/k} A^{1/k} \log x}$ many $N \le (x/AB)^{1/k}$ having $P_T(N) \le q$. 
The sum over $A$ is $\le \prod_p (1+\sum_{v \ge k+1} p^{-v/k}) \ll 1$, so that the total contribution of all tuples $(B, N, A)$ with $A \le x^{1/2}$ 
is $O(x^{1/k} (\log_2 x)^T/\log x)$. The second formula in \eqref{eq:PTnleq Negligible} now follows from \eqref{prop:fnqcoprimecount}.   

In all of Theorems \ref{thm:RestrictedInputGenMod1} to \ref{thm:HigherPolyControl1}, we may assume $q$ to be sufficiently large, for otherwise these results follow directly from Theorem N and \eqref{eq:PTnleq Negligible}. These formulae also 
show the equality of the second and third expressions in \eqref{eq:WUDPRn} and \eqref{eq:WUDPRnk}, so it remains to show the first equality in either. Recall that for this theorem, we have $\epsilon \coloneqq 1$ and $y = \exp(\sqrt{\log x})$ in the framework developed in section \ref{sec:Convenientn}. Now any convenient $n$ has $P_J(n_k)>y$ and hence is counted in the left hand sides of both \eqref{eq:WUDPRn} and \eqref{eq:WUDPRnk}. By Theorem \ref{thm:ConvenientMAINTerm}, it suffices to show that the contributions of the inconvenient $n$ to the left hand sides of \eqref{eq:WUDPRn} and \eqref{eq:WUDPRnk} are negligible 
compared to $\phi(q)^{-K} \#\{n \le x: (f(n), q)=1\}$. 
In fact, by \eqref{eq:Smoothfnq=1} and \eqref{eq:Large power of prime exceeding y negligible}, it remains to show the bounds \eqref{eq:MainRemainingForRestrInpGenModThm}(i) and (ii) below to establish subparts (a) and (b) of the theorem, respectively:
\begin{align}\label{eq:MainRemainingForRestrInpGenModThm}
\text{(i)} ~ ~ ~  \largesum\nolimits^*_{\substack{n: ~P_R(n)>q}} 1 ~ ~ \ll \xonekphiqKlogxOMalpha, ~ ~ ~ ~ ~ ~ ~ ~ ~  \text{(ii)} ~ ~ ~ \largesum\nolimits^*_{\substack{n: ~ P_{KD+1}(n_k)>q}} 1 ~ ~ \ll \xonekphiqKlogxOMalpha.     
\end{align} 
Here and in the rest of the manuscript, any sum of the form $\sum\nolimits_n^*$ denotes a sum over positive integers $n \le x$ that are not $z$-smooth, not divisible by the $(k+1)$-th power of a prime exceeding $y$, have $P_{Jk}(n)\le y$ and satisfy $\finaimodq$ for all $i \in [K]$. Other conditions imposed on this sum are additional to these.  

Defining $\omegaParalleln \coloneqq \#\{p>q: p^k \parallel n\} = \#\{n \le x: p \parallel n_k\}$ and $\omegaStn \coloneqq \#\{p>q: p^{k+1} \mid n\}$, we first show the following three bounds: 
\begin{equation}\label{eq:OmegaParallelOmegaStReducn}
\sumnolimitsSt_{\substack{n: ~\omegaParalleln \ge KD+1}} ~ 1, ~~ \sumnolimitsSt\mathop{}_{\mkern-15mu\substack{n: ~\omegaParalleln = KD\\\omegaStn \ge 1}} ~  1, ~~  
 \sum_{\substack{n \le x: ~(f(n), q)=1\\\omegaStn \ge Kk, ~ P_{Jk}(n) \le y, ~ P(n)>z\\\nopgrystpkOnedivn}} 1 ~ ~ \ll ~ ~\xonekphiqKlogxOMalpha.
\end{equation}
Any $n$ counted in the first sum is of the form $m (P_{KD+1} \cdots P_1)^k$, where $P_{Jk}(m) \le y$, where $P_1, \dots, P_{KD+1}$ are primes exceeding $q$ satisfying $P_1 \coloneqq P(n) > z$ and $q< P_{KD+1} < \cdots < P_1$, and where $f_i(n) = f_i(m) \prod_{j=1}^{KD+1} f_i(P_j^k) = f_i(m) \prod_{j=1}^{KD+1} \Wik(P_j)$. The conditions $\finaimodq$ can be rewritten as $(P_1, \dots, P_{KD+1}) \bmod q ~ \in ~ \mathcal V_{KD+1, K}^{(k)}\big(q; \aifimfam\big).$ 
Given $m$, $(v_1, \dots, v_{KD+1}) \in \mathcal V_{KD+1, K}^{(k)}\big(q; \aifimfam\big)$, 
and $P_2, \dots, P_{KD+1}$, the number of $P_1$ in $(q$, $x^{1/k}\big/m^{1/k} P_2 \cdots P_{KD+1}]$ satisfying $P_1 \equiv v_1 \pmod q$ is $\ll x^{1/k} \log_2 x\big/m^{1/k} P_2 \cdots P_{KD+1} \phi(q) \log x$, by Brun-Titchmarsh. We sum this over all possible $P_2, \dots, P_{KD+1}$, making use of the bound $\sum_{\substack{q < p \le x \\ p\equiv v\pmod{q}}} 1/p$ $\ll {\log_2{x}}\big/{\phi(q)}$ uniformly in $v \in U_q$ (this follows from Brun--Titchmarsh and partial summation). We deduce that the number of possible $(P_1, \dots, $ $ P_{KD+1})$ satisfying $P_j \equiv v_j \pmod q$ for each $j \in [KD+1]$ is no more than
\begin{equation}\label{eq:P1...PKD+1CountInCongrClasses}
\largesum_{\substack{q<P_{KD+1} < \cdots < P_2 \le x\\\foralljPjvjmodq}} ~ ~ \largesum_{\substack{z < P_1 \le x^{1/k}\big/m^{1/k} P_2 \cdots P_{KD+1}\\P_1 \equiv v_1 \pmod q}} 1 ~ ~ \ll ~ ~ \frac1{\phi(q)^{KD+1}} \cdot \frac{x^{1/k}(\log_2 x)^{O(1)}}{m^{1/k} \log x}.
\end{equation}
Define $V_{r, K}' \coloneqq \max\left\{\#\mathcal V_{r, K}^{(k)}\big(q; \wifam\big): w_1, \dots, w_K \in U_q\right\}$. Summing \eqref{eq:P1...PKD+1CountInCongrClasses} over all $(v_1, \dots,$ $v_{KD+1}) \in \mathcal V_{KD+1, K}^{(k)}\big(q; \aifimfam\big)$ and then over all $m$ via \eqref{eq:mInconvRecip} shows that 
\begin{align}\label{eq:OmegaParallelgeKD+1_BasicBound} \allowdisplaybreaks  
\sumnolimitsSt_{\substack{n: ~\omegaParalleln \ge KD+1}} 1 ~ ~ &\ll \frac{V_{KD+1, K}'}{\phi(q)^{KD+1}} \cdot \frac{x^{1/k}}{(\log x)^{1-\alpha_k/2}} \cdot \expOlogThreexSqPluslogTwoqSq.
\end{align} 
Applying \eqref{eq:VNKkCount_largeN} with $N \coloneqq KD+1$, we get $V_{KD+1, K}'/\phi(q)^{KD+1} \ll \phi(q)^{-K} \prod_{\ell \mid q} (1+O(\ell^{-1/D})) \ll \phi(q)^{-K} \exp\big(O((\log q)^{1-1/D})\big)$. This yields the first bound in \eqref{eq:OmegaParallelOmegaStReducn}.

Next, any $n$ counted in the second sum in \eqref{eq:OmegaParallelOmegaStReducn} can be written in the form $m p^c (P_{KD} \cdots P_1)^k$ for some $m, c$ and distinct primes $p, P_1, \dots, P_{KD}$ exceeding $q$, which satisfy the conditions $P_1 = P(n)>z$, $q<P_{KD}< \cdots < P_1$, $P_{Jk}(m) \le y$, $c \ge k+1$ and $f_i(n) = f_i(m)f_i(p^c) \prod_{j=1}^{KD} \Wik(P_j)$, so that $(P_1, \dots, P_{KD})$ mod $q$ $\in \mathcal V_{KD, K}^{(k)}\big(q; (a_i f_i(mp^c)^{-1})_{i=1}^K\big)$. Given $m, p, c$ and $(v_1, \dots, v_{KD}) \in \mathcal V_{KD, K}^{(k)}\big(q; (a_i f_i(mp^c)^{-1})_{i=1}^K\big)$, the arguments leading to  \eqref{eq:P1...PKD+1CountInCongrClasses} show that the number of possible $(P_1, \dots, P_{KD})$ satisfying $(P_j)_{i=1}^{KD} \equiv (v_j)_{i=1}^{KD} \pmod q$ is $\ll {x^{1/k}(\log_2 x)^{O(1)}}\Big/{\phi(q)^{KD} m^{1/k} p^{c/k} \log x}$. Summing this successively over all $(v_1, \dots, v_{KD})$, $c \ge k+1$, $p>q$ and all possible $m$, shows that the second of the three sums in \eqref{eq:OmegaParallelOmegaStReducn} is $\ll \frac{V_{KD, K}'}{q^{1/k} \phi(q)^{KD}} \cdot \frac{x^{1/k}}{(\log x)^{1-2\alpha_k/3}}$.
(Here we have noted that $\sum_{p>q, ~ c \ge k+1} ~ p^{-c/k} \ll \sum_{p>q} ~ p^{-1-1/k} \ll q^{-1/k}$.) By \eqref{eq:VNKkCount_smallN}, we have ${V_{KD, K}'}\Big/{q^{1/k} \phi(q)^{KD}} \ll 1/{q^K}$, proving the second inequality in \eqref{eq:OmegaParallelOmegaStReducn}.

Lastly, any $n$ counted in the third sum in \eqref{eq:OmegaParallelOmegaStReducn} still has $P(n)>z$ and $P(n)^k \parallel q$, and thus can be written in the form $m p_1^{c_1} \cdots p_{Kk}^{c_{Kk}} P^k$ for some distinct primes $p_1, \dots, p_{Kk}, P$ exceeding $q$ and some integers $m, c_1, \dots, c_{Kk}$, which satisfy $P = P(n)>z$, $P_{Jk}(m) \le y$, $c_j \ge k+1$ for all $j \in [Kk]$, and $\gcd(f(m), q)=1$. Given $m, p_1, \dots, p_{Kk}, c_1, \dots, c_{Kk}$, the number of possible $P>z$ satisfying $P^k \le x/m p_1^{c_1} \cdots p_{Kk}^{c_{Kk}}$ is $\ll x^{1/k}\big/(m p_1^{c_1} \cdots p_{Kk}^{c_{Kk}})^{1/k} \log z$. Summing this over all $c_1, \dots, c_{Kk} \ge k+1$, and then over all $p_1, \dots, p_{Kk}, m$, shows the third bound in \eqref{eq:OmegaParallelOmegaStReducn}.

In what follows, note that $R$ as in the statement of the theorem is the least integer exceeding  
$$\max\left\{k(KD+1)-1, k\left(1+(k+1)\left(K-\frac1D\right)\right)\right\} = \begin{cases}
 k(KD+1)-1, & \text{ if } k<D\\
 k\left(1+(k+1)\left(K-1/D\right)\right) & \text{ if } k \ge D.
\end{cases}$$
\subsection*{Completing the proof of Theorem  \ref{thm:RestrictedInputGenMod1}(a)} Since $q$ is sufficiently large, the $q$-rough part of any $n$ satisfying $\gcd(f(n), q)=1$ is $k$-full (by Lemma \ref{lem:kfreepartbdd}). As such, any $n$ with $\omegaStn = 0$ counted in \eqref{eq:MainRemainingForRestrInpGenModThm}(i) must have $\omegaParalleln \ge \lfloor R/k \rfloor \ge KD+1$, and hence is counted in the first sum in \eqref{eq:OmegaParallelOmegaStReducn}. Moreover, any $n$ with $\omegaParalleln = KD$ counted in \eqref{eq:MainRemainingForRestrInpGenModThm}(i) must also have $\omegaStn \ge R - k \omegaParalleln \ge k(KD+1) - kKD \ge 1$, and hence is counted in the second sum in \eqref{eq:OmegaParallelOmegaStReducn}. By \eqref{eq:OmegaParallelOmegaStReducn}, it thus remains to show that the contribution of $n$ having $\omegaParalleln \in [KD-1]$ and $\omegaStn \in [Kk-1]$ to the left hand side of \eqref{eq:MainRemainingForRestrInpGenModThm} is absorbed in the right hand side. This would follow once we show that for any fixed $r \in [KD-1]$ and $s \in [Kk-1]$, the contribution $\Sigma_{r, s}$ of all $n$ with $\omegaParalleln = r$ and $\omegaStn = s$ to the left hand side of \eqref{eq:MainRemainingForRestrInpGenModThm}(i) is absorbed in the right hand side. 

Now any $n$ counted in $\Sigma_{r, s}$ is of the form $m p_1^{c_1} \cdots p_s^{c_s} P_1^k \cdots P_r^k$ for some distinct primes $p_1, \dots, p_s,$ $ P_1,$ $\dots, P_r$ and integers $m, c_1, \dots, c_s$, which satisfy the following conditions: \textbf{(i)} $P(m) \le q$; \textbf{(ii)} $P_1 \coloneqq P(n)>z$; $q<P_r < \cdots < P_1$; \textbf{(iii)} $p_1, \dots, p_s>q$; \textbf{(iv)} $c_1, \dots, c_s \ge k+1$ and $c_1 + \cdots + c_s \ge R-kr$; 
\textbf{(v)} $m$, $p_1, \dots, p_s, P_1, \dots, P_r$ are all pairwise coprime, so that $f_i(n) = f_i(m) f(p_1^{c_1}) \cdots f(p_s^{c_s}) \prod_{j=1}^r \Wik(P_j)$ for each $i \in [K]$. Here, property (i) holds because the $q$-rough part of any $n$ satisfying $\gcd(f(n), q)=1$ is $k$-full, whereas $\omegaParalleln = r$, $\omegaStn = s$ . 

With $\tau_i \coloneqq \min\{c_i, R-kr\}$, it is easy to see that the integers $\tau_1, \dots, \tau_s \in [k+1, R-kr]$ satisfy $\tau_1 \le c_1, \dots, \tau_s \le c_s$ and $\tau_1 + \cdots + \tau_s \ge R-kr$. (Here it is important that $R \ge k(KD+1)$, $r \le KD-1$ and $c_1 + \cdots + c_s \ge R-kr$.) 
Turning this around, we find that 
\begin{equation}\label{eq:Sigma_rs_FirstBound_GenModGen(a)}
\Sigma_{r, s} \le \sum_{\substack{\tau_1, \dots, \tau_s \in [k+1, R-kr]\\\tau_1+ \cdots + \tau_s \ge R-kr}} ~ ~ \mathcal N_{r, s}(\tau_1, \dots, \tau_s),
\end{equation} 
where $\mathcal N_{r, s}(\tau_1, \dots, \tau_s)$ denotes the contribution of all $n$ counted in \eqref{eq:MainRemainingForRestrInpGenModThm}(i) which can be written in the form $m p_1^{c_1} \cdots p_s^{c_s} P_1^k \cdots P_r^k$ for some distinct primes $p_1, \dots, p_s, P_1,$ $\cdots, P_r$ and integers $m, c_1, \dots, c_s$ satisfying the conditions (i)-(v) above, along with the condition $c_1 \ge \tau_1, \dots, c_s \ge \tau_s$. We will show that for each tuple $(\tau_1, \dots, \tau_s)$ occurring in \eqref{eq:Sigma_rs_FirstBound_GenModGen(a)}, we have 
\begin{equation}\label{eq:NrsIneq_GenModGen}
\mathcal N_{r, s}(\tau_1, \dots, \tau_s) ~ ~ \ll ~ ~  \xonekqKloglogxlogx.
\end{equation}
Consider an arbitrary such tuple $(\tau_1, \dots, \tau_s)$, and write $n$ in the form $m p_1^{c_1} \cdots p_s^{c_s} P_1^k \cdots P_r^k$ as above. The conditions $\finaimodq$ lead to $(P_1, \dots, P_r) \bmod q \in$ $\mathcal V_{r, K}^{(k)}\big(q; (a_i f_i(mp_1^{c_1} \cdots $ $ p_s^{c_s})^{-1})_{i=1}^K\big).$ Given $m, p_1, \dots, p_s, c_1, \dots, c_s$ and $(v_1, \dots, v_r)$ $\in \mathcal V_{r, K}^{(k)}\big(q; (a_i f_i(mp_1^{c_1} \cdots p_s^{c_s})^{-1})_{i=1}^K\big)$, the arguments leading to
\eqref{eq:P1...PKD+1CountInCongrClasses} show that the number of possible $P_1, \dots, P_r$ satisfying $P_j \equiv v_j$ mod $q$ for each $j \in [r]$, is 
$\ll {x^{1/k} (\log_2 x)^{O(1)}}\Big/{\phi(q)^r m^{1/k} p_1^{c_1/k} \cdots p_s^{c_s/k} \log x}.$
With $V_{r, K}' = \max_{(w_i)_i \in U_q^K} \#\mathcal V_{r, K}^{(k)}\big(q; \wifam\big)$ as before, 
the bounds $\sum_{p_i>q: ~ c_i \ge \tau_i} p_i^{-c_i/k} \ll q^{-(\tau_i/k-1)}$ yield 
\begin{equation}\label{eq:NrsGeneralBoundOld}
\mathcal N_{r, s}(\tau_1, \dots, \tau_s) \ll \frac1{q^{(\tau_1+\cdots+\tau_s)/k-s}} ~ \frac{V_{r, K}'}{\phi(q)^r} \cdot \frac{x^{1/k} (\log_2 x)^{O(1)}}{\log x} \largesummqSmfmqOne. 
\end{equation} 
Proceeding as in the argument for \eqref{eq:mInconvRecip}, we write any $m$ in the above sum as $BM$ where $B$ is $k$-free and $M$ is $k$-full, so that $B = O(1)$ and $P(M) \le q$. We find that 
\begin{equation}\label{eq:mP(m)leq}
\largesummqSmfmqOne \ll \largesum_{\substack{M \le x: ~ P(M) \le q\\M\iskfull}} \frac1{M^{1/k}} \le \prod_{p \le q} \left(1+\frac1p + O\left(\frac1{p^{1+1/k}}\right)\right) \ll \exp\left(\largesum_{p \le q} \frac1p\right) \ll \log q.    
\end{equation} 
Inserting this into \eqref{eq:NrsGeneralBoundOld}, we obtain
\begin{equation}\label{eq:NrsGeneralBound}
\mathcal N_{r, s}(\tau_1, \dots, \tau_s) \ll \frac1{q^{(\tau_1+\cdots+\tau_s)/k-s}} ~ \frac{V_{r, K}'}{\phi(q)^r} \cdot \frac{x^{1/k} (\log_2 x)^{O(1)}}{\log x}. 
\end{equation}
Now since $1 \le r \le KD-1$, an application of \eqref{eq:VNKkCount_smallN} with $N \coloneqq r$ now yields 
\begin{equation}\label{eq:Nrstau1...taus FinalBound}
\mathcal N_{r, s}(\tau_1, \dots, \tau_s) \ll \frac{\exp\big(O(\omega(q))\big)}{q^{(\tau_1+\cdots+\tau_s)/k-s+r/D}} \cdot \xkonekloglogxlogx \ll \frac{\exp\big(O(\omega(q))\big)}{q^{\max\{s/k, R/k-r-s\} +r/D}} \cdot \xkonekloglogxlogx,
\end{equation}
where in the last equality we have recalled that $\tau_1, \dots, \tau_s \ge k+1$ and $\tau_1 + \cdots + \tau_s \ge R-kr$. We claim that $\max\{s/k, R/k-r-s\} +r/D>K$. This is tautological if $s/k+r/D > K$, so suppose $s/k+r/D \le K$. Then $r \le D(K-s/k) \le DK-D/k$, and $s \le k(K-r/D)$ so that $R/k-r-s+r/D \ge R/k - Kk + ((k+1)/D-1)r$. If $k<D$, then $(k+1)/D-1 \le 0$, so for all $1 \le r \le DK-D/k$, we have $R/k - Kk + ((k+1)/D-1)r \ge R/k - Kk + ((k+1)/D-1)(DK-D/k)$ and this exceeds $K$ since $R \ge k(KD+1)$. 
If on the other hand, we had $k \ge D$, then $k+1>D$ and the minimum value of $R/k - Kk + ((k+1)/D-1)r$ is attained at $r=1$, giving us $R/k - Kk + ((k+1)/D-1)r \ge R/k - Kk + ((k+1)/D-1)$ which also exceeds $K$ since $R> k\big(1+(1+k)(K-1/D)\big)$. This shows our claim, 
so that \eqref{eq:Nrstau1...taus FinalBound} leads to \eqref{eq:NrsIneq_GenModGen}. Summing \eqref{eq:NrsIneq_GenModGen} over the $O(1)$ many possible tuples $(\tau_1, \dots, \tau_s)$ occurring in the right hand side of \eqref{eq:Sigma_rs_FirstBound_GenModGen(a)} yields $\Sigma_{r, s} \ll x^{1/k} (\log_2 x)^{O(1)}\big/q^K \log x$, which (as argued before) establishes Theorem \ref{thm:RestrictedInputGenMod1}(a). 

\subsection*{Completing the proof of Theorem \ref{thm:RestrictedInputGenMod1}(b):} 
Define $\omega_k(n) \coloneqq \#\{p>q: p^2\mid n_k\}$. Any $n$ with $\omega_k(n) = 0$ counted in \eqref{eq:MainRemainingForRestrInpGenModThm}(ii)  also has $\omegaParalleln \ge KD+1$ (since $P_{KD+1}(n_k)>q$), and hence any such $n$ is counted in the first sum in \eqref{eq:OmegaParallelOmegaStReducn}. Likewise, any $n$ with $\omegaParalleln = KD$ counted in \eqref{eq:MainRemainingForRestrInpGenModThm}(ii) has $\sum_{p>q: ~ v_p(n_k)>1} ~ v_p(n_k) \ge (KD+1)-\omegaParalleln \ge 1$, so that any such $n$ also has $\omegaStn \ge \omega_k(n) \ge 1$ and is counted in the second sum in \eqref{eq:OmegaParallelOmegaStReducn}. 
It thus remains to show that for each $r \in [KD-1]$ and $s \in [Kk-1]$, the contribution $\Sigmatilrs$ of all $n$ with $\omegaParalleln = r$ and $\omegakn = s$ to the left hand side of \eqref{eq:MainRemainingForRestrInpGenModThm}(ii) is absorbed in the right. 

Any $n$ counted in $\Sigmatilrs$ has $n_k$ of the form $m' p_1^{c_1} \cdots p_s^{c_s} P_1 \cdots P_r$ for some distinct primes $p_1, \dots ,p_s,$ $P_1, \dots, P_r$ and integers $m', c_1, \dots, c_s$, which satisfy conditions (i)--(v): \textbf{(i)} $P(m') \le q$; \textbf{(ii)} $P_1 \coloneqq P(n_k) = P(n)>z$, $q<P_r < \cdots < P_1$; \textbf{(iii)} $p_1, \dots, p_s>q$;
\textbf{(iv)} $c_1, \dots, c_s \ge 2$ and $c_1 + \cdots + c_s \ge KD+1-r$; \textbf{(v)} $m'$, $p_1, \dots, p_s, P_1, \dots, P_r$ are all pairwise coprime. Hence, $n$ is of the form $m  p_1^{c_1 k} \cdots p_s^{c_s k} P_1^k \cdots P_r^k$, where $p_1, \dots, p_s, P_1, \dots, P_r$ are as above, and: \textbf{(vi)} $P_{Jk}(m) \le y$; \textbf{(vii)} $f_i(n) = f_i(m) f_i(p_1^{c_1 k}) \cdots f_i(p_s^{c_s k}) \prod_{j=1}^r \Wik(P_j)$ for each $i \in [K]$.

Now since $r \le KD-1$, the integers $\tau_j \coloneqq \min\{c_j, KD+1-r\}$ ($j \in [s]$) 
satisfy $\tau_1, \dots, \tau_s \in [2, KD+1-r]$ and $\tau_1 + \cdots + \tau_s \ge KD+1-r$. We now obtain the following analogue of \eqref{eq:Sigma_rs_FirstBound_GenModGen(a)}
\begin{equation}\label{eq:Sigmatil_rs_FirstBound_GenModGen(b)}
\Sigmatilrs \le \sum_{\substack{\tau_1, \dots, \tau_s \in [2, KD+1-r]\\\tau_1+ \cdots + \tau_s \ge KD+1-r}} ~ ~ \NtiltauOnetaus,
\end{equation} 
where $\NtiltauOnetaus$ denotes the number of $n$ which can be written in the form $m  p_1^{c_1 k} \cdots p_s^{c_s k} P_1^k$ $\cdots P_r^k$ with $m,  p_1, \dots, p_s, c_1, \dots, c_s, P_1, \dots, P_r$ satisfying the conditions (ii), (iii), (vi), (vii) above, and with $c_1 \ge \tau_1, \dots, c_s \ge \tau_s$. 
We show that for each $\tau_1, \dots, \tau_s$ counted above, 
\begin{equation}\label{eq:NrsTilIneq_GenModGen}
\NtiltauOnetaus ~~ \ll ~~ \xonekphiqKlogxOMalpha.   
\end{equation}
The argument is analogous to that given for \eqref{eq:NrsIneq_GenModGen}, so we only sketch it. We write any $n$ counted in $\Sigmatilrs$ in the form $m  p_1^{c_1 k} \cdots p_s^{c_s k} P_1^k$ $\cdots P_r^k$, with $m,  p_1, \dots, p_s, c_1, \dots, c_s, P_1, \dots, P_r$ satisfying the conditions (ii), (iii), (vi), (vii) above, and with $c_1 \ge \tau_1, \dots, c_s \ge \tau_s$, so that $(P_1, \dots, P_r)$ mod $q$ $\in \mathcal V_{r, K}^{(k)}\big(q; (a_i f_i(m p_1^{c_1 k} \cdots p_s^{c_s k})^{-1})_{i=1}^K\big)$. Thus, given $m,  p_1, \dots, p_s, c_1, \dots, c_s$ and $(v_1, \dots, v_r)\in \mathcal V_{r, K}^{(k)}\big(q; (a_i f_i(m p_1^{c_1 k} \cdots p_s^{c_s k})^{-1})_{i=1}^K\big)$, the number of possible $P_1, \dots, P_r$ satisfying $P_j \equiv v_j \pmod q$ for each $j \in [r]$, is $\ll {x^{1/k} (\log_2 x)^{O(1)}}\big/{\phi(q)^r m^{1/k} p_1^{c_1} \cdots p_s^{c_s} \log x}.$ Hence
\begin{equation}\label{eq:NrstilGeneralBound}
\NtiltauOnetaus \ll \frac1{q^{\tau_1+\cdots+\tau_s-s}} ~ \frac{V_{r, K}'}{\phi(q)^r} \cdot \frac{x^{1/k}}{(\log x)^{1-\alpha_k/2}} \expOlogThreexSqPluslogTwoqSq.   
\end{equation}
Now applying \eqref{eq:VNKkCount_smallN} and using the fact that $\tau_1+ \cdots + \tau_s \ge \max\{2s, KD+1-r\}$, we find that ${V_{r, K}'}\big/{\phi(q)^r q^{\tau_1+\cdots+\tau_s-s}} \ll {\exp\big(O(\omega(q))\big)}\big/{q^{\max\{s, KD+1-r-s\}+r/D}} \ll \phi(q)^{-K}$, since from $D \ge 2$, it is easily seen that $\max\{s, KD+1-r-s\}+r/D > K$. This establishes \eqref{eq:NrsTilIneq_GenModGen}, so that \eqref{eq:Sigmatil_rs_FirstBound_GenModGen(b)} yields $\Sigmatilrs \ll {x^{1/k}}\big/{\phi(q)^K (\log x)^{1-2\alpha_k/3}}$, completing the proof of Theorem \ref{thm:RestrictedInputGenMod1}(b). \hfill \qedsymbol
\section{Final preparatory step for Theorem \ref{thm:RestrictedInputSqfreeMod1}: Counting points on varieties}\label{sec:FinalAlgStep} 
To establish Theorem \ref{thm:RestrictedInputSqfreeMod1}, we will need the following partial improvements of 
Corollary \ref{cor:VNKCountSqfrMegaGen}. 
In this section, we again deviate from the general notation set up for Theorems \ref{thm:UnrestrictedInput_SqfrPoly} to \ref{thm:HigherPolyControl1}, so the notation set up in this section will be relevant in this section only. 
\begin{prop}\label{prop:V21V3Count_Sqfreeq}
Let $F \in \Z[T]$ be a fixed nonconstant polynomial which is not squarefull.
\begin{enumerate}
\item[\textbf{(a)}] Define $\VTwoOne \coloneqq \{(v_1, v_2) \in U_\ell^2: F(v_1) F(v_2) \equiv w \pmod \ell\}$. 
Then $\#\VTwoOne \le \phi(\ell) \left(1 + O\left({\ell^{-1/2}}\right)\right)$, uniformly for primes $\ell$ and coprime residues $w$ mod $\ell$. 
\item[\textbf{(b)}] Let $G \in \Z[T]$ be any fixed polynomial such that $\{F, G\} \subZT$ are multiplicatively independent. Let $
\VThreeTwo$ be the set of $(v_1, v_2, v_3) \in U_\ell^3$ satisfying the two congruences $F(v_1) F(v_2) F(v_3) \equiv u  \pmod \ell$ and $~ G(v_1) G(v_2) G(v_3) \equiv w \pmod \ell$. 
Then $\#\VThreeTwo \ll_{F, G} \phi(\ell)$, uniformly in primes $\ell$ and coprime residues $u, w$ mod $\ell$.
\end{enumerate}
\end{prop} 
Our starting idea will be to look at $\VTwoOne$ and $\VThreeTwo$ as subsets of the sets of $\F_\ell$-rational points of certain varieties over the algebraic closure $\Fellbar$ of $\F_\ell$. 
\begin{prop}\label{prop:Variety_PointCounting} 
Let $V$ be a variety defined over $\F_\ell$ and $V(\F_\ell) \coloneqq V \cap \F_\ell$. 
\begin{enumerate} 
\item[\textbf{(a)}] 
If $V$ is an absolutely irreducible affine plane curve, 
then $\#V(\F_\ell) \le \ell+O(\sqrt{\ell})$, where the implied constant depends only on the degree of $V$.  
\item[\textbf{(b)}] Let $d$ be the positive integer such that $V \subset (\overline{\F}_\ell)^d$. 
We have $\#V(\F_\ell) \ll \ell^{\dim V}$, where $\dim V$ is the dimension of $V$ as a variety, and the implied constant depends at most on $d$ and on the number and degrees of the polynomials defining $V$. 
\end{enumerate}
\end{prop} 
Subpart(a) is a consequence of \cite[Corollary 2b]{LY94}, while subpart (b) is a weaker version of \cite[Claim 7.2]{DKL12}  
but in fact goes back to work of Lang and Weil \cite[Lemma 1]{LangWeil}.
To make use of the aforementioned results, we will also be needing the following observations. 
\begin{lem}\label{lem:Irredicibility and NonDivisibility}
Let $F, G \in \Z[T]$ be fixed multiplicatively independent polynomials such that $F$ is not squarefull. There exist constants $\kappa_0(F)$ and $\kappa_1(F, G)$ such that:
\begin{enumerate}
\item[\textbf{(a)}] For any $N \ge 2$, $\ell > \kappa_0(F)$ and $w \in \F_\ell^\times$, the polynomial $\prod_{i=1}^N F(X_i) - w$ is absolutely irreducible over $\F_\ell$, that is, it is irreducible in the ring $\Fellbar[X_1, \dots, X_N]$.
\item[\textbf{(b)}] For any $\ell>\kappa_1(F, G)$ and $u, w \in \F_\ell^\times$, the polynomial $F(X) F(Y) F(Z) - u$ is irreducible and doesn't divide the polynomial $G(X) G(Y) G(Z) - w$ in the ring $\Fellbar[X, Y, Z]$. 
\end{enumerate}
\end{lem}
\begin{proof} 
Write $F \coloneqq r\prod_{j=1}^M G_j^{b_j}$ for some $r \in \Z$, $b_j \in \NatNos$, and pairwise coprime irreducibles $G_j \in \Z[T]$, so that by the nonsquarefullness of $F$ in $\Z[T]$, we have $b_j = 1$ for some $j \in [M]$. By the observations at the start of the proof of Proposition \ref{prop:OrdDerivInfo}, there exists a constant $\kappa_0(F)$ such that for any prime $\ell > \kappa_0(F)$, $\ell$ doesn't divide the leading coefficient of $F$ 
and $\prod_{j=1}^M G_j$ is separable in $\F_\ell[T]$. This forces $\prod_{\substack{\theta \in \Fellbar\\F(\theta) = 0}} (T-\theta)^2 \nmid F(T)$ in $\Fellbar[T]$. 

\textit{Proof of (a).} 
We will show that for any $\ell>\kappa_0(F)$ and $U, V \in \Fellbar[X_1, \dots, X_N]$ satisfying
\begin{equation}\label{eq:ProdFi-w EqualsUV}
\prod_{i=1}^N F(X_i) - w = U(X_1, \dots, X_N) V(X_1, \dots, X_N),   
\end{equation}
one of $U$ or $V$ must be constant. First note that for any root $\theta \in \Fellbar$ of $F$, we have $-w = U(X_1, \dots, X_{N-1}, \theta) V(X_1, \dots,$ $X_{N-1}, \theta)$, forcing $U(X_1, \dots, X_{N-1}, \theta)$ and $V(X_1, \dots, X_{N-1}, \theta)$ to be constant in the ring $\Fellbar[X_1, \dots, X_N]$. Writing $U(X_1, \dots, X_N)$, $V(X_1, \dots, X_N)$ as
$$\sum_{\substack{i_1, \dots, i_{N-1} \ge 0\\i_1 \le R_1, \dots, i_{N-1} \le R_{N-1}}} u_{i_1, \dots, i_{N-1}}(X_N) ~ X_1^{i_1} \cdots X_{N-1}^{i_{N-1}}, ~ ~ \sum_{\substack{j_1, \dots, j_{N-1} \ge 0\\j_1 \le T_1, \dots, j_{N-1} \le T_{N-1}}} v_{j_1, \dots, j_{N-1}}(X_N) ~ X_1^{j_1} \cdots X_{N-1}^{j_{N-1}}$$ 
respectively (where $u_{i_1, \dots, i_{N-1}}, v_{j_1, \dots, j_{N-1}} \in \Fellbar[X_N]$ and neither $u_{R_1, \dots, R_{N-1}}$ nor $v_{T_1, \dots, T_{N-1}}$ is identically zero), we thus find that $u_{i_1, \dots, i_{N-1}}(\theta) = v_{j_1, \dots, j_{N-1}}(\theta) = 0$ for any $(i_1, \dots, i_{N-1}) \ne (0, \dots, 0)$, $(j_1, \dots, j_{N-1}) \ne (0, \dots, 0)$, and any $\theta$ as above. Thus, if the tuples $(R_1, \dots, R_{N-1})$ and $(T_1, \dots, T_{N-1})$ are both nonzero, then $\prod_{\substack{\theta \in \Fellbar\\F(\theta) = 0}} (X_N-\theta)$ divides $u_{R_1, \dots, R_{N-1}}(X_N)$ and $v_{T_1, \dots, T_{N-1}}(X_N)$ in $\Fellbar[X_N]$. But then, if $\alpha \in \Z$ is the leading coefficient of $F$, then comparing the monomials (in $X_1, \dots, X_{N-1}$) with maximal total degree in \eqref{eq:ProdFi-w EqualsUV}, we find that $\alpha^{N-1} F(X_N) = u_{R_1, \dots, R_{N-1}}(X_N) ~ v_{T_1, \dots, T_{N-1}}(X_N) 
\equiv 0 \pmod{\prod_{\substack{\theta \in \Fellbar\\F(\theta) = 0}} (X_N-\theta)^2}$, which is impossible by the observations in the first paragraph of the proof. 
This forces one of $(R_1, \dots, R_{N-1})$ or $(T_1, \dots, T_{N-1})$ to be $(0, \dots, 0)$, say the latter. Then $V(X_1, \dots, X_N) = v_{0, \dots, 0}(X_N)$ and 
since $N \ge 2$, plugging $X_1 \coloneqq \theta$ for some root $\theta \in \Fellbar$ of $F$ into \eqref{eq:ProdFi-w EqualsUV} yields $-w = U(\theta, X_2, \dots, X_N) v_{0, \dots, 0}(X_N)$, forcing $V$ to be identically constant. 

\textit{Proof of (b).} 
We claim that for all primes $\ell \gg_{F, G} 1$, if the rational function $F^a G^b$ is constant in the ring $\Fellbar(T)$ for some integers $a, b$, then $a \equiv b \equiv 0 \pmod\ell$.\footnote{It is not difficult to see that this also forces $a=b=0$, but we won't need that.} The argument for this is a simple variant of that given for the inequality ``$\ord_\ell(\Ftil) \le \bbm_{\ell \le C_1} C_1$" in the proof of Proposition \ref{prop:OrdDerivInfo}(b), so we only sketch the outline. Since $\{F, G\} \subset \Z[T]$ are multiplicatively independent, the polynomials $\{F' G, F G'\} \subset \Z[T]$ are $\Q$-linearly independent, hence so are the columns of the matrix $M_1$ listing the coefficients of $F'G$ and $FG'$ in two columns. Hence we can find invertible matrices $P_1$ 
and $Q_1$ (where $Q_1$ is a $2 \times 2$ matrix) such that $P_1 M_1 Q_1  = \diag(\beta_1, \beta_2)$ for some $\beta_1, \beta_2 \in \Z\sm\{0\}$ satisfying $\beta_1 \mid \beta_2$. Let $\ell>|\beta_2|$ be any prime not dividing the leading coefficients of $F$, $G$, $F' G$ or $F G'$. If $F^a G^b$ is identically constant in $\F_\ell[T]$, then $a F' G + b F G' \equiv 0$ in $\FellT$, so $M_1 (a ~ b)^\top \equiv 0 \pmod\ell$. Hereafter, familiar calculations yield $(a ~ b)^\top \equiv 0 \pmod\ell$. 

Collecting our observations, we have shown that there exists a constant $\kappa_1(F, G)$ such that for all primes $\ell>\kappa_1(F, G)$, the following three properties hold:
\begin{enumerate}
    \item[(i)] $\ell > \kappa_0(F)$, so that 
$\prod_{\substack{\theta \in \Fellbar\\F(\theta) = 0}} (T-\theta)^2 \nmid F(T)$ in $\Fellbar[T]$; \item[(ii)] $\ell$ doesn't divide the leading coefficient of $F$ or $G$; and, \item[(iii)] For any $a, b \in \Z$ for which $F^a G^b$ is identically constant in $\Fellbar(T)$, we have $\ell \mid a$ and $\ell \mid b$.
\end{enumerate}

We will now show that any such constant $\kappa_1(F, G)$ satisfies the property in subpart (b) of the lemma. By subpart (a), $\FXFYFZminusu$ is already irreducible in $\Fellbar[X, Y, Z]$ for any $u \in \F_\ell^\times$. Assume by way of contradiction that for some $\ell>\kappa_1(F, G)$ and $u, w \in \F_\ell^\times$, we have 
\begin{equation}\label{eq:DivisibilityIdentity}
G(X) G(Y) G(Z) - w = H_0(X, Y, Z) ~ (\FXFYFZminusu) ~ \text{ for some }H_0 \in \Fellbar[X, Y, Z]. 
\end{equation}
Write $H_0(X, Y, Z) \eqqcolon \sum_{\substack{0 \le i_1 \le r_1\\0 \le i_2 \le r_2}} h_{i_1, i_2}(X) Y^{i_1} Z^{i_2}$ for some $h_{i_1, i_2} \in \Fellbar[X]$ with $h_{r_1, r_2}$ not identically zero. If $(r_1, r_2) = (0, 0)$, 
then substituting a root of $F$ and $G$ in place of $Y$ and $Z$ respectively, we see that $H_0$ must be a constant $\lambda_0 \in \Fellbar\sm\{0\}$ satisfying $w = \lambda_0 u$. Thus  
$G(X) G(Y) G(Z) = \lambda_0 F(X) F(Y) F(Z)$. Now substituting some $\eta \in \Fellbar$ which is not a root of $FG$ in place of both $Y$ and $Z$ leads to $F(X) G(X)^{-1} = \lambda_0^{-1} F(\eta)^{-2} G(\eta)^2$, a nonzero constant. But since $(1, -1) \not\equiv (0, 0) \pmod\ell$, this violates condition (iii) in the definition of $\kappa_1(F, G)$. Hence $(r_1, r_2) \ne (0, 0)$. 

Let $\alpha, \beta \in \Z$ denote the leading coefficients of $F$ and $G$ respectively. Comparing the monomials in $Y$ and $Z$ of maximal total degree 
in \eqref{eq:DivisibilityIdentity} yields $\beta^2 G(X) = \alpha^2 F(X) h_{r_1, r_2}(X)$ in $\Fellbar[X]$, so that (since either side of this identity is nonzero), we get $F \mid G$ in $\Fellbar[X]$. Write $G = F^m H$ for some $m \ge 1$ and $H \in \Fellbar[X]$ such that $F \nmid H$ in $\Fellbar[X]$. An easy finite induction shows that with $G_t(X, Y, Z) \coloneqq F(X)^{m-t} F(Y)^{m-t}F(Z)^{m-t} H(X) H(Y) H(Z) - u^{-t}w$ and $\FhatXYZ \coloneqq \FXFYFZminusu$, we have $\Fhat \mid G_t$ for each $t \in \{0, 1, \dots, m\}$. Indeed, the case $t=0$ is just \eqref{eq:DivisibilityIdentity}, and if $\Fhat \mid G_t$ for some $t \le m-1$, then writing $G_t = Q_t\Fhat$ shows that $F(X)F(Y)F(Z)$ $\mid (Q_t(X, Y, Z) - u^{-(t+1)}w)$. With $Q_{t+1}$ defined by $Q_t(X, Y, Z) - u^{-(t+1)}w = F(X)F(Y)F(Z) Q_{t+1}(X, Y, Z)$, we obtain $G_{t+1} = Q_{t+1}\Fhat$ completing the induction. 

Applying this last observation with $t \coloneqq m$ shows that $\Fhat(X, Y, Z)$ divides $H(X) H(Y) H(Z) - u^{-m} w$ in $\Fellbar[X, Y, Z]$. We claim that this forces $H$ to be constant. Indeed if not, then letting $\gamma \in \Fellbar\sm\{0\}$ be the leading coefficient of $H$, \footnote{Here $\gamma \ne 0$ in $\Fellbar$ because $\ell$ doesn't divide the leading coefficient of $G = F^m H$.} writing $H(X) H(Y) H(Z) - u^{-m} w = (\FXFYFZminusu) ~ \sum_{\substack{0 \le i_1 \le b_1\\0 \le i_2 \le b_2}} g_{i_1, i_2}(X) Y^{i_1} Z^{i_2}$
for some $g_{i_1, i_2} \in \Fellbar[X]$ with $g_{b_1, b_2} \ne 0$, and comparing the monomials in $Y$ and $Z$ of maximal degree, we obtain $\gamma^2 H(X) = \alpha^2 F(X) g_{b_1, b_2}(X)$. This leads to $F \mid H$, contrary to hypothesis. Hence $H$ must be constant, so the identity $F^{-m} G = H$ in $\Fellbar(X)$ violates condition (iii) in the definition of $\kappa_1(F, G)$, as $(-m, 1) \not\equiv (0, 0) \pmod\ell$. This shows that $\Fhat$ cannot divide $G(X)G(Y)G(Z)-w$, completing the proof. 
\end{proof}
Given a commutative ring $R$ and an $R$-module $M$, we say that $x \in R$ is an \textsf{$M$-regular element} if $x$ is not a zero-divisor on $M$, that is, if $xz = 0$ for some $z \in M$ implies $z=0$. A sequence $x_1, \dots, x_n$ of elements of $R$ is said to be \textsf{$M$-regular} if  
$x_1$ is an $M$-regular element, 
each $x_i$ is an $M/(x_1, \dots, x_{i-1})M$-regular element, \textit{and}  
$M/(x_1, \dots, x_n) M \ne 0$. It is well-known (see \cite[Proposition 1.2.14]{BH98}) that for any proper ideal $I$ in a Noetherian ring $R$, the height of $I$ is at least the length of the longest $R$-regular sequence contained in $I$. 
\begin{proof}[Proof of Proposition \ref{prop:V21V3Count_Sqfreeq}.] With $\kappa_0(F)$ and $\kappa_1(F, G)$ as in Lemma \ref{lem:Irredicibility and NonDivisibility}, 
the affine plane curve $\{(X, Y) \in \Fellbar^2: ~ F(X) F(Y) - w = 0\}$  
is absolutely irreducible for any $\ell > \kappa_0(F)$, so that Proposition \ref{prop:Variety_PointCounting}(a) yields Proposition \ref{prop:V21V3Count_Sqfreeq}(a).  
For (b), it suffices to show that for any prime $\ell > \kappa_1(F, G)$, the variety 
$V_\ell \subset \Fellbar^3$ defined by the polynomials $\FhatXYZ \coloneqq F(X) F(Y) F(Z) - u$ and $\GhatXYZ \coloneqq G(X) G(Y) G(Z) - w$ has $\ll_{F, G} \ell$ many $\F_\ell$-rational points. 
Consider the ideal $I(V_\ell)$ of the ring $R \coloneqq \Fellbar[X, Y, Z]$ consisting of all polynomials vanishing at all the points of $V_\ell$, so that $(\Fhat, \Ghat) \subset I(V_\ell)$. 
If $I(V_\ell) = R$, then $V_\ell = \emptyset$, so suppose $I(V_\ell) \subsetneq R$. Lemma \ref{lem:Irredicibility and NonDivisibility}(b) shows that the sequence $\Ghat, \Fhat \in I(V_\ell)$ is $R$-regular, so by \cite[Proposition 1.2.14]{BH98}, $I(V_\ell)$ has height at least $2$. By \cite[Chapter 11, Exercise 7]{AM69}, the Krull-dimension of $R$ is $3$, whence that of $R/I(V_\ell)$ is at most $3-2=1$ (by, say, \cite[p. 31]{Matsu06}). Thus $\dim(V_\ell) \le 1$,  and Proposition \ref{prop:Variety_PointCounting} completes the proof.
\end{proof}
\section{Restricted inputs to squarefree moduli: Proof of Theorem \ref{thm:RestrictedInputSqfreeMod1}}\label{sec:ResInpSqfreePf}
Returning to the notation set up in the introduction, we start with the same initial reductions as in section \ref{sec:Thm_RestrictedInputGenModGen_Proof}. 
As such, to establish subpart (a) of the theorem, it suffices to show the bound (i) below with $k \ge 2$ and with the respective values of $R$ defined in the statement, -- while in order to establish subpart (b), it suffices to show (ii) below, with the $2K+1$ replaced by $2$ in the case when $K=1$ and $W_k = W_{1, k}$ is not squarefull:
\begin{align}\label{eq:MainRemainingForRestrInpSqfreeModThm}
\text{(i)} ~ ~ ~  \largesum\nolimits^*_{\substack{n: ~P_R(n)>q}} 1 ~ ~ \ll \xonekphiqKlogxOMalpha, ~ ~ ~ ~ ~ ~ ~ ~ ~  \text{(ii)} ~ ~ ~ \largesum\nolimits^*_{\substack{n: ~ P_{2K+1}(n_k)>q}} 1 ~ ~ \ll \xonekphiqKlogxOMalpha.     
\end{align} 
 Here 
we again have $\epsilon = 1$ and $y = \exp(\sqrt{\log x})$ in the framework developed in section \ref{sec:Convenientn}. We will also retain the notation $\omegaParalleln = \#\{p>q: p^k \parallel n\} = \#\{n \le x: p \parallel n_k\}$, $\omegaStn = \#\{p>q: p^{k+1} \mid n\}$, and $\omegakn = \#\{p>q: p^2 \mid n_k\}$ from section \ref{sec:Thm_RestrictedInputGenModGen_Proof}.

For technical reasons, we first give a separate proof of 
all the above bounds in the case $K=1$ (so that $f=f_1$). These bounds would follow once we show that 
\begin{equation}\label{eq:K=1Sqfreeq_Inconvn_Remaining}
\sumnolimitsSt_{\substack{n: ~P_{tk+1}(n)>q}} 1 ~ ~ \ll \xonekphiqlogxOMalpha, 
\end{equation}
with $t \coloneqq 1$ if $W_k$ is not squarefull, and with $t \coloneqq 2$ in general.   
Indeed any $n$ with $P_2(n_k)>q$ automatically has $P_{k+1}(n) \ge P_{2k}(n) > q$, and any $n$ with $P_3(n_k)>q$ automatically has $P_{2k+1}(n) \ge P_{3k}(n) > q$, so \eqref{eq:MainRemainingForRestrInpSqfreeModThm}(ii), as well as its analogue with $2K+1$ replaced by $2$, would also follow once we show \eqref{eq:K=1Sqfreeq_Inconvn_Remaining}. 

To show \eqref{eq:K=1Sqfreeq_Inconvn_Remaining}, we start by estimating the contribution of the $n$'s which are divisible by the $(k+1)$-th power of a prime exceeding $q$. Any such $n$ can be written in the form $m p^c P^k$ for some positive integers $m, c$ and primes $p, P$, satisfying $P=P(n)>z$, $q<p<P$, $c \ge k+1$, $P_{Jk}(m) \le y$ and $f(n) = f(m) f(p^c) W_k(P)$. 
Recalling that 
$\#\{u \in U_q: W_k(u) \equiv b \pmod q\} \ll D^{\omega(q)}$ uniformly in $b \in \Z$, 
the argument given for the second bound in \eqref{eq:OmegaParallelOmegaStReducn} shows that the contribution of such $n$ is $\ll \frac{D^{\omega(q)}}{q^{1/k}\phi(q)} \cdot \frac{x^{1/k}}{(\log x)^{1-2\alpha_k/3}} ~ ~ \ll \xonekphiqlogxOMalpha$. 
On the other hand, for any $n$ counted in \eqref{eq:K=1Sqfreeq_Inconvn_Remaining} which is not divisible by the $(k+1)$-th power of any prime exceeding $q$, 
the condition $P_{tk+1}(n)>q$ forces $\omegaParalleln \ge t+1$ (again since $q$ is sufficiently large and the $q$-rough part of $n$ is $k$-full). Thus $n = m (P_{t+1} \cdots P_1)^k$, for some $m$ and primes $P_1, \dots, P_{t+1}$ satisfying $P_1 \coloneqq P(n)>z$, $q<P_{t+1} < \cdots < P_1$, $P_{Jk}(m) \le y$ and $f(n) = f(m) \prod_{j=1}^{t+1} W_k(P_j)$. 
The arguments leading to \eqref{eq:OmegaParallelgeKD+1_BasicBound} show that the contribution of such $n$ is 
\begin{equation}\label{eq:K=1Nok+1thPower}
\ll \frac{V_{t+1, 1}'}{\phi(q)^{t+1}} \cdot \xoneklogxOMalphaHalf \expOlogThreexSqPluslogTwoqSq.
\end{equation}
Now when $W_k$ is not squarefull (so that $t+1=2$), Proposition \ref{prop:V21V3Count_Sqfreeq}(a) shows that 
$V_{2, 1}'/\phi(q)^2 \ll \phi(q)^{-1} \exp(O(\sqrt{\log q}))$, inserting which into \eqref{eq:K=1Nok+1thPower} yields \eqref{eq:K=1Sqfreeq_Inconvn_Remaining}. 
In general (when $t+1=3$), we may invoke \eqref{eq:VNKkCountGen_Sqfreeq_LargeN} (with $K=L=1$ and $G_{1, 1} \coloneqq W_k$ $= W_{1, k}$) to see that ${V_{3, 1}'}/{\phi(q)^3}$ $\ll \phi(q)^{-1} \exp(O(\sqrt{\log q}))$,  
once again showing \eqref{eq:K=1Sqfreeq_Inconvn_Remaining}. This proves Theorem \ref{thm:RestrictedInputSqfreeMod1} for $K=1$.

We may therefore assume in the rest of the proof that $K \ge 2$. 
By replicating the arguments given for the first two bounds in \eqref{eq:OmegaParallelOmegaStReducn} (and replacing the use of Proposition \ref{prop:Vqwi_ReducnToBddModulus} by Corollary \ref{cor:VNKCountSqfrMegaGen}), we arrive at the following analogous of these two bounds: 
\begin{equation}\label{eq:OmegaParallelOmegaStReducn_Sqfreeq_Aliter}
\sumnolimitsSt_{\substack{n: ~\omegaParalleln \ge 2K+1}} 1, ~~~ \sumnolimitsSt\mathop{}_{\mkern-15mu \substack{n: ~\omegaParalleln = 2K\\\omegaStn \ge 1}} ~ 1 ~~~ \ll ~~~ \xonekphiqKlogxOMalpha,   
\end{equation}  
with the respective values of $R$. 
\subsection*{Completing the proof of Theorem \ref{thm:RestrictedInputSqfreeMod1}(a)} 
Let $R$ be any one of the two values defined in the statement (so we will only be assuming that $R \ge k(Kk+K-k)+1$ until stated otherwise). If $\omegaStn = 0$, then $k\omegaParalleln \ge R \ge k(Kk+K-k)+1$, so that 
$\omegaParalleln \ge Kk+K-k+1 \ge 2K+1$, with the last inequality being true since \textit{both} $K, k \ge 2$. As such, any $n$ with $\omegaStn = 0$ counted in \eqref{eq:MainRemainingForRestrInpSqfreeModThm}(i) is automatically counted in the first sum in \eqref{eq:OmegaParallelOmegaStReducn_Sqfreeq_Aliter}. Likewise, since \textit{both} $K, k \ge 2$, the condition $\omegaParalleln = 2K$ forces $\sum_{p>q: ~ p^{k+1} \mid n} v_p(n) \ge R - k\omegaParalleln \ge k(Kk+K-k) + 1 - 2Kk = k((K-1)(k-1)-1)+1 \ge 1$. Thus $\omegaStn \ge 1$, showing that any $n$ with $\omegaParalleln = 2K$ contributing to \eqref{eq:MainRemainingForRestrInpSqfreeModThm}(i) is counted in the second sum in \eqref{eq:OmegaParallelOmegaStReducn_Sqfreeq_Aliter}. Furthermore, by the third bound in \eqref{eq:OmegaParallelOmegaStReducn}, the contribution of all $n$ having $\omegaStn \ge Kk$ to the left hand side of \eqref{eq:MainRemainingForRestrInpSqfreeModThm}(i) is absorbed in the right hand side. It thus  
suffices to show that for any $r \in [2K-1]$ and $s \in [Kk-1]$, the contribution $\Sigmars$ of all $n$ with $\omegaParalleln = r$ and $\omegaStn = s$ to the left hand side of  \eqref{eq:MainRemainingForRestrInpSqfreeModThm}(i) is absorbed in the right hand side.  

Recall that any $n$ counted in $\Sigma_{r, s}$ is of the form $m p_1^{c_1} \cdots p_s^{c_s} P_1^k \cdots P_r^k$ for some distinct primes $p_1, \dots, p_s,$ $P_1, \dots, P_r$ and integers $m, c_1, \dots, c_s$, which satisfy the conditions (i)--(v) in the proof of Theorem \ref{thm:RestrictedInputGenMod1}(a), but with either of the current values of $R$. 
Once again, the integers $\tau_1, \dots, \tau_s$ defined by $\tau_j \coloneqq \min\{c_j, R-kr\}$ satisfy $\tau_j \in [k+1, R-kr]$, $\tau_j \le c_j$ and $\tau_1 + \cdots + \tau_s \ge R-kr$. (Here to have $R-kr \ge k+1$, it is important that $r \le 2K-1$ and $K, k \ge 2$.) As such, 
\begin{equation}\label{eq:Sigma_rs_FirstBound_SqModGen(a)}
\Sigma_{r, s} \le \sum_{\substack{\tau_1, \dots, \tau_s \in [k+1, R-kr]\\\tau_1+ \cdots + \tau_s \ge R-kr}} ~ ~ \mathcal N_{r, s}(\tau_1, \dots, \tau_s),
\end{equation}
where $\mathcal N_{r, s}(\tau_1, \dots, \tau_s)$ denotes the contribution of all $n$ counted in the left hand side of \eqref{eq:MainRemainingForRestrInpSqfreeModThm}(i) which can be written in the form $m p_1^{c_1} \cdots p_s^{c_s}$ $P_1^k \cdots P_r^k$ for some distinct primes $p_1, \dots, p_s, P_1,$ $\cdots, P_r$ and integers $m, c_1, \dots, c_s$ satisfying $c_1 \ge \tau_1, \dots, c_s \ge \tau_s$ and the conditions (i)--(v) in the proof of Theorem \ref{thm:RestrictedInputGenMod1}(a) (but with either of the current values of $R$). We will show that for each tuple $(\tau_1, \dots, \tau_s)$ occurring in \eqref{eq:Sigma_rs_FirstBound_SqModGen(a)}, we have 
\begin{equation}\label{eq:NrsIneq_SqModGen}
\mathcal N_{r, s}(\tau_1, \dots, \tau_s) ~ ~ \ll ~ ~  \xonekqKloglogxlogx \expOsqrtlogq. 
\end{equation}
Now the bound \eqref{eq:NrsGeneralBound} continues to hold, so we have 
\begin{equation}\label{eq:NrsGeneralBound_RewrittenForSqfreeq}
\mathcal N_{r, s}(\tau_1, \dots, \tau_s) \ll \frac1{q^{(\tau_1+\cdots+\tau_s)/k-s}} ~ \frac{V_{r, K}'}{\phi(q)^r} \cdot \xkonekloglogxlogx 
\end{equation}
with the current values of $r, s, \tau_1, \dots, \tau_s$ and with $V_{r, K}'$ defined in the usual manner. As such, applying \eqref{eq:VNKkCountGen_Sqfreeq_SmallN} with $L \coloneqq 1$ and $(G_{i, 1})_{i=1}^K \coloneqq (W_{i, k})_{i=1}^K$, we find that  
$$\mathcal N_{r, s}(\tau_1, \dots, \tau_s) \ll \frac{\exp\big(O(\omega(q))\big)}{q^{(\tau_1+\cdots+\tau_s)/k-s+r/2}} \cdot \xkonekloglogxlogx \ll \frac{\exp\big(O(\omega(q))\big)}{q^{\max\{s/k+r/2, ~ R/k-r/2-s\}}} \cdot \xkonekloglogxlogx.$$ 
Now $\max\{s/k+r/2, R/k-r/2-s\} > K$ whenever either \textbf{($A_1$)} $k \ge 3$, $r \ge 3$, \textit{or} \textbf{($A_2$)} $k=2$, $r \ge 4$ holds: 
indeed, if $s/k+r/2 \le K$, then $s \le k(K-r/2)$, so that (as $R \ge k(Kk+K-k)+1$) we have $R/k-r/2-s \ge K + (k-1)(r/2-1) - 1 + 1/k$.  This last quantity strictly exceeds $K$ precisely under ($A_1$) or ($A_2$), establishing \eqref{eq:NrsIneq_SqModGen} under one of these two conditions. It thus only remains to tackle the cases $r \in \{1, 2\}$, and the case $k=2, r=3$. 

The case $r=1$ is dealt with easily by inserting into \eqref{eq:NrsGeneralBound_RewrittenForSqfreeq} the trivial bound $V_{r, K}' = V_{1, K}' \ll \Dmin^{\omega(q)}$. 
It is to deal with the case $r=2$ and the case $k=2, r=3$ that we need the dichotomy in the statement of the theorem. 

\textit{When at least one of the $\Wikset$ is not squarefull $\dots$}\\
First assume that one of the polynomials $\Wikset \subZT$ is not squarefull, say $W_{1, k}$ (this is the first time in the argument for $K \ge 2$ that we are making this assumption).  If $r=2$,  
then Proposition \ref{prop:V21V3Count_Sqfreeq}(a) yields $\#\mathcal V_{2, K}(q; \wifam)/\phi(q)^2 \le \#\mathcal V_{2, 1}(q; w_1)/\phi(q)^2  \ll \phi(q)^{-1} \exp(O(\sqrt{\log q}))$, uniformly for $\wifam \in U_q^K$. Inserting this into \eqref{eq:NrsGeneralBound_RewrittenForSqfreeq}, we find that  
\begin{equation}\label{eq:N2s_WhenW1kNotSqfull}
\mathcal N_{2, s}(\tau_1, \dots, \tau_s) \ll \frac1{q^{\max\{s/k+1, R/k-1-s\}}} \cdot \xkonekloglogxlogx \expOsqrtlogq. 
\end{equation}
Since $\max\{s/k+1, R/k-1-s\}$ is always at least $K$, this establishes  \eqref{eq:NrsIneq_SqModGen} in the case when $r=2$ and one of $\Wikset$ is not squarefull.

For $k=2, r=3$, 
the multiplicative independence of $\{W_{1, k}, W_{2, k}\}$ allows us to use Proposition \ref{prop:V21V3Count_Sqfreeq}(b) to get  
${\#\mathcal V_{3, K}^{(k)}(q;\wifam)}\big/{\phi(q)^3} \ll {\exp\big(O(\omega(q))\big)}\big/{\phi(q)^2}$ uniformly for $\wifam$. By \eqref{eq:NrsGeneralBound_RewrittenForSqfreeq},  $\mathcal N_{3, s}(\tau_1, \dots, \tau_s) \ll \frac{\exp\big(O(\omega(q))\big)}{q^{\max\{s/2+2, ~ R/2-1-s\}}} \cdot \xkonekloglogxlogx$, and it is easily checked that $\max\{s/2+2, ~ R/2-1-s\}>K$. 
This shows \eqref{eq:NrsIneq_SqModGen} when one of $\Wikset$ is not squarefull. 

\textit{When all of the $\Wikset$ may be squarefull $\dots$}\\
In general (i.e., without any nonsquarefullness assumption on $\Wikset$), we can still use the second assertion of Corollary \ref{cor:VNKCountSqfrMegaGen} for $r=2$ and $3$, in place of their improved versions in Proposition \ref{prop:V21V3Count_Sqfreeq} (both of these values of $r$ are at most $2K$ as $K \ge 2$). Coming to the case $r=2$ (and $k \ge 2$), we invoke \eqref{eq:VNKkCountGen_Sqfreeq_SmallN} 
to obtain $\mathcal N_{2, s}(\tau_1, \dots, \tau_s) \ll \frac{\exp\big(O(\omega(q)\big)}{q^{\max\{s/k+1, ~ R/k-1-s\}}} \cdot \xkonekloglogxlogx,$
and we need the exponent of $q$ in the denominator to \textit{strictly} exceed $K$, in order to get a power saving of $q^{-K}$. (Compare this with \eqref{eq:N2s_WhenW1kNotSqfull} where owing to the absence of the factor $\exp\big(O(\omega(q))\big)$ we only needed the same exponent to be at least $K$.) This is where we use, for the first time (under the case $K \ge 2$) that $R = k(Kk+K-k+1)+1$. Indeed, this value of $R$ guarantees that $\max\{s/k+1, ~ R/k-1-s\}$ always exceeds $K$,  
establishing \eqref{eq:NrsIneq_SqModGen} for $r=2$. 

Finally, we turn to the case $k=2$, $r=3$. Here we only need that $R = k(Kk+K-k+1)+1 = 6K-1 \ge 6K-2$. Inserting the bound coming from \eqref{eq:VNKkCountGen_Sqfreeq_SmallN} into \eqref{eq:NrsGeneralBound_RewrittenForSqfreeq}, we get
$\mathcal N_{3, s}(\tau_1, \dots, \tau_s) \ll \frac{\expOomegaq}{q^{\max\{s/2+3/2, 3K-5/2-s\}}} \cdot \xkonekloglogxlogx$, and it is again easily seen that $\max\{s/2+3/2, 3K-5/2-s\}>K$.
This completes the proof of \eqref{eq:NrsIneq_SqModGen}, summing which over all the $O(1)$ tuples $(\tau_1, \dots, \tau_s)$ occurring in \eqref{eq:Sigma_rs_FirstBound_SqModGen(a)} establishes Theorem \ref{thm:RestrictedInputSqfreeMod1}(a).

\subsection*{Completing the proof of Theorem \ref{thm:RestrictedInputSqfreeMod1}(b)}
By arguments analogous to those given for subpart (a), it suffices to show that for any $r \in [2K-1]$ and $s \in [Kk-1]$, the contribution $\Sigmatilrs$ of all $n$ with $\omegaParalleln = r$ and $\omegakn = s$ to the left hand side of \eqref{eq:MainRemainingForRestrInpSqfreeModThm}(ii) satisfies 
\begin{equation}\label{eq:Sigmatilrs_ForSqfreeModGen(b)}
\Sigmatilrs \ll \xonekphiqKlogxOMalpha.    
\end{equation}
Any $n$ counted in $\Sigmatilrs$ has $n_k$ of the form $m' p_1^{c_1} \cdots p_s^{c_s} P_1 \cdots P_r$ for some distinct primes $p_1, \dots ,p_s,$ $P_1, \dots, P_r$ and integers $m', c_1, \dots, c_s$, which satisfy conditions (i)--(v) in the proof of Theorem \ref{thm:RestrictedInputGenMod1}(b), but with ``$KD+1-r$" replaced by ``$2K+1-r$". Hence again 
$n$ is of the form $m  p_1^{c_1 k} \cdots p_s^{c_s k} P_1^k \cdots P_r^k$, where $p_1, \dots, p_s, P_1, \dots, P_r$ are as above, $P_{Jk}(m) \le y$, and $f_i(n) = f_i(m) f_i(p_1^{c_1 k}) \cdots f_i(p_s^{c_s k}) \prod_{j=1}^r \Wik(P_j)$ for each $i \in [K]$. Defining $\tau_j \coloneqq \min\{c_j, 2K+1-r\}$ for all $j \in [s]$, we see that $\tau_j \ge 2$ (since $r \le 2K-1$) and that $\tau_1 + \cdots + \tau_s \ge 2K+1-r$. Thus 
\begin{equation}\label{eq:Sigmatil_rs_FirstBound_SqModGen(b)}
\Sigmatilrs \le \sum_{\substack{\tau_1, \dots, \tau_s \in [2, 2K+1-r]\\\tau_1+ \cdots + \tau_s \ge 2K+1-r}} ~ ~ \NtiltauOnetaus,
\end{equation} 
where (exactly as in the proof of Theorem \ref{thm:RestrictedInputGenMod1}(b)), $\NtiltauOnetaus$ denotes the number of $n$ counted in \eqref{eq:MainRemainingForRestrInpSqfreeModThm}(ii) that can be written in the form $m  p_1^{c_1 k} \cdots p_s^{c_s k} P_1^k \cdots P_r^k$ with $m,  p_1, \dots, p_s, c_1,$ $\dots, c_s, P_1, \dots, P_r$ being pairwise coprime and satisfying $P_1>z$; $q<P_r<\dots<P_1$; $p_1, \dots, p_s>q$; $P_{Jk}(m) \le y$; $c_1 \ge \tau_1, \dots, c_s \ge \tau_s$. Combining \eqref{eq:NrstilGeneralBound}, \eqref{eq:VNKkCountGen_Sqfreeq_SmallN}, and the fact that $\max\{s+r/2, ~ 2K+1-(s+r/2)\}>K$, we get 
$\NtiltauOnetaus ~~ \ll ~~ x^{1/k} \big/ \phi(q)^K (\log x)^{1-2\alpha_k/3}$, for each $\tau_1, \dots, \tau_s$ counted in \eqref{eq:Sigmatil_rs_FirstBound_SqModGen(b)}.   
This yields \eqref{eq:Sigmatilrs_ForSqfreeModGen(b)}, concluding the proof of Theorem \ref{thm:RestrictedInputSqfreeMod1}. \hfill \qedsymbol    
\subsection{Optimality in the conditions $P_{k(Kk+K-k)+1}(n)>q$ and $P_{2K+1}(n_k)>q$.}\label{subsec:Optimality_RestrictedSqfree} We will now show that the smaller value of $R$ given in Theorem \ref{thm:RestrictedInputSqfreeMod1}(a) is optimal and that the value ``$2K+1$" in (b) is nearly optimal. 
We retain the setting in subsection \cref{subsec:UnrestrictedInput_Optimality} we had used to show optimality in Theorem \ref{thm:UnrestrictedInput_SqfrPoly}(ii). To recall: fix an arbitrary $k \in \NatNos$ and $d>1$, and define $\Wik(T) \coloneqq \prod_{j=1}^{d} (T-2j) + 2(2i-1)$, so that 
$\prodik \Wik$ is separable (over $\Q$). Let $\CZeroTilde>4KD$ be any constant (depending only on $\Wikset$) exceeding the size of the (nonzero) discriminant of $\prodik \Wik$, and such that any $\CZeroTilde$-rough $k$-admissible integer lies in $\WUDkAdmSet$. Fix a prime $\ell_0>C_0$ and nonconstant polynomials $\{W_{i, v}\}_{\substack{1 \le i \le K\\1 \le v < k}} \subZT$ with all coefficients divisible by $\ell_0$. Let $\qlelogxKZ$ be any squarefree integer having $P^-(q) = \ell_0$, so that as before $q \in \WUDkAdmSet$.  
Recall also that $(2(2i-1))_{i=1}^K \in U_q^K$, that any prime $P$ satisfying $\prod_{j=1}^{d} (P-2j) \equiv 0 \pmod q$ also satisfies $f_i(P^k) \equiv 2(2i-1) \pmod q$, and that 
the congruence $\prod_{j=1}^{d} (v-2j) \equiv 0 \pmod q$ has exactly $d^{\omega(q)}$ distinct solutions $v \in U_q$.

\textit{Optimality in Theorem \ref{thm:RestrictedInputSqfreeMod1}(a).} First, we show that the condition ``$R = k(Kk+K-k)+1$" in Theorem \ref{thm:RestrictedInputSqfreeMod1}(a) cannot be weakened to ``$R = k(Kk+K-k)$". To this end, let $f_1, \dots, f_K \colon \NatNos \rightarrow \Z$ be any multiplicative functions such that $f_i(p^v) \coloneqq \Wiv(p)$ and $f_i(p^{k+1}) \coloneqq 1$ for all primes $p$, all $i \in [K]$ and $v \in [k]$. Consider $n$ of the form $(p_1 \cdots p_{k(K-1)})^{k+1} P^k \le x$ where $P, p_1, \dots, p_{k(K-1)}$ are primes satisfying the conditions $P \coloneqq P(n) > x^{1/3k}$, $q<p_{k(K-1)} < \cdots < p_1 < x^{1/4Kk^2}$, and $\prod_{1 \le j \le d} (P-2j) \equiv 0 \pmod q$. Then 
$P_{k(Kk+K-k)}(n) = p_{k(K-1)}>q$ and $f_i(n) = f_i(P^k) \prod_{j=1}^{k(K-1)} f_i(p_j^{k+1}) \equiv 2(2i-1) \pmod q$ for each $i \in [K]$. Given $p_1, \dots, p_{k(K-1)}$, the number of primes $P$ satisfying $x^{1/3k} < P \le x^{1/k}\big/(p_1 \cdots p_{k(K-1)})^{1+1/k}$ is, by the Siegel--Walfisz Theorem, $\gg {d^{\omega(q)} x^{1/k}}\big/{\phi(q) (p_1 \cdots p_{k(K-1)})^{1+1/k} \log x}$,
where we have noted that $(p_1 \cdots p_{k(K-1)})^{1+1/k} \le x^{(K-1)(k+1)/4Kk^2} \le x^{1/2k}$. Dividing by $k!$ allows us to replace the condition $p_{k(K-1)} < \cdots < p_1$ by a distinctness condition, giving us 
\begin{align}\label{eq:Kk+K-kCountereg_Calculn}
\largesum_{\substack{n \le x\\P_{k(Kk+K-k)}(n)>q\\(\forall i) f_i(n) \equiv 2(2i-1) \pmod q}} 1 ~ 
\gg ~ \frac{d^{\omega(q)} x^{1/k}}{\phi(q) \log x} \left(\mathcal T_1 - \mathcal T_2\right),
\end{align}
where $\mcTOne$ denotes the sum ignoring the distinctness condition on the $p_1, \dots , p_{k(K-1)}$, and $\mcTwo$ denotes the 
sum over all the tuples $(p_1, \dots, p_{k(K-1)})$ for which $p_i = p_j$ for some $i \ne j \in [k(K-1)]$. 
Now $\mcTOne = \prod_{1 \le j \le k(K-1)} \left(\sum_{q<p_j \le x^{1/4Kk^2}} ~ {p_j^{-(1+1/k)}}\right) \gg 1/{q^{K-1} (\log q)^{k(K-1)}}$ 
while $\mcTwo \ll \left(\sum_{p>q}  {p^{-(2+2/k)}}\right) ~ \left(\sum_{p>q} {p^{-(1+1/k)}}\right)^{k(K-1)-2} \ll 1/{q^K}$.
Consequently, the expression on the right hand side of \eqref{eq:Kk+K-kCountereg_Calculn} is $\gg {d^{\omega(q)} x^{1/k}}\big/{\phi(q)^K (\log_2 x)^{k(K-1)+1} \log x}$, which
by Proposition \ref{prop:fnqcoprimecount}, grows strictly faster than $\phi(q)^{-K} \#\{n \le x: \gcd(f(n), q)=1\}$ as soon as $d^{\omega(q)} > (\log x)^{(1+\epsilon)\alpha_k}$. We have already constructed such $q$ in subsection \cref{subsec:UnrestrictedInput_Optimality}. 
Hence, the condition $P_{k(Kk+K-k)+1}(n)>q$ in Theorem \ref{thm:RestrictedInputSqfreeMod1}(a) is optimal. 

\textit{Optimality in Theorem \ref{thm:RestrictedInputSqfreeMod1}(b).} We now address the optimality of the input restrictions in Theorem \ref{thm:RestrictedInputSqfreeMod1}(b). For $K=1$, we are assuming $P_2(n_k)>q$ when $W_k = W_{1, k}$ is not squarefull, and this is optimal for it cannot be replaced by the condition $P(n_k)>q$, as shown in \eqref{eq:PknGreaterThanqInsuff}. Turning to the condition $P_{2K+1}(n_k) > q$, we claim that it cannot be replaced by $P_{2K-1}(n_k) > q$ for any $K \ge 1$, even if $\prodik \Wik$ is assumed to be separable. 
Having already shown this above for $K=1$, we assume that $K \ge 2$.

To show our claim above, we continue with the same definitions of $\Wikset$, $\CZeroTilde$, $\ell_0$ and $q$. Let $f_1, \dots, f_K \colon \NatNos \rightarrow \Z$ be any multiplicative functions satisfying $f_i(p^v) \coloneqq \Wiv(p)$ and $f_i(p^{2k}) \coloneqq 1$ for all primes $p$, all $i \in [K]$ and $v \in [k]$. Consider any $n \le x$ of the form $(p_1 \cdots p_{K-1})^{2k} P^k$ with $P, p_1, \dots, p_{K-1}$ being primes satisfying $P \coloneqq P(n)>x^{1/3k}$, $q<p_{K-1} < \cdots < p_1 < x^{1/4Kk}$ and $\prod_{1 \le j \le d} (P-2j) \equiv 0 \pmod q$. Then $n_k = (p_1 \cdots p_{K-1})^2 P$, $P_{2K-1}(n_k) = p_{K-1} > q$ and $f_i(n) = \Wik(P) \equiv 2(2i-1) \pmod q$ for each $i \in [K]$. Given $p_1, \dots, p_{K-1}$, the number of possible $P$ is $\gg {d^{\omega(q)} x^{1/k}}\big/{\phi(q)(p_1 \cdots p_{K-1})^2 \log x}$, since $(p_1 \cdots p_{K-1})^2 \le x^{2(K-1)/4Kk}$ $\le x^{1/2k}$. 
Proceeding exactly as above, we find that the number of $n \le x$ having $P_{2K-1}(n_k)>q$ and satisfying $f_i(n) \equiv 2(2i-1) \pmod q$ for all $i$, is $\gg {d^{\omega(q)} x^{1/k}}\big/{\phi(q)^K (\log_2 x)^K \log x}$. The same $q$ as mentioned before satisfy $d^{\omega(q)} > (\log x)^{(1+\epsilon)\alpha_k}$, making this last expression grow strictly faster than $\phi(q)^{-K} \#\{n \le x: \gcd(f(n), q)=1\}$. 
The condition $P_{2K+1}(n_k)>q$ in Theorem \ref{thm:RestrictedInputSqfreeMod1}(b) is thus nearly optimal in that it cannot be replaced by $P_{2K-1}(n_k)>q$. 
\section{Restricted inputs with higher polynomial control: Proof of Theorem \ref{thm:HigherPolyControl1}}
By the same initial reductions as in the proofs of Theorems \ref{thm:RestrictedInputGenMod1} and  \ref{thm:RestrictedInputSqfreeMod1}, it suffices to show that, with the respective values of $R$ in the two subparts, we have 

\begin{align}\label{eq:MainRemaining_HigherPolyCont}\allowdisplaybreaks
\sumnolimitsSt_{\substack{n: ~P_R(n)>q}} 1 ~ ~ \ll \xonekphiqKlogxOMalpha,
\end{align}
The subsequent calculations will hold for either value of $R$ until stated explicitly. Note that (for the first time in our proofs), we will allow our implied constants to depend on $V$, and on the full set of polynomials $\Wivfullset$. 

We will first show that in either of the two subparts of the theorem, the contribution to the left hand side of \eqref{eq:MainRemaining_HigherPolyCont} from the $n$'s which are divisible by the $(V+1)$-th power of a prime exceeding $q$ can be absorbed in the right hand side. Any such $n$ can be written in the form $m p^c P^k$, where $P \coloneqq P(n)>z$, $p \in (q, P)$ is prime, $c \ge V+1$, $P_{Jk}(m) \le y$ and $P$ mod $q$ $\in \mathcal V_{1, K}^{(k)}\big(q; (a_i f_i(m p^c)^{-1})\big)$. Proceeding as in the proof of the second bound in \eqref{eq:OmegaParallelOmegaStReducn}, we see that the contribution of such $n$ is $\ll \frac{V_{1, K}'}{\phi(q) q^{(V+1)/k-1}} \cdot \frac{x^{1/k}}{(\log x)^{1-2\alpha_k/3}}$. For general $q$, an application of \eqref{eq:KonyaSub} (with $H$ being a polynomial among $\Wikset$ having least degree) shows that the expression above is $\ll {x^{1/k}}\big/{q^K(\log x)^{1-2\alpha_k/3}}$, 
since $(V+1)/k - 1 + 1/\Dmin > K$ by the hypothesis of Theorem \ref{thm:HigherPolyControl1}(a). On the other hand, if $q$ is squarefree, then from $V_{1, K}' \ll \Dmin^{\omega(q)}$ and $V \ge Kk$, 
it follows that the contribution of such $n$ is once again $\ll {x^{1/k}}\big/{q^K(\log x)^{1-2\alpha_k/3}}$.  

To prove \eqref{eq:MainRemaining_HigherPolyCont}, it thus only remains to consider the contribution of the $n$'s for which $v_p(n) \le V$ for any prime $p>q$. We may further restrict to those $n$ which have $\omegaStn \in [Kk-1]$ and $\omegaParalleln \in [KD]$ (resp. $\omegaParalleln \in [2K]$ if $q$ is squarefree). 
This is because the contribution of the $n$ having $\omegaParalleln \ge KD+1$ (resp. $\omegaParalleln \ge 2K+1$) has already been bounded in the first bound in \eqref{eq:OmegaParallelOmegaStReducn} (resp. \eqref{eq:OmegaParallelOmegaStReducn_Sqfreeq_Aliter}), while the contribution of the $n$ having $\omegaStn \ge Kk$ has already been bounded in the third bound of \eqref{eq:OmegaParallelOmegaStReducn}, and finally since any $n$ for which $\omegaStn = 0$ must anyway have $\omegaParalleln \ge KD+1$ (resp. $\omegaParalleln \ge 2K+1$) as $R \ge k(KD+1)$ (resp. $R \ge k(2K+1)$). It thus remains to show that for a given $r \in [KD]$ (resp. $r \in [2K]$) and $s \in [Kk-1]$, we have 
\begin{equation}\label{eq:Mrstil_RemainingBound}
\Mrstil \ll \xonekqKloglogxlogx,  
\end{equation}
where $\Mrstil$ denotes the contribution to the left hand side of \eqref{eq:MainRemaining_HigherPolyCont} from all the $n$ having $\omegaParalleln = r$, $\omegaStn = s$, and $k \le v_p(n) \le V$ for all $p>q$ dividing $n$. For given $r$ and $s$, 
\begin{equation}\label{eq:Mrstil_Split}
\Mrstil \le \largesum_{\substack{c_1, \dots, c_s \in [k+1, V]\\c_1 + \cdots + c_s \ge R-kr}} \MrstilcOnecs,    
\end{equation}
with $\MrstilcOnecs$ denoting the count of $n$ counted in $\Mrstil$ which can be written in the form $m p_1^{c_1} \cdots p_s^{c_s} P_1^k$ $\cdots P_r^k$, with $p_1, \dots, p_s, P_1, \dots, P_r$ being distinct primes exceeding $q$, $P(m) \le q$, $P_1 = P(n)>z$, $P_r< \cdots < P_1$, \textit{and} 
$f_i(n) = 
f_i(m) \prod_{l=1}^s W_{i, c_l}(p_l) \cdot \prod_{j=1}^r W_{i, k}(P_j)$. 
With $\Vsrkqcjw$ being the set of tuples $(u_1, \dots, u_s, v_1, \dots, v_r) \in U_q^{s+r}$ satisfying the congruences $\prod_{l=1}^s W_{i, c_l}(u_l) \cdot \prod_{j=1}^r W_{i, k}(v_j) \equiv w_i \pmod q$ for each $i \in [K]$, the conditions $\finaimodq$ amount to  
$(p_1, \dots, p_s, $ $P_1, \dots, P_r) \bmod q ~ \in ~ \mathcal V_{r+s, K}\big(q; ~ \cjfam; ~ \aifimfam)$. 

Given $m$ and $(u_1, \dots, u_s, v_1, \dots, v_r) \in \mathcal V_{r+s, K}\big(q; ~ \cjfam; ~ \aifimfam)$, we bound the number of possible $(p_1, \dots, p_s, P_1, \dots ,P_r)$ satisfying $(p_1, \dots, p_s, P_1, \dots ,P_r)$ $\equiv (u_1, \dots, u_s, v_1, \dots, v_r)$ mod $q$. First, given $(p_1, \dots, p_s)$, the number of possible $(P_1, \dots, P_r)$ is, by the arguments leading to \eqref{eq:P1...PKD+1CountInCongrClasses}, 
$\ll {x^{1/k}(\log_2 x)^{O(1)}}\big/{\phi(q)^r p_1^{c_1/k} \cdots p_s^{c_s/k}m^{1/k} \log x}$. We sum this over possible $p_1, \dots, p_s > q$, making use of the observation that for fixed $\varepsilon_1 > 0$, we have 
$\sum_{\substack{n>q\\n \equiv u \pmod q}} 1/{n^{1+\theta}}$ $ ~ \ll_{\varepsilon_1} ~ ~ 1/{q^{1+\theta}}$,
uniformly in residue classes $u$ mod $q$, and uniformly in $\theta > \varepsilon_1$.  
We find that the number of possible $(p_1, \dots, p_s, P_1, \dots ,P_r)$ 
is $\ll {x^{1/k} (\log_2 x)^{O(1)}} \big/{\phi(q)^r q^{(c_1 + \cdots + c_s)/k} m^{1/k} \log x}$. 
Finally summing the above expression over all possible $(u_1, \dots, u_s, v_1, \dots, v_r)$ and then over all $m$ via \eqref{eq:mP(m)leq}, we obtain 
\begin{equation}\label{eq:Mrstil_cOnecs_GeneralBound}
\MrstilcOnecs\\ \ll \frac1{q^{(c_1 + \cdots + c_s)/k-s}} \cdot \frac{\VrsKPrqck}{\phi(q)^{r+s}} \cdot \xkonekloglogxlogx, 
\end{equation} 
where $\VrsKPrqck \coloneqq \max \left\{\#\Vsrkqcjw : \wifam \in U_q^K\right\}$.

\subsection*{Completing the proof of Theorem \ref{thm:HigherPolyControl1}(a)}
We specialize to $R \coloneqq \max\{k(KD+1), (Kk-1)D_0+2\}$, and apply Proposition \ref{prop:VNKCountmegageneral}(b) with $\GirsetOKLmat$ being the system $(W_{i, v})_{\substack{1 \le i \le K\\k \le v \le V}}$, so that $G_{i, r} \coloneqq W_{i, k+r-1}$ and $\sumik \deg G_{i, r} = D_{k+r-1}$. We also set $N \coloneqq r+s$, and define $\FijsetOKN$ by setting (for all $i \in [K]$) $F_{i, j} \coloneqq W_{i, c_j}$ for $j \in [s]$ and $F_{i, j} \coloneqq W_{i, k}$ for $s+1 \le j \le s+r$, so that $\VNKtilqwi = \Vsrkqcjw$.   

If $r+s \ge KD_0+1$, then \eqref{eq:VNKkCount_smallN_Gen} (applied to $N \coloneqq r+s$ $\in [KD_0+1, KD+Kk-1]$ \footnote{ Here we are of course assuming that such $r$ and $s$ exist in the first place, which amounts to having $KD_0+1 \le KD+Kk-1$})  yields $\VrsKPrqck/\phi(q)^{r+s} \ll q^{-K} \expOomegaq$. Inserting this into \eqref{eq:Mrstil_cOnecs_GeneralBound} and using that $(c_1+ \cdots + c_s)/k-s \ge s/k \ge 1/k$, we obtain $\MrstilcOnecs \ll x^{1/k} (\log_2 x)^{O(1)}\big/q^K \log x$. On the other hand, if $r+s \le KD_0$, then \eqref{eq:VNKkCount_smallN_Gen} and \eqref{eq:Mrstil_cOnecs_GeneralBound} lead to 
$$\MrstilcOnecs ~ \ll ~ \frac{\big(\prod_{\ell^e \parallel q} e\big) \expOomegaq}{q^{\max\{s/k + (r+s)/D_0, ~ R/k - (1-1/D_0)(r+s)\}}} \cdot \xkonekloglogxlogx,$$
where we have recalled that $(c_1+ \cdots + c_s)/k-s \ge \max\{s/k, R/k-r-s\}$. Since $R>(Kk-1)D_0+1$, it is easy to check that the exponent of $q$ above exceeds $K$.  
This proves that $\MrstilcOnecs \ll x^{1/k} (\log_2 x)^{O(1)}\big/q^K \log x$ for any tuple $(c_1, \dots, c_s)$ counted in the sum \eqref{eq:Mrstil_Split}, and since there are $O(1)$ many such tuples, we obtain the desired bound \eqref{eq:Mrstil_RemainingBound}. 
\subsection*{Completing the proof of Theorem \ref{thm:HigherPolyControl1}(b)} This time we use Corollary \ref{cor:VNKCountSqfrMegaGen} in place of Proposition \ref{prop:VNKCountmegageneral}. If $r+s \ge 2K+1$, then \eqref{eq:VNKkCountGen_Sqfreeq_SmallN} yields $\VrsKPrqck \ll q^{-K} \expOomegaq$. Inserting this into \eqref{eq:Mrstil_cOnecs_GeneralBound} and again using $(c_1+ \dots + c_s)/k-s \ge s/k \ge 1/k$ shows that $\MrstilcOnecs \ll x^{1/k} (\log_2 x)^{O(1)}\big/q^K \log x$ in this case. On the other hand, if $r+s \le 2K$, then \eqref{eq:VNKkCountGen_Sqfreeq_SmallN} yields  
$$\MrstilcOnecs ~ \ll ~ \frac{\expOomegaq}{q^{\max\{s/k + (r+s)/2, ~ R/k - (r+s)/2\}}} \cdot \xkonekloglogxlogx,$$
and it is easy to see that the exponent of $q$ above always exceeds $K$. 
\hfill \qedsymbol

This finally establishes 
Theorems \ref{thm:UnrestrictedInput_SqfrPoly} to \ref{thm:HigherPolyControl1}. As such, we shall no longer continue with the set-up for these results. In the next section, we shall prove Theorems \ref{thm:MIHNecessity} and \ref{thm:IFHNecessity}, and thus shall only be assuming the hypotheses mentioned explicitly in their respective statements.
\section{Necessity of the multiplicative independence and invariant factor hypotheses: Proofs of Theorems \ref{thm:MIHNecessity} and \ref{thm:IFHNecessity}} 
We first give a lower bound that will be useful in both the theorems. Until we specialize to each theorem, we will not assume anything about $\Wikset \in \Z[T]$ beyond that they are nonconstant, and our estimates will be uniform in all $\qlelogxKZ$ and $\aifam \in U_q^K$. 

Let $y \coloneqq \exp(\sqrt{\log x})$ and given any fixed $R \ge 1$, we let $V_q' \coloneqq \VRKkqai = \{(v_1, \dots, v_R)$ $\in U_q^R: (\forall i \in [K]) ~ \prod_{j=1}^R \Wik(v_j) \equiv a_i \pmod q\}$. 
Consider any $N \le x$ of the form $N = (P_1 \cdots P_R)^k$, where $P_1, \dots, P_R$ are primes satisfying $y<P_R< \cdots < P_1$, and $(P_1, \dots, P_R) \bmod q ~ \in ~ V_q'$. Then $P_{Rk}(N)>y>q$ and $f_i(N) = \prod_{j=1}^R \Wik(P_j) \equiv a_i \pmod q$, 
so that estimating the count of such $N$ by the arguments leading to \eqref{eq:convnSplitForm}, we obtain for some constant $K_1>0$,
\begin{align*}\allowdisplaybreaks 
\sum_{\substack{n \le x: ~P_{Rk}(n)>q\\ \forallifinaimodq}} 1 ~ ~ \ge ~ \frac{V_q'}{\phi(q)^R} \cdot \frac1{R!} \largesum_{\substack{P_1, \dots, P_R > y\\P_1 \cdots P_R \le x^{1/k}\\ P_1, \dots, P_R \text{ distinct }}} 1 ~ - ~ x^{1/k} \exp(-K_1 (\log{x})^{1/4}). 
\end{align*}
The sum in the main term is exactly the count of squarefree $y$-rough integers $m \le x^{1/k}$ having $\Omega(m)=R$. Ignoring this squarefreeness condition with a negligible error of $O(x^{1/k}/y)$, we thus find that the main term equals $\#\{m \le x^{1/k}: P^-(m)>y, ~ \Omega(m)=R\}$, which is 
$\gg x^{1/k} (\log_2 x)^{R-1}/\log x$ by a straightforward induction on $R$ (via Chebyshev's estimates). 
So 
\begin{equation}\label{eq:ThmMIH,IFH_BasicLowerBound}
\sum_{\substack{n \le x: ~P_{Rk}(n)>q\\ \forallifinaimodq}} 1 ~ ~ \gg ~ \frac{V_q'}{\phi(q)^R} \cdot \frac{x^{1/k} (\log_2 x)^{R-1}}{\log x} ~ - ~ x^{1/k} \exp(-K_1 (\log{x})^{1/4}).
\end{equation}
\subsection*{Completing the proof of Theorem \ref{thm:MIHNecessity}} We now restrict to the $\Wikset$ and $\aifam$ considered in Theorem \ref{thm:MIHNecessity}, so $K \ge 2$, $\{\Wik\}_{1 \le i \le K-1} \subZT$ are multiplicatively independent, $W_{K, k} = \prod_{i=1}^{K-1} \Wik^{\lambda_i}$ for some tuple $(\lambda_i)_{i=1}^{K-1} \ne (0, \dots, 0)$ of nonnegative integers, and $(a_i)_{i=1}^{K} \in U_q^K$ satisfy $a_K \equiv \prod_{i=1}^{K-1} a_i^{\lambda_i} \pmod q$. The key observation is that relations assumed between the $\Wikset$ and $\aifam$ guarantee that $V_q' = \VRKkqai = \VRKMinusOnekqai$, with the set $\VRKMinusOnekqai$ defined by the congruences $\prod_{j=1}^R \Wik(v_j) \equiv a_i \pmod q$, $i \in [K-1]$. 

Define $D_1 \coloneqq \sum_{i=1}^{K-1} \degWik>1$ and  
let ``$C$" in the statement of the theorem be 
any constant $C^* \coloneqq C^*(W_{1, k},$ $\cdots, W_{K-1, k})$ exceeding $(32 D_1)^{2D_1+2}$, the sizes of the leading and constant coefficients of $\{W_{i, k}\}_{i=1}^K$, and the constant $C_1^* \coloneqq C_1(W_{1, k}, \dots, W_{K-1, k})$ coming from an application of Proposition \ref{prop:OrdDerivInfo} to the family $\{W_{i, k}\}_{i=1}^{K-1}$ of nonconstant multiplicatively independent polynomials. To show the lower bound in Theorem \ref{thm:MIHNecessity}, we may assume that $R > 4KD_1 (D_1+1)$.  
We shall carry out some of the arguments of Proposition \ref{prop:VNKCountmegageneral}; note that 
$\alpha_k(q) = \frac1{\phi(q)} \#\{u \in U_q: \prod_{i=1}^{K-1} \Wik(u) \in U_q\} \ne 0.$ 
For each prime $\ell \mid q$, we have $\gcd(\ell-1, \beta(W_{1, k},$ $\cdots, W_{K-1, k}))=1$ and $\ell>C^*$ $> C_1^*$. Thus the hypothesis $IFH(W_{1, k}, \dots, W_{K-1, k}; 1)$ holds true, and so does the corresponding analogue of the inequality \eqref{eq:CharacSum_LargePrime_Explicit}. We find that 
\begin{equation}\label{eq:ThmMIH_CharacSum_LargePrime_Explicit}
\frac1{\big(\alpha_k(\ell)\phi(\ell^e) \big)^R}\largesum_{(\chi_1, \dots, \chi_{K-1}) \ne \chiZelltuplist \bmod \ell^e} ~ \left|Z_{\ell^e; ~ \chi_1, \dots, \chi_{K-1}} (W_{1,k}, \dots, W_{K-1, k})\right|^R ~ \le ~ \frac{2 (4D_1)^R}{\ell^{R/D_1 - K}},
\end{equation}
where as usual $Z_{\ell^e; ~ \chi_1, \dots, \chi_{K-1}}(W_{1,k}, \dots, W_{K-1, k}) = \sum_{u \bmod{\ell^e}} \chi_{0, \ell}(u) \prod_{i=1}^{K-1} \chi_i(W_{i, k}(u))$. Now since $R \ge 4KD_1(D_1+1)$ and $\ell> C^* > (32 D_1)^{2D_1+2}$, we see that $\ell^{R/D_1 - K} \ge \ell^{R/(D_1+1)} \ge \ell^{R/(2D_1+2)} \cdot (C^*)^{R/(2D_1+2)} \ge \ell^2(32D_1)^R$, showing that the right hand expression in \eqref{eq:ThmMIH_CharacSum_LargePrime_Explicit} is at most $1/4\ell^2$. 
Invoking the corresponding analogue of \eqref{eq:OrthogPrimePower}, we see for each prime power $\ell^e \parallel q$ that ${\#\mathcal V_{R, K-1}^{(k)}(\ell^e; (a_i)_{i=1}^{K-1})}\big/{\phi(\ell^e)^R} ~ \ge ~ ({\alpha_k(\ell)^R}/{\phi(\ell^e)^{K-1}}) \cdot \left(1-1/{2\ell^2}\right)$. But since $\prod_{\ell \mid q} (1-1/2\ell^2) \ge 1-\frac12 \sum_{\ell \ge 2} 1/\ell^2 \ge 1/2$, we obtain ${V_q'}\big/{\phi(q)^R} = {\mathcal V_{R, K-1}^{(k)}\big(q; (a_i)_{i=1}^{K-1}\big)}\big/{\phi(q)^R} \ge {\alpha_k(q)^R}\big/{2\phi(q)^{K-1}}$, which holds true  
uniformly in $q$ having $P^-(q)> C^*$. Inserting this bound into \eqref{eq:ThmMIH,IFH_BasicLowerBound} and recalling that $\alpha_k(q) \gg 1/(\log_2 (3q))^D$, we are done. \hfill \qedsymbol 
\subsection*{Completing the proof of Theorem \ref{thm:IFHNecessity}} Again, it suffices to consider the case $R>18KD(D+1)$ to prove \eqref{eq:IFHNecess_Estim}. We start by choosing ``$C$" in the statement of the theorem to be a constant $C_2 \coloneqq C_2(\Wiklist)$ exceeding 
$(32D)^{6D+6}$, the sizes of the leading and constant coefficients of $\{W_{i, k}\}_{i=1}^K$, and the constant 
$C_1(\Wiklist)$ obtained by applying 
Proposition \ref{prop:OrdDerivInfo} 
to the family $\Wikset$ of multiplicatively independent polynomials. 
The analogue of \eqref{eq:Zell^eBoundFor e_0>1} continues to hold for each $\ell \mid q$, and the computation leading to \eqref{eq:CharacSum_LargePrime_Explicit} yields
\begin{equation}\label{eq:ThmIFH_Zell^eBoundFor e_0>1_TotalBound}
\frac1{(\alpha_k(\ell) \phi(\ell^e))^R} ~ \largesum_{\substack{\chiituplist \bmod \ell^e\\\lcmcondchilist \in \{\ell^2, \dots, \ell^e\}}} ~ |Z_{\ell^e; ~ \chi_1, \dots, \chi_K} (\Wiklist)|^R ~ \le ~ \frac{2 (4D)^R}{\ell^{R/D-K}} ~ \le ~ \frac1{4\ell^2},
\end{equation}
where in the last inequality, we have recalled that $R>4KD(D+1)$ and $\ell>C_2 \ge (32D)^{6D+6}$.

If $\chiituplist$ is a tuple of characters mod $\ell^e$ having $\lcmcondchilist = \ell$, then with $\psi_\ell$ being a generator of the character group mod $\ell$, we have $\chi_i = \psi_\ell^{A_i}$ for some unique $(A_1, \dots, A_K) \in [\ell-1]^K$ satisfying $(A_1, \dots, A_K) \not\equiv (0, \dots, 0) \pmod{\ell-1}$. 
Recall from the arguments leading to \eqref{eq:Zell^eBoundFor e_0=1} that if $\prod_{i=1}^K W_{i, k}^{A_i}$ is \textit{not} of the form $c \cdot G^{\ell-1}$ in $\FellT$, then $|Z_{\ell^e; ~ \chi_1, \dots, \chi_K} (\Wiklist)| \le D\ell^{e-1/2}$.
On the other hand, if $\prod_{i=1}^K W_{i, k}^{A_i}$ \textit{is} of that form (with $G$ monic, say), then since each $\Wik$ is monic, we must have $\prod_{i=1}^K W_{i, k}^{A_i} = G^{\ell-1}$. Since $G(v)$ is a unit mod $\ell$ iff $\prodik \Wik(v)$ is, it follows that $Z_{\ell^e; ~ \chi_1, \dots, \chi_K} (\Wiklist) = \ell^{e-1} \largesum_{v \bmod \ell} \psi_{\ell}\left((v G(v))^{\ell-1}\right) = \alpha_k(\ell) \phi(\ell^e)$.
Combining these observations with \eqref{eq:ThmIFH_Zell^eBoundFor e_0>1_TotalBound} and using that $\prodik \overline\chi_i(a_i)$ $= 1$ for any characters $\chiituplist$ mod $\ell^e$ with $\lcmcondchilist = \ell$ (as $a_i \equiv 1 \bmod \ell$), we get 
\begin{equation}\label{eq:IFHThm_Vell^eBound1}
\frac{\#\mathcal V_{R, K}^{(k)}\big(\ell^e; (a_i)_{i=1}^K\big)}{\phi(\ell^e)^R} ~ \ge ~ \frac{\alpha_k(\ell)^R}{\phi(\ell^e)^K} \left(1 + \mathcal B_\ell -\frac1{2\ell^2}\right), 
\end{equation}
where $\mathcal B_\ell$ denotes the number of tuples $(A_1, \dots, A_K) \in [\ell-1]^K \sm \{(0, \dots, 0)\}$ for which $\prodik \Wik^{A_i}$ is a perfect $(\ell-1)$-th power in $\FellT$. 

Now recalling the definition of the constant $C_1 = C_1(\Wiklist)$ from the proof of Proposition \ref{prop:OrdDerivInfo}, we know that for any $\ell > C_1$, the pairwise coprime irreducible factors of the product $\prodik \Wik$ in $\Z[T]$ continue to be separable and pairwise coprime in the ring $\FellT$. By the arguments given in the proof of Proposition \ref{prop:OrdDerivInfo}(a), $\prodik \Wik^{A_i}$ is a perfect $(\ell-1)$-th power in $\FellT$ precisely when $E_0 (A_1 \cdots A_K)^\top \equiv (0 \cdots \cdots 0)^\top \pmod{\ell-1}$, where $E_0 = E_0(\Wiklist)$ is the exponent matrix. Thus, $\mathcal B_\ell$ is exactly the number of nonzero vectors $X \in (\Z/(\ell-1)\Z)^K$ satisfying the matrix equality $E_0 X=0$ over the ring $\Z/(\ell-1)\Z$.  

Recall that $E_0$ has $\Q$-linearly independent columns and non-zero last invariant factor $\beta = \beta(\Wiklist) \in \Z$. 
By \cite[Theorem 6.4.17]{payne}, the matrix equation $E_0 X = 0$ has a nontrivial solution in the ring $\Z/(\ell-1)\Z$ precisely when 
some nonzero element of $\Z/(\ell-1)\Z$ annihilates all the $K \times K$ minors of the matrix $E_0$. 
But if $\gcd(\ell-1, \beta) \ne 1$, then the canonical image of $d \coloneqq (\ell-1)/\gcd(\ell-1, \beta)$ in $\Z/(\ell-1)\Z$ clearly does this, since $d\beta \equiv 0 \pmod{\ell-1}$ and since $\beta$ divides the gcd of the $K \times K$ minors of $E_0$ (in $\Z$). 
We thus obtain $\mathcal B_\ell \ge 1$ for each prime prime $\ell \mid q$ satisfying $\gcd(\ell-1, \beta) \ne 1$, which from \eqref{eq:IFHThm_Vell^eBound1} yields 
${V_q'}\big/{\phi(q)^R} \ge ~ 2^{\#\{\ell \mid q: ~(\ell-1, \beta) \ne 1\}} \alpha_k(q)^R\big/{2\phi(q)^K}$. Inserting this into \eqref{eq:ThmMIH,IFH_BasicLowerBound} establishes \eqref{eq:IFHNecess_Estim}. \hfill \qedsymbol 

\textbf{Remark:} If $K=1$ and $W_{1, k}$ is a constant $c$, then the $k$-admissibility of $q$ forces $\gcd(q, c)=1$, which by \eqref{eq:ThmMIH,IFH_BasicLowerBound} gives $\#\{n \le x: P_{Rk}(n)>q, f(n) \equiv c^R \pmod q\} \gg x^{1/k}(\log_2 x)^{R-1}/\log x$. 
\subsection{Explicit Examples.}\label{subsec:MIHIFH_Egs} We now construct examples where the lower bounds in Theorems \ref{thm:MIHNecessity} and \ref{thm:IFHNecessity} grow strictly faster than the expected quantity $\phi(q)^{-K} \#\{n \le x: (f(n), q)=1\}$.  
\subsection*{Failure of joint weak equidistribution upon violation of multiplicative independence hypothesis (example for Theorem \ref{thm:MIHNecessity})} By Proposition \ref{prop:fnqcoprimecount}, it is clear that the lower bound in Theorem \ref{thm:MIHNecessity} grows strictly faster once $q$ grows fast enough compared to $\log x$. 
For a concrete example, we start with any $\{W_{i, k}\}_{1 \le i \le K-1} \subset \Z[T]$ for which $\beta^* = \beta(W_{1, k}, \dots, W_{K-1, k})$ is odd (for instance, $\Wik \coloneqq H_i^{b_i}$ for some pairwise coprime irreducibles $H_1, \dots, H_{K-1} \in \Z[T]$ and odd integers $b_i>1$ satisfying $b_i \mid b_{i+1}$ for each $i<K-1$). Fix nonnegative integers $(\lambda_i)_{i=1}^{K-1} \ne (0, \dots, 0)$ and nonzero integers $\aifam$ satisfying $a_K = \prod_{i-1}^{K-1} a_i^{\lambda_i}$ (in $\Z$), and let $W_{K, k} = \prod_{i=1}^{K-1} W_{i, k}^{\lambda_i}$. Consider a constant $\widetilde C > \max\{C^*, \prodik |a_i|\}$, such that any $\widetilde C$-rough $k$-admissible integer lies in $\WUDkAdmSet$. Here $C^*$ as in the proof of Theorem \ref{thm:MIHNecessity}, so that $\widetilde C > D_1 + 1 = \sum_{i=1}^{K-1} \degWik + 1$. Let $\ell_0$ be the least prime exceeding $\widetilde C$ and satisfying $\ell_0 \equiv -1$ mod $\beta^*$. \footnote{Our arguments go through for any $c^* \in U_{\beta^*}$ for which $c^*-1 \in U_{\beta^*}$, in place of the residue $-1$ mod $\beta^*$.} Let $\{W_{i, v}\}_{\substack{1 \le i \le K\\1 \le v \le k-1}} \subZT$ be nonconstant polynomials with all coefficients divisible by $\ell_0$, and let $q \coloneqq \prod_{\substack{\ell_0 \le \ell \le Y\\\ell \equiv -1 \pmod{\beta^*}}} \ell$, with $Y$ any parameter lying in $(4|\beta^*| \log_2 x, (K_0/2) \log_2 x)$. Since $\alpha_k(\ell) \ge 1-D_1/(\ell-1) > 0$ for $\ell>\widetilde C$, we see that $\qlelogxKZ$ is $k$-admissible and hence lies in $\WUDkAdmSet$.  
As $\beta^*$ is odd and $\ell \equiv -1 \pmod{\beta^*}$ for all $\ell \mid q$, we have 
$\gcd(\ell-1, \beta^*) 
=1$ for all such $\ell$. Further, $q = \exp\Big(\sum_{\substack{\ell_0 \le \ell \le Y\\ \ell \equiv -1 \pmod{\beta^*}}} ~ \log \ell \Big) \ge \exp\left(Y/{2|\beta^*|}\right) ~ \ge ~ \log^2 x$, so 
the lower bound in Theorem \ref{thm:MIHNecessity} grows strictly faster than $\phi(q)^{-K} \#\{n \le x: (f(n), q)=1\}$.  
\subsection*{Failure of joint weak equidistribution upon violation of Invariant Factor Hypothesis (example for Theorem \ref{thm:IFHNecessity})} Define $\Wik(T) \coloneqq T-i$ for each $i \in [K-1]$ and $W_{K, k}(T) \coloneqq (T-K)^d$, for some fixed $d \in \{2, \dots, K\}$. Then $\Wikset$ are nonconstant, monic and pairwise coprime (hence multiplicatively independent); also $E_0(\Wiklist) = \diag(1, \dots, 1, d)$ so $\beta \coloneqq \beta(\Wiklist) = d$. Note that $\alpha_k(\ell) = 1-K/(\ell-1)>0$ for any prime $\ell>K+1$. 
Let $C_3 \coloneqq C_3(\Wiklist)$ be a constant exceeding the constant $C_2$ in the proof of Theorem \ref{thm:IFHNecessity}, such that any $k$-admissible $C_3$-rough integer lies in $\WUDkAdmSet$; note that $C_3>D+1 \ge K+2$. 
Let $\ell_0$ be the least prime exceeding $C_3$ and satisfying $\ell_0 \equiv 1 \pmod d$, let $\{W_{i, v}\}_{\substack{1 \le i \le K\\1 \le v < k}} \subZT$ be nonconstant polynomials all of whose coefficients are divisible by $\ell_0$, and let $q \coloneqq \prod_{\substack{\ell_0 \le \ell \le Y\\\ell \equiv 1 \pmod d}} \ell$, with $Y \le (K_0/2) \log_2 x$ a parameter to be chosen later. 

Then $\qlelogxKZ$, $P^-(q) > C_3$ and $q \in \WUDkAdmSet$. By Theorem \ref{thm:IFHNecessity} and Proposition \ref{prop:fnqcoprimecount}, it follows that the residues $a_i \equiv 1 \pmod q$ are overrepresented 
if $\#\{\ell \mid q: ~ (\ell-1, \beta) \ne 1\} \ge 4 \alpha_k \log_2 x$. But $\#\{\ell \mid q: ~ (\ell-1, \beta) \ne 1\} =$ $\sum_{\substack{\ell_0 \le \ell \le Y\\\ell \equiv 1 \pmod d}} 1 \ge Y/2\phi(d) \log Y$, whereas (since $K \ge \phi(d)$), we have $\alpha_k \le K_3/\log Y$ for some constant $K_3>0$ depending at most on $C_3$, $K$ and $d$, so we only need 
$8 K_3 \phi(d) \log_2 x < Y < (K_0/2) \log_2 x$. 

Therefore, our multiplicative independence and invariant factor hypotheses are both necessary for achieving uniformity in $\qlelogxKZ$ in our main results, and neither of them can be bypassed by restricting to inputs $n$ with sufficiently many prime factors exceeding $q$.
\section{Concluding Remarks}
It is interesting to note that despite the extensive amount of `multiplicative machinery' known in analytic number theory, there does not seem to be any estimate in the literature, a direct application of which can replace our arguments in section \ref{sec:scourfield}. For instance, Hal\'asz's Theorem only yields an upper bound on the character sums that is not precise enough, 
while a direct application of the (known forms of) the Landau-Selberg-Delange method, -- one of the most precise estimates on the mean values of multiplicative functions known in literature, -- seems to give an extremely small range of uniformity in $q$.

Theorem \ref{thm:RestrictedInputSqfreeMod1} suggests a few directions of improvement. First, 
we are still ``one step away" from optimality in the $K \ge 2$, $k=1$ case in subpart (b): we proved that ``$2K+1$" is sufficient while ``$2K-1$" is not, so the question is whether the optimal value is ``$2K$" or ``$2K+1$". If it is the former, then we will need a sharper bound on $V_{2K, K}'$ than what comes from our methods in section \ref{sec:ResInpSqfreePf}. 
One can also ask whether it is possible to weaken the nonsquarefullness condition in subpart (a). 
Theorem 
\ref{thm:HigherPolyControl1} 
also suggests other avenues for improvement, for instance, by optimizing the values of $R$ and $V$. 
We hope to return to these questions in future papers.
\section*{Acknowledgements}
This work was done in partial fulfillment of my PhD at the University of Georgia. As such, I would like to thank my advisor, Prof. Paul Pollack, for the past joint research and fruitful discussions that have led me to think about this question, as well as for his continued support and encouragement. I would also like to thank the Department of Mathematics at UGA for their support and hospitality. 

\textbf{Data Availability} The manuscript has no associated data. 

\section*{Declarations}
\textbf{Conflict of Interest} On behalf of all authors, the corresponding author states that there is no conflict of interest.

\providecommand{\bysame}{\leavevmode\hbox to3em{\hrulefill}\thinspace}
\providecommand{\MR}{\relax\ifhmode\unskip\space\fi MR }
\providecommand{\MRhref}[2]{
  \href{http://www.ams.org/mathscinet-getitem?mr=#1}{#2}
}
\providecommand{\href}[2]{#2}

\end{document}